\documentclass[11pt,a4paper]{article}

\usepackage{amsmath,amsthm}
\PassOptionsToPackage{backend=bibtex}{biblatex}
\usepackage{cgeometry}
\allowdisplaybreaks
\hypersetup{
  colorlinks=true,
  linkcolor=blue!55!black,
  citecolor=green!40!black,
  urlcolor=blue!60!black,
}

\setlist[enumerate]{label=(\arabic*),leftmargin=2.3em}

\newcommand{\an}{\mathrm{an}}  
\newcommand{\divv}{\mathrm{div}}    
\newcommand{\vtriv}{v_{\mathrm{triv}}}
\newcommand{\PSHhom}{\PSH_{\mathrm{hom}}}
\newcommand{\EP}{\operatorname{EP}}
\newcommand{\Norm}{\operatorname{N}}
\newcommand{\Av}{\operatorname{Av}}

\title{The envelope property for unibranch projective varieties,\\
\large with an appendix on the singular Yau--Tian--Donaldson conjecture}
\author{Yanbo Fang, Jiyuan Han and Mingchen Xia}
\date{}

\begin{document}

\maketitle

\begin{abstract}
We prove the envelope property for every real numerical class on an integral unibranch projective variety over an algebraically closed trivially valued field based on the recent work of L.~Wu. Together with Boucksom--Jonsson's converse, this characterizes unibranchness by the plurisubharmonic envelope property. The Appendix applies the result to the uniform Yau--Tian--Donaldson conjecture for projective varieties with klt singularities.
\end{abstract}
\section{Introduction}

\subsection{The envelope conjecture}

Let $k$ be algebraically closed with the trivial absolute value, and let $X$ be an integral projective $k$-variety. Boucksom--Jonsson \cite{BJ22} associate with each $\theta\in N^1(X)_{\mathbf R}$ a space $\PSH(X,\theta)$ of plurisubharmonic (psh) functions on $X^{\an}$, the Berkovich analytification of $X$.

Various foundational claims in non-Archimedean pluripotential theory relies crucially on the variational solution of Monge--Amp\`ere equations depend on the \emph{envelope property}: for a fixed class $\theta$ one asks that for every nonempty family $(\varphi_i)_{i\in I}\subseteq\PSH(X,\theta)$ bounded from above, the upper semicontinuous (usc) regularization of the pointwise supremum,
\[
  \sups_i\varphi_i,
\]
again lies in $\PSH(X,\theta)$.  We abbreviate this statement as $\EP_X(\theta)$.  It is the analogue of the classical fact that the usc regularization of the upper envelope of a locally uniformly bounded above family of psh functions on a domain in $\mathbf C^n$ is psh.

Boucksom and Jonsson proved the envelope property for every numerical class on a smooth projective variety when $\operatorname{char}k=0$ or $\dim X\leq2$; see \cite[Theorem~5.20]{BJ22} for nef classes and \cite[Theorem~A]{BJaddendum} in general. Conversely, if $\EP_X(\omega)$ holds for one ample class $\omega$, then $X$ is unibranch \cite[Theorem~5.12]{BJ22}. Boucksom and Jonsson conjectured that unibranchness is sufficient for the envelope property in every class \cite[Conjecture~5.14]{BJ22}.

The purpose of this paper is to prove the conjecture in full generality.

\begin{maintheorem}\label{thm:intro-main}
Let $X$ be an integral unibranch projective variety over an algebraically closed trivially valued field $k$, and let $\theta\in N^1(X)_{\mathbf{R}}$.  Then every nonempty bounded-above family $(\varphi_i)\subseteq\PSH(X,\theta)$ satisfies
\[
    \sups_i\varphi_i \in\PSH(X,\theta).
\]
\end{maintheorem}

The theorem applies in every dimension and characteristic. Together with \cite[Theorem~5.12]{BJ22}, it yields the following characterization.

\begin{samepage}
\begin{corollary}\label{cor:intro-characterization}
For an integral projective variety $X$ over $k$, the following are equivalent:
\begin{enumerate}
\item $X$ is unibranch;
\item $\EP_X(\theta)$ holds for every $\theta\in N^1(X)_{\mathbf{R}}$;
\item $\EP_X(\omega)$ holds for some ample class $\omega$.
\end{enumerate}
\end{corollary}
\end{samepage}

Indeed $(1)\Rightarrow(2)$ is \cref{thm:intro-main}, $(2)\Rightarrow(3)$ follows by choosing an ample class, and $(3)\Rightarrow(1)$ is \cite[Theorem~5.12]{BJ22}.

\subsection{Proof strategy}

We first prove the envelope property for nef classes on smooth varieties (\cref{thm:smooth-EP}) following the strategy of \cite{BJaddendum}. In positive characteristic, asymptotic test ideals on the smooth product $X\times\mathbf A^1$ replace multiplier ideals on a resolution of a test configuration.

For a normal variety over a countable field $k$, we use the Stein factorization $Y\to T\to X$ of a smooth alteration to descend the envelope property from $Y$ to $X$. The maximum principle of Gubler--Rabinoff \cite[Theorem~6.1]{GublerRabinoff25} makes the envelope on $Y$ constant on the fibres of the birational morphism $Y\to T$. Local regularization then shows that the descended function is psh (\cref{cor:normal-weak-global}), and weighted field norms descend the envelope property through the finite map $T\to X$ (\cref{thm:finite-normal-descent}). Flat base change of base ideals reduces the problem over an arbitrary field to a countable subfield (\cref{thm:normal-ample-EP}).

The local regularization argument follows Wu \cite{Wu26}. On the open set where the modification is an isomorphism, regularizability follows directly from the global psh function upstairs. Polynomial hulls and norms along finite projections then extend it across the complement. This uses only the envelope property of projective space, already supplied by the smooth case.

Finally, we descend through the normalization: Frobenius powers suffice in positive characteristic, while Roberts--Singh additive powers give a decreasing approximation after a fixed twist in characteristic zero (\cref{thm:normalization-reflection}).

\subsection{Related work}

The smooth characteristic-zero case and the positive-characteristic surface case are due to Boucksom--Jonsson \cite{BJ22,BJaddendum}. Wu \cite{Wu26} developed a local theory of non-Archimedean psh functions and applied it to the envelope property of normal varieties in residue characteristic zero. \Cref{sec:normal-local-regularization} follows his approach; we indicate there where his arguments have to be corrected or extended to positive characteristic. Our positive-characteristic argument in the smooth case uses asymptotic test ideals \cite{HaraYoshida,MustataTest,BMS08}, and the normalization argument uses subintegrality \cite{Swan80,RobertsSingh}. Alterations \cite{deJong} and the maximum principle for classically psh functions \cite{GublerRabinoff25} enter the passage from smooth to normal varieties.

\subsection{Conventions}\label{subsec:conventions}
Throughout, $k$ is an algebraically closed field endowed with the trivial absolute value.  A \emph{variety} is an integral separated scheme of finite type over $k$; unless stated otherwise, all varieties are projective. The notation $N^1(X)_{\mathbf{R}}$ denotes the space of numerical classes of real Cartier divisors, and a \emph{rational line bundle} is an element of $\Pic(X)_{\mathbf{Q}}$.  Numerical classes are real unless the contrary is stated. Countability of $k$ is assumed only in \cref{sec:normal-local-regularization} and in Step~1 of the proof of \cref{thm:normal-ample-EP}.

\subsection*{Acknowledgments}
The initial proof of the envelope property was assisted by the \href{https://github.com/frenzymath/Danus}{Danus agent} (improved by Felix Ye), and the gpt-5.6-sol model. The authors then simplified the proof and improved the presentation, and are fully responsible for all assertions in this paper. The Appendix, written by the authors, slightly generalizes Boucksom's argument in \cite[Appendix~A]{T26App}. A sorry-free Lean~4 formalization of the main body of the paper and a list of mistakes in the literature we found during the formalization process are available at \url{https://github.com/MingchenXia/envelope-property-lean}.

We would like to thank Sébastien Boucksom and Mattias Jonsson for relevant discussions.

The second-named author is supported by grant 2024YFA1014800 of the National Key R\&D Program of China, NSFC-12671075 and XHD23A0101. The third-named author is supported by the National Key R\&D Program of China under Grant No.~2025YFA1018200.

\section{Preliminaries}\label{sec:preliminaries}

Let $X$ be an integral projective variety over $k$. We fix notation for its analytification and psh functions, and record the behavior of Fubini--Study classes under finite morphisms.

\subsection{The analytification and its semivaluations}
\label{subsec:prelim-berkovich}

We regard $k$ as a trivially valued non-Archimedean field and write $X^{\an}$ for the Berkovich analytification of $X$ \cite{Berkovich90}; see \cite[\S1]{BJ22} for a detailed discussion of the trivially valued case.  A point $v\in X^{\an}$ is a \emph{semivaluation}: it is given by a point $\xi\in X$, the \emph{support} of $v$, together with a real valuation of the residue field $\kappa(\xi)$, trivial on $k$.  We use additive notation, so $v(f)\in\mathbf{R}\cup\{+\infty\}$ for a local function $f$, with $v(f)=+\infty$ exactly when $f$ vanishes along the support.  Since $X$ is proper, every $v\in X^{\an}$ admits a unique \emph{center} $c(v)\in X$, the point such that $v\geq0$ on $\mathcal{O}_{X,c(v)}$ and $v>0$ on the maximal ideal.

The analytification is a compact Hausdorff topological space, functorial in $X$; for a morphism $f\colon X\to X'$ we write $f^{\an}$ for the induced map, and we normally abbreviate $f^{\an}(v)$ to $f(v)$.

Two classes of points play a special role.  The \emph{trivial valuation} $\vtriv$ is the zero valuation on $k(X)^\times$.  The set $X^{\divv}\subseteq X^{\an}$ of \emph{divisorial points} consists of $\vtriv$ together with all positive rational multiples $c\ord_E$, where $E$ is a prime divisor on some normal projective birational model of $X$ and $c\in\mathbf{Q}_{>0}$.  Points of $X^{\divv}$ have support the generic point of $X$; the support of a general point of $X^{\an}$ may be any point of $X$.

The multiplicative group $\mathbf{R}_{>0}$ acts continuously on $X^{\an}$ by scaling: $t\cdot v$ is the semivaluation $tv$ with the same support.  The induced \emph{scaling action} on functions is
\begin{equation}\label{eq:scaling-action}
  (t\mathbin\cdot u)(v)\coloneqq t\,u(t^{-1}v),\qquad t>0,
\end{equation}
and a function $u$ is \emph{homogeneous} if $t\mathbin\cdot u=u$ for all $t>0$.

Finally, let $S=X\times\mathbf{A}^1$ be the trivial test configuration, with coordinate $t$ on $\mathbf{A}^1$.  Every $v\in X^{\an}$ has a \emph{Gauss extension} $\sigma(v)$: the unique $\mathbf{G}_m$-invariant semivaluation on $S$ whose restriction to functions from $X$ is $v$, normalized by $\sigma(v)(t)=1$; explicitly, $\sigma(v)\bigl(\sum_\lambda s_\lambda t^\lambda\bigr) =\min_\lambda\{v(s_\lambda)+\lambda\}$ for local sections $s_\lambda$ of $\mathcal{O}_X$.  By \cite[\S1.3]{BJ22}, $\sigma$ is a bijection from $X^{\an}$ (resp.\ $X^{\divv}$) onto the set of $\mathbf{G}_m$-invariant semivaluations (resp.\ divisorial valuations) $w$ on $S$ with $w(t)=1$; in particular $\sigma(X^{\divv})\subseteq S^{\divv}$.  For a nonzero (equivariant) coherent ideal $I\subseteq\mathcal{O}_S$ we set
\[
  \phi_I(v)\coloneqq-\sigma(v)(I),
  \qquad v\in X^{\an},
\]
where $\sigma(v)(I)$ is the minimum of $\sigma(v)$ over local generators of $I$ at the center.  We explicitly allow $\phi_I(v)=-\infty$, which happens when the support of $\sigma(v)$ lies inside the zero locus of $I$. The support of $\sigma(v)$ is the generic point of $\operatorname{supp}(v)\times\mathbf A^1$, which does not lie on the central fibre $X\times\{0\}$. Hence, if $\mu\colon\mathcal X\to S$ is a projective birational morphism that is an isomorphism outside the central fibre, $\sigma(v)$ lifts uniquely to a semivaluation on $\mathcal X$, and we define $\phi_I$ in the same way for ideals $I\subseteq\mathcal O_{\mathcal X}$. For a Cartier divisor $D$ on $\mathcal X$ supported on the central fibre we put $\phi_D(v)\coloneqq\sigma(v)(D)$, the value of $\sigma(v)$ on a local equation of $D$; thus $\phi_{\mathcal O_{\mathcal X}(-D)}=-\phi_D$ for $D\geq0$.

\subsection{Plurisubharmonic functions and the envelope property}
\label{subsec:prelim-psh}

For a rational line bundle $L\in\Pic(X)_\mathbf{Q}$, a \emph{rational, generically finite Fubini--Study function} for $L$, or simply a \emph{Fubini--Study function}, is a function of the form
\begin{equation}\label{eq:FS-function}
  \varphi=\frac1m\max_j\{\log|s_j|+\lambda_j\},
\end{equation}
where $m\in\mathbf{Z}_{>0}$ is such that $mL$ is an honest line bundle, the $s_j\in H^0(X,mL)$ form a finite family of sections, not all zero, and $\lambda_j\in\mathbf{Q}$.  Here $\log|s|(v)=-v(s)$ is computed in a local trivialization at the center of $v$; the choice of trivialization does not matter.  We denote by $\mathcal{H}^{\mathrm{gf}}_{\mathbf{Q}}(X,L)$, or $\mathcal{H}^{\mathrm{gf}}_{\mathbf{Q}}(L)$ when $X$ is understood, the set of such functions; this is the space denoted $\mathcal H^{\mathrm{gf}}_{\mathbf{Q}}(L)$ in \cite{BJ22,BJaddendum}.

When applying \cite[Definition~4.1]{BJ22}, we also allow real weights in \eqref{eq:FS-function}; the notation $\mathcal H^{\mathrm{gf}}_{\mathbf Q}$ remains reserved for rational weights. This gives the same psh class. Indeed, for each member $f_j$ of a decreasing real-weight net and each $\epsilon>0$, round the weights upwards to obtain a rational-weight function $r$ with $f_j+\epsilon\le r\le f_j+2\epsilon$. Order the choices by increasing $j$, decreasing $r$ and decreasing $\epsilon$. They form a directed set: after taking a common later original index, choose twice the new $\epsilon$ smaller than both old ones. The inequalities are preserved even at common zeros, since rounding does not change a pole set. The resulting rational-weight net has the same limit and class limit.

Since the sections may have common zeros, a Fubini--Study function may take the value $-\infty$; it is finite at $\vtriv$, and more generally at every point whose support is not contained in the common zero locus.  We write $\mathcal{H}_{\mathrm{hom}}(X,L)\subseteq\mathcal{H}^{\mathrm{gf}}_{\mathbf{Q}}(X,L)$ for the subset of \emph{homogeneous} Fubini--Study functions, namely those of the form \eqref{eq:FS-function} with all $\lambda_j=0$.

For $\theta\in N^1(X)_\mathbf{R}$, Boucksom--Jonsson define the class $\PSH(X,\theta)$ of $\theta$-psh functions $u\colon X^{\an}\to\mathbf{R}\cup\{-\infty\}$, $u\not\equiv-\infty$; when $\theta=c_1(L)$ we also write $\PSH(X,L)$.  We further write
\[
  \PSHhom(X,\theta)\coloneqq
  \{u\in\PSH(X,\theta): u\text{ is homogeneous}\}.
\]
We use the following properties of $\PSH(X,\theta)$.  Recall that $\sups$ denotes the upper semicontinuous (usc) regularization of a supremum.

We summarize a few results we need in the sequel from \cite{BJ22}.
\begin{proposition}\label{prop:BJ-package}
Let $X$ be an integral projective variety and $\theta\in N^1(X)_{\mathbf{R}}$.
\begin{enumerate}
\item \textup{(PSH2)} If $\theta_i\to\theta$ and $u_i\in\PSH(X,\theta_i)$ decreases pointwise, its limit belongs to $\PSH(X,\theta)$ unless it is identically $-\infty$, provided that $\theta$ is nef.
\item  Every $\theta$-psh function is usc and real valued on $X^{\divv}$.  If $u\in\PSH(X,\theta)$ and $h\colon X^{\an}\to\mathbf{R}\cup\{-\infty\}$ is usc, then $u\leq h$ on $X^{\an}$ if and only if $u\leq h$ on $X^{\divv}$.  In particular, every $\theta$-psh function is determined by its restriction to $X^{\divv}$.
\item  Every $\varphi\in\mathcal{H}^{\mathrm{gf}}_{\mathbf{Q}}(X,L)$ belongs to $\PSH(X,c_1(L))$.  Conversely, every $u\in\PSH(X,\theta)$ is the pointwise limit of a decreasing sequence $\varphi_i\in\mathcal{H}^{\mathrm{gf}}_{\mathbf{Q}}(X,L_i)$ with rational classes $c_1(L_i)\to\theta$; if $u$ is homogeneous, the $\varphi_i$ can be taken in $\mathcal{H}_{\mathrm{hom}}(X,L_i)$.
\item  If $f\colon X'\to X$ is a surjective morphism of integral projective varieties, then $f^{\an}$ is surjective, and $f^*\coloneqq(\cdot)\circ f^{\an}$ maps $\PSH(X,\theta)$ into $\PSH(X',f^*\theta)$.
\item  Finite maxima of $\theta$-psh functions are $\theta$-psh, and so is $u+c$ for $u\in\PSH(X,\theta)$ and $c\in\mathbf{R}$.  Ordinary scalar multiplication satisfies $a u\in\PSH(X,a\theta)$ for $a>0$, and $\PSH(X,\theta)+\PSH(X,\eta)\subseteq\PSH(X,\theta+\eta)$. In particular, if $\eta-\theta$ is nef, then $\PSH(X,\theta)\subseteq\PSH(X,\eta)$.
\item  Every $u\in\PSH(X,\theta)$ attains its supremum at the trivial valuation: $\sup_{X^{\an}}u=u(\vtriv)\in\mathbf{R}$.  Every homogeneous $\theta$-psh function vanishes at $\vtriv$, hence is nonpositive.
\end{enumerate}
\end{proposition}

\begin{lemma}\label{lem:homogeneous-decreasing-stability}
Let $X$ be an integral projective variety. If $\theta_i\to\theta$ and a net $u_i\in\PSHhom(X,\theta_i)$ decreases pointwise to $u$, then $u\in\PSHhom(X,\theta)$.
\end{lemma}

\begin{proof}
\emph{Step 1: uniform comparison.} Let $p_j$ be a decreasing net of homogeneous Fubini--Study functions with pointwise limit $p$, and let $q$ be a homogeneous Fubini--Study function with $p\le q$. We claim that, for every rational $0<a<1$,
\[
 p_j\le aq\quad\text{on }X^{\an}\text{ for all sufficiently large }j.
\]
On the compact set $K\coloneqq\{q=-1\}$, we have $p_j\downarrow p\le-1<-a$. Thus the increasing open sets $\{p_j<-a\}$ cover $K$, and compactness gives one index after which $p_j<-a$ throughout $K$. Homogeneity extends the comparison to $\{-\infty<q<0\}$; on $\{q=0\}$ it follows from $p_j\le0$.

It remains to consider $\{q=-\infty\}=Z^{\an}$, where $Z$ is the common zero locus of the sections defining $q$. At the trivial semivaluation supported at the generic point of an irreducible component $W$ of $Z$, each $p_j$ has value $0$ or $-\infty$. Since its limit there is $-\infty$, all the sections defining $p_j$ eventually vanish on $W$. Hence $p_j=-\infty$ on $W^{\an}$ for all sufficiently large $j$. There are finitely many components, so a common index works on all of $Z^{\an}$. This proves the claim.

\emph{Step 2: a decreasing Fubini--Study approximation.} By \cref{prop:BJ-package}(3), for each $i$ choose a decreasing sequence $p_{i,j}\downarrow u_i$ of homogeneous Fubini--Study functions, with numerical classes $\eta_{i,j}\to\theta_i$ as $j\to\infty$. Fix a norm on $N^1(X)_{\mathbf R}$. For each triple $\alpha=(i,j,b)$ with rational $0<b<1$, put
\[
 F_\alpha\coloneqq bp_{i,j},\qquad
 e_\alpha\coloneqq\|\eta_{i,j}-\theta_i\|+1-b.
\]
Then $F_\alpha$ is homogeneous Fubini--Study in the class $b\eta_{i,j}$, and $F_\alpha\ge u_i\ge u$, since $p_{i,j}\le0$. Declare $\alpha'\succeq\alpha$ when $i'\ge i$, $F_{\alpha'}\le F_\alpha$ on $X^{\an}$, and $e_{\alpha'}\le e_\alpha$.

To prove directedness, take two triples $\alpha_k=(i_k,j_k,b_k)$, $k=1,2$, and choose $0<\delta<\min\{e_{\alpha_1},e_{\alpha_2}\}$. Choose $i\ge i_1,i_2$, as late as desired, and rational $b$ with
\[
 \max\{b_1,b_2\}<b<1,\qquad 1-b<\delta/2.
\]
Since $u_i\le u_{i_k}\le p_{i_k,j_k}$, Step~1, applied with $a=b_k/b$, gives
\[
 p_{i,j}\le\frac{b_k}{b}\,p_{i_k,j_k}\qquad(k=1,2)
\]
for all sufficiently large $j$. Increase $j$ further so that $\|\eta_{i,j}-\theta_i\|<\delta/2$. The resulting triple is a common successor, with $F_\alpha\le F_{\alpha_1},F_{\alpha_2}$ and $e_\alpha<\delta$.

Thus the $F_\alpha$ form a decreasing net. Since $i$ can be chosen arbitrarily late and $e_\alpha$ arbitrarily small, their classes $b\eta_{i,j}$ tend to $\theta$. Their pointwise infimum is $u$: given $v\in X^{\an}$ and a real number $A>u(v)$, first choose $i,j$ with $p_{i,j}(v)<A$, then choose $b<1$ sufficiently close to $1$ that $F_\alpha(v)<A$. This argument also applies when $u(v)=-\infty$. Finally, $u(\vtriv)=\lim_i u_i(\vtriv)=0$, so \cite[Definition~4.1]{BJ22} gives $u\in\PSH(X,\theta)$. Homogeneity passes to the pointwise limit.
\end{proof}

Recall from the introduction that $\EP_X(\theta)$ denotes the envelope property of $\theta$. It has an increasing-net formulation which is convenient for descent arguments.

\begin{proposition}\label{prop:EP-criterion}
The envelope property $\EP_X(\theta)$ is equivalent to the following assertion: for every bounded-above increasing net $(\varphi_i)$ in $\PSH(X,\theta)$ there is a \textup(necessarily unique\textup) $\varphi\in\PSH(X,\theta)$ with
\[
  \varphi=\sup_i\varphi_i
  \quad\text{pointwise on }X^{\divv}.
\]
In this case $\varphi=\sups_i\varphi_i$; in particular $\varphi$ is homogeneous if all $\varphi_i$ are.
\end{proposition}

If the envelope property holds, \cite[Theorem~5.6]{BJ22} gives $(\sup_i\varphi_i)^*=\sup_i\varphi_i$ on $X^{\divv}$, which proves the asserted implication. Conversely, it suffices to consider the increasing net of finite maxima of an arbitrary bounded-above family. By \cref{prop:BJ-package}(2), $\varphi\ge\varphi_i$ everywhere, so $\varphi\ge\sups_i\varphi_i$, and $\varphi\le\sups_i\varphi_i$ since this holds on $X^{\divv}$. The scaling action is by homeomorphisms of $X^{\an}$, so it commutes with pointwise suprema and usc regularization; this gives the last assertion.

\subsection{Envelopes in ample classes}\label{subsec:prelim-envelopes}

For an ample $L\in\Pic(X)_{\mathbf Q}$ and a bounded function $\phi\colon X^{\an}\to\mathbf R$, the \emph{psh envelope} is the pointwise supremum
\[
  \mathrm{P}_L(\phi)\coloneqq\sup\{\tau\in\PSH(X,L):\tau\leq\phi\},
\]
without usc regularization. It is increasing in $\phi$, and $\mathrm P_L(\phi+c)=\mathrm P_L(\phi)+c$ for $c\in\mathbf R$.

Rational piecewise linear (PL) functions on $X^{\an}$ are described by test configurations \cite[\S2]{BJ22}. For the model description in this paragraph, assume that $X$ is normal. Let $L$ be a line bundle and $(\mathcal X,\mathcal L)$ a test configuration for $(X,L)$ with $\mathcal L$ a line bundle. Passing to a higher test configuration, we may assume that $\mathcal X$ is normal and dominates $S=X\times\mathbf A^1$ through a projective birational morphism $\mu\colon\mathcal X\to S$ that is an isomorphism outside the central fibre. Then $\mathcal L=\mu^*L+D$ for a Cartier divisor $D$ supported on the central fibre, and the associated PL function is $\phi_{\mathcal L}=\phi_D$. By \cite[Theorem~2.31]{BJ22}, after multiplying the class and the function by a common positive integer, every rational PL function admits such a description. Replacing $\mathcal L$ by $\mathcal L+c\,\mathcal X_0$, with $\mathcal X_0=\operatorname{div}(t)$ the central fibre and $c\in\mathbf Z$, adds $c$ to $\phi_{\mathcal L}$.

\begin{proposition}\label{prop:ample-envelopes}
Let $L\in\Pic(X)_{\mathbf Q}$ be ample.
\begin{enumerate}
\item $\EP_X(c_1(L))$ holds if and only if $\mathrm P_L(\phi)$ is $L$-psh for every rational PL function $\phi$ on $X^{\an}$.
\item Let $L$ be a line bundle and $(\mathcal X,\mathcal L)$ a test configuration for $(X,L)$ with $\mathcal L$ a line bundle, and let $\mathfrak a_m\subseteq\mathcal O_{\mathcal X}$ be the base ideal of $m\mathcal L$, generated by $H^0(\mathcal X,m\mathcal L)$. Then
\[
\mathrm P_L(\phi_{\mathcal L})=\phi_{\mathcal L}+\sup_{m\ge1}\frac1m\phi_{\mathfrak a_m}=\phi_{\mathcal L}+\lim_{m\to\infty}\frac1m\phi_{\mathfrak a_m}\qquad\text{pointwise on }X^{\an}.
\]
\end{enumerate}
\end{proposition}

Part (2) is \cite[Lemma~5.22]{BJ22} (see also \cite[Lemmas~1.5 and~1.6]{BJaddendum}). For (1), we give the PL-obstacle reduction in \cite[Lemma~1.3]{BJaddendum} explicitly. If the envelope property holds, then the usc regularization of $\mathrm P_L(g)$ is $L$-psh and remains bounded above by the continuous obstacle $g$; hence it equals $\mathrm P_L(g)$. Conversely, assume that $\mathrm P_L(g)$ is $L$-psh for every rational PL function $g$. Let $f$ be the supremum of a nonempty bounded-above family of $L$-psh functions, and let $\mathcal G$ be the rational PL majorants of $f$. This family is nonempty and stable under finite minima. Since rational PL functions are uniformly dense in $C^0(X^{\an})$ \cite[Theorem~2.2]{BJ22}, and continuous majorants recover the usc regularization on the compact Hausdorff space $X^{\an}$, we have $\inf_{g\in\mathcal G}g=f^*$. Ordered by decreasing obstacles, $\mathrm P_L(g)$ is a decreasing net of $L$-psh functions, bounded below by every member of the original family. PSH2 in the ample class $c_1(L)$ gives a psh limit $u$. Since $f\le u\le f^*$ and $u$ is usc, $u=f^*$. This proves (1).

\subsection{Fubini--Study classes under finite morphisms}
\label{subsec:prelim-finite-classes}

Let $h\colon T\to X$ be a finite surjective morphism of integral projective varieties. The following lemma represents Fubini--Study functions on $T$ in line bundles pulled back from $X$.

\begin{lemma}\label{lem:finite-class-promotion}
Let $h\colon T\to X$ be a finite surjective morphism of integral projective varieties and let $\theta\in N^1(X)_{\mathbf R}$.  Suppose
\[
 \varphi_i\in\mathcal{H}^{\mathrm{gf}}_{\mathbf{Q}}(T,M_i),\qquad
 c_1(M_i)\longrightarrow h^*\theta
\]
is a sequence.  After discarding finitely many initial terms, there are $L_i\in\Pic(X)_{\mathbf Q}$ with $c_1(L_i)\to\theta$ such that the same functions belong to $\mathcal{H}^{\mathrm{gf}}_{\mathbf{Q}}(T,h^*L_i)$, presented by sections of $h^*(mL_i)$ with $mL_i$ an honest line bundle on $X$.
\end{lemma}

\begin{proof}
Fix an ample line bundle $A$ on $X$.  Since $h^*A$ is ample, choose positive rational $\epsilon_i\to0$ sufficiently slowly that $c_1(M_i)-h^*\theta=o(\epsilon_i)$.  Choose rational line bundles $L_i$ with $c_1(L_i)=\theta+\epsilon_i c_1(A)+o(\epsilon_i)$.  The space $N^1(T)_{\mathbf R}$ is finite dimensional, and the ample cone is open by Kleiman's criterion. Putting $P_i\coloneqq h^*L_i-M_i$, we have $\epsilon_i^{-1}c_1(P_i)\to h^*c_1(A)$, so $P_i$ is ample for all large $i$.

Fix $i$ and take $r$ sufficiently divisible that $rP_i$ is an honest globally generated line bundle, with generators $(t_a)$, and that $rL_i$ is an honest line bundle on $X$.  At every semivaluation, all their local coefficients are regular at the center and some coefficient is a unit.  Thus $r^{-1}\max_a\log|t_a|=0$ on $T^{\an}$.  Writing $\varphi_i=m^{-1}\max_j\{\log|s_j|+\lambda_j\}$ gives
\[
 \varphi_i=\frac1{rm}\max_{j,a}
 \bigl\{\log|s_j^{\otimes r}\otimes t_a^{\otimes m}|+r\lambda_j\bigr\},
 \qquad
 s_j^{\otimes r}\otimes t_a^{\otimes m}
 \in H^0(T,rmh^*L_i).
\]
This equality also holds at common zeros of the $s_j$, since the $t_a$ have no common zero, and it proves the claim.
\end{proof}

\section{The smooth case}\label{sec:smooth}

In this section $X$ is a smooth projective variety over $k$ and $n=\dim X$.

\begin{theorem}\label{thm:smooth-EP}
Let $X$ be a smooth projective variety over $k$. Then every nef class $\theta\in N^1(X)_{\mathbf{R}}$ has the envelope property.
\end{theorem}

For $\operatorname{char}k=0$ or $n\leq2$, this is \cite[Theorem~5.20]{BJ22}. We therefore assume $\operatorname{char}k=p>0$. The proof occupies the rest of this section. Starting from the base-ideal formula for a PL envelope, we construct a decreasing approximation using asymptotic test ideals on the smooth scheme $S\coloneqq X\times\mathbf A^1$. A relative version of Musta\c{t}\u{a}'s generation theorem controls the classes of the approximants (\cref{lem:smooth-uniform-twist}), and a Frobenius-trace bound identifies their divisorial limit (\cref{lem:divisorial-trace-estimate}). For the approach using test ideals on a resolution of the test configuration, see \cite[\S2.1]{BJaddendum} and \cite{GJKM19}.

\subsection*{Reduction}

Choose rational ample classes $\theta_j\to\theta$ with $\theta_j-\theta$ ample, which is possible since the ample cone is open and rational classes are dense. By \cite[Lemma~5.9]{BJ22}, it suffices to prove $\EP_X(\theta_j)$ for all $j$. We may therefore assume $\theta=c_1(L)$ with $L\in\Pic(X)_{\mathbf Q}$ ample, and by \cref{prop:ample-envelopes}(1) it suffices to show that $\mathrm P_L(\phi)$ is $L$-psh for every rational PL function $\phi$. Since $\mathrm P_{mL}(m\phi)=m\mathrm P_L(\phi)$ and $\mathrm P_L(\phi+c)=\mathrm P_L(\phi)+c$, the discussion before \cref{prop:ample-envelopes} allows us to assume that $L$ is a very ample line bundle and that $\phi=\phi_D$ for a test configuration
\[
 \mu\colon\mathcal{Y}\longrightarrow S,
 \qquad
 \mathcal{L}=\mu^*L+D,
 \qquad D\leq0,
\]
where $\mathcal Y$ is normal, $\mu$ is projective, birational and an isomorphism outside the central fibre, and $D$ is a Cartier divisor supported on the central fibre. Line bundles on $X$ and their pullbacks to $S$ are denoted by the same symbols. Let $\mathfrak a_m\subseteq\mathcal O_{\mathcal Y}$ be the base ideal of $m\mathcal L$. By \cref{prop:ample-envelopes}(2),
\begin{equation}\label{eq:smooth-BJ-envelope}
 \mathrm Q\coloneqq\mathrm P_L(\phi_D)=\phi_D+\sup_{m\ge1}\frac1m\phi_{\mathfrak a_m}=\phi_D+\lim_{m\to\infty}\frac1m\phi_{\mathfrak a_m}.
\end{equation}
In particular $\mathrm Q\le\phi_D$. We shall construct a decreasing sequence of psh functions whose limit agrees with $\mathrm Q$ on $X^{\divv}$.

\subsection*{The graded ideals \texorpdfstring{$J_m$}{Jm}}

Since $S$ is normal and $\mu$ is proper and birational, $\mu_*\mathcal O_{\mathcal Y}=\mathcal O_S$. As $D\le0$, this gives $H^0(\mathcal Y,m\mathcal L)\subseteq H^0(\mathcal Y,\mu^*(mL))=H^0(S,mL)$. Let $J_m\subseteq\mathcal O_S$ be the ideal such that $J_m\otimes\mathcal O_S(mL)$ is the subsheaf of $\mathcal O_S(mL)$ generated by $H^0(\mathcal Y,m\mathcal L)$. Then
\[
J_mJ_r\subseteq J_{m+r},
\]
every power $J_q^r\otimes\mathcal O_S(qrL)$ is globally generated, and $J_m\neq0$: for $0\ne s\in H^0(X,mL)$ and $N$ large, $t^Ns\in H^0(\mathcal Y,m\mathcal L)$, since $-mD$ is supported on the central fibre, where $t$ vanishes. Both $(J_m\mathcal O_{\mathcal Y})\otimes\mu^*\mathcal O_S(mL)$ and $\mathfrak a_m\otimes\mathcal O_{\mathcal Y}(m\mathcal L)$ are the subsheaf of $\mu^*\mathcal O_S(mL)$ generated by $H^0(\mathcal Y,m\mathcal L)$. Hence
\[
 J_m\cdot\mathcal{O}_{\mathcal{Y}}=\mathfrak a_m\cdot\mathcal{O}_{\mathcal{Y}}(mD)
\]
as ideals of $\mathcal O_{\mathcal Y}$, where $\mathcal{O}_\mathcal{Y}(mD)\subseteq\mathcal{O}_\mathcal{Y}$ is the ideal of the effective divisor $-mD$. Evaluating $\sigma(v)$ on both sides gives $\frac1m\phi_{J_m}=\phi_D+\frac1m\phi_{\mathfrak a_m}$ on $X^{\an}$, and \eqref{eq:smooth-BJ-envelope} becomes
\begin{equation}\label{eq:smooth-Q-limit}
 \mathrm Q=\sup_{m\ge1}\frac1m\phi_{J_m}=\lim_{m\to\infty}\frac1m\phi_{J_m}.
\end{equation}

\subsection*{Asymptotic test ideals}

The smooth scheme $S$ is $F$-finite since $k$ is perfect. For a nonzero coherent ideal $I\subseteq\mathcal O_S$ and $c\in\mathbf Q_{\geq0}$, its test ideal $\tau(I^c)$ is given by the stable Frobenius-trace image
\begin{equation}\label{eq:test-def}
 \operatorname{Tr}^e\!\left(F_*^e\bigl(I^{\lceil cp^e\rceil}\otimes\omega_S\bigr)\right)
       =\tau(I^c)\otimes\omega_S\qquad(e\gg0),
\end{equation}
where $F$ is Frobenius and $\operatorname{Tr}^e:F_*^e\omega_S\to\omega_S$ is its Grothendieck trace \cite[\href{https://stacks.math.columbia.edu/tag/0FK6}{Tag 0FK6}]{Stacks}; see \cite[\S3]{MustataTest} for this description and its stabilization.

For $m\geq1$, the \emph{asymptotic test ideal} $B_m\coloneqq\tau(J_\bullet^m)$ equals $\tau(J_q^{m/q})$ for every sufficiently divisible $q$, and
\begin{equation}\label{eq:smooth-test-inclusions}
                  J_m\subseteq B_m,
 \qquad
                  B_{m\ell}\subseteq B_m^\ell\qquad(\ell\ge1),
\end{equation}
the second inclusion being the subadditivity of asymptotic test ideals. See Hara--Yoshida \cite{HaraYoshida} and Blickle--Musta\c{t}\u{a}--Smith \cite{BMS08} for the local theory, \cite{SchwedeTucker} for a survey, \cite[\S3, Proposition~3.1]{MustataTest} for the asymptotic construction in the form used here, and \cite[\S4]{GJKM19} for the compatibility between them. Put
\[
 u_m\coloneqq\frac1m\phi_{B_m}\qquad(m\ge1).
\]
Since $I\subseteq J$ implies $\phi_I\leq\phi_J$, \eqref{eq:smooth-test-inclusions} gives $\frac1{m\ell}\phi_{J_{m\ell}}\le u_{m\ell}\le u_m$ for $\ell\ge1$. Letting $\ell\to\infty$ and using \eqref{eq:smooth-Q-limit} and $B_m\subseteq\mathcal O_S$, we obtain
\begin{equation}\label{eq:smooth-Q-test-sandwich}
                         \mathrm{Q}\leq u_m\leq0
                         \qquad\text{on }X^{\an}.
\end{equation}

\subsection*{A uniform twist}

Fix a very ample line bundle $H$ on $X$ and put
\[
 A\coloneqq K_X+(n+1)H.
\]
The following relative version of Musta\c{t}\u{a}'s uniform generation theorem \cite[Theorem~A]{MustataTest} supplies a twist independent of $m$.

\begin{lemma}\label{lem:smooth-uniform-twist}
For every $m\geq1$, the sheaf $B_m(mL+A)$ is globally generated relative to $\mathbf A^1$, and
\begin{equation}\label{eq:smooth-FS-class}
 u_m\in\mathcal{H}^{\mathrm{gf}}_{\mathbf{Q}}\left(X,L+\frac1mA\right).
\end{equation}
\end{lemma}

\begin{proof}
Fix $m$ and choose $q$ sufficiently divisible that $B_m=\tau(J_q^{m/q})$.  Put $r_e=\lceil mp^e/q\rceil$. For $e\gg0$, \eqref{eq:test-def} gives a surjection
\[
 F_*^e(J_q^{r_e}\otimes\omega_S)
       \twoheadrightarrow B_m\otimes\omega_S.
\]
Since $J_q^{r_e}(qr_eL)$ is globally generated, it is a quotient of a finite sum of copies of $\mathcal O_S$.  Apply the exact functor $F_*^e$ to the corresponding surjection onto $J_q^{r_e}\otimes\omega_S$, compose with trace, and tensor by $mL+(n+1)H$.  The projection formula gives a surjection
\begin{equation}\label{eq:smooth-root-source}
 \bigoplus F_*^e\!\left(\mathcal O_S
 \bigl(K_S+(mp^e-qr_e)L+p^e(n+1)H\bigr)\right)
 \twoheadrightarrow B_m(mL+A).
\end{equation}
Here $K_S=\operatorname{pr}_X^*K_X$ since $K_{\mathbf A^1}$ is trivial, so the displayed twists are exact line-bundle identities.

For $1\leq i\leq n$, the higher direct image of a summand in \eqref{eq:smooth-root-source} twisted by $-iH$ is
\[
 F_{\mathbf A^1,*}^e
 \left(H^i\bigl(X,K_X+a_eL+p^e(n+1-i)H\bigr)\otimes_k k[t]\right),
 \qquad
 a_e\coloneqq mp^e-qr_e\in\{1-q,\ldots,0\}.
\]
This uses $\rho\circ F_S^e=F_{\mathbf A^1}^e\circ\rho$ for the projection $\rho\colon S\to\mathbf A^1$ and the cohomology formula for the product projection.  Since $n+1-i\geq1$, Serre vanishing \cite[Chapter~III, Theorem~5.2]{Hartshorne} for the finitely many residual bundles $K_X+a_eL$ kills these groups for all sufficiently large $e$; for $i>n$ they vanish by dimension.  Relative Castelnuovo--Mumford regularity \cite[Example~1.8.24]{LazarsfeldI} shows that the source, and hence its quotient, is relatively globally generated.  The necessary $e$ may depend on $m$ and $q$, but $A$ is fixed.

Over the affine base, choose finitely many global generators $g_j=\sum_\lambda s_{j\lambda}t^\lambda$ of $B_m(mL+A)$, where $s_{j\lambda}\in H^0(X,mL+A)$ and $\lambda\in\mathbf N$. The Gauss-extension formula yields
\[
 u_m=\frac1m\max_{j,\lambda}
             \{\log|s_{j\lambda}|-\lambda\},
\]
which is \eqref{eq:smooth-FS-class}, also at common zeros.
\end{proof}

\subsection*{A divisorial bound for test ideals}

Since $S$ is smooth, the \emph{log discrepancy} of a nontrivial divisorial valuation $w=r\ord_E$ on $S$, with $E$ a prime divisor on a normal birational model $g\colon V\to S$ and $r\in\mathbf Q_{>0}$, is
\[
 A_S(w)\coloneqq r\bigl(1+\ord_E(K_{V/S})\bigr)\in\mathbf{Q}_{>0};
\]
it is independent of the model and homogeneous of degree one in $w$. We only need a valuative bound at divisorial valuations. It is a special case of the discrepancy bounds in \cite{Canton}, but it also admits the following direct proof using the Frobenius trace on top differential forms; only smoothness of $S$ over the perfect field $k$ is used.

\begin{lemma}\label{lem:divisorial-trace-estimate}
Let $w$ be a nontrivial divisorial valuation with center on $S$, let $I\subseteq\mathcal{O}_S$ be a nonzero ideal and let $c\in\mathbf Q_{\geq0}$.  Then
\begin{equation}\label{eq:divisorial-test-bound}
 w\bigl(\tau(I^c)\bigr)\geq c\,w(I)-A_S(w).
\end{equation}
\end{lemma}

\begin{proof}
Fix a trivialization $\eta$ of $\omega_S$ near the center and write $\operatorname{Tr}^e(f\eta)=\Phi_e(f)\eta$.  We first prove that, for every $e\geq1$ and every local function $f$,
\begin{equation}\label{eq:divisorial-trace-bound}
 w(\Phi_e(f))\geq\frac{w(f)}{p^e}
               -\left(1-\frac1{p^e}\right)A_S(w).
\end{equation}
First take $w=\ord_E$ and choose a normal birational model $g\colon V\to S$ carrying $E$.  Then $R=\mathcal O_{V,E}$ is a DVR, so $V$ is regular at the generic point of $E$ and hence, the regular locus being open and $k$ perfect, smooth over $k$ near it; the whole computation takes place in $R$.  If $\eta$ generates $\omega_S$ near the center, write
\[
 g^*\eta=h\xi,\qquad
 a\coloneqq\ord_E(h)=A_S(\ord_E)-1\geq0,
\]
with $\xi$ a generator of top differential forms at $R$; the inequality holds because $g^*\eta$ is a regular top form there.  For a rational top form $\alpha$ we write $\ord_E(\alpha)$ for $\ord_E(\alpha/\xi)$, so that $\ord_E(fg^*\eta)=\ord_E(f)+a$.  On rational top forms Frobenius trace is the iterated Cartier operator: being compatible with localization, the traces of $S$ and of $V$ both restrict at the common generic point to the $e$-th Cartier operator of $\Omega^{\dim S}_{K/k}$, $K=k(S)=k(V)$, which depends on $K$ alone; and $\operatorname{Tr}^e\colon F^e_*\omega_R\to\omega_R$ exists because $R$ is regular and $F$-finite, so regular top forms are preserved. See \cite[\S3]{MustataTest} and \cite{SchwedeTucker} for the trace--Cartier description and its local coordinate formula.

Let $z$ be a uniformizer and put $\ell=\lfloor(\ord_E(f)+a)/p^e\rfloor$.  Then $fg^*\eta=z^{p^e\ell}\omega$ for a regular top form $\omega$ at $R$. Semilinearity and regularity of trace give
\[
 \operatorname{Tr}^e(fg^*\eta)=z^\ell\operatorname{Tr}^e(\omega),
 \qquad \ord_E\operatorname{Tr}^e(fg^*\eta)\geq\ell.
\]
Since $\operatorname{Tr}^e(fg^*\eta)=\Phi_e(f)\,g^*\eta$ has $\ord_E$ equal to $\ord_E\Phi_e(f)+a$, this yields
\[
 \ord_E\Phi_e(f)
 \geq\left\lfloor\frac{\ord_E(f)+a}{p^e}\right\rfloor-a
 \geq\frac{\ord_E(f)-(p^e-1)(a+1)}{p^e},
\]
which is \eqref{eq:divisorial-trace-bound} for $w=\ord_E$ because $(p^e-1)(a+1)=(p^e-1)A_S(\ord_E)$.  Both sides are homogeneous of degree one in $w$, which extends the bound to every nontrivial divisorial $w$.

Finally, fix $e\gg0$ and put $N=\lceil cp^e\rceil$. By \eqref{eq:test-def}, $\tau(I^c)=\{\Phi_e(f):f\in I^N\}$ near the center of $w$. Every $f\in I^N$ satisfies $w(f)\geq N\,w(I)\geq cp^e\,w(I)$, since $w(I)\geq0$, so \eqref{eq:divisorial-trace-bound} and $A_S(w)>0$ give $w(\Phi_e(f))\geq c\,w(I)-A_S(w)$. This proves \eqref{eq:divisorial-test-bound}.
\end{proof}

\subsection*{Conclusion of the proof}

Take $m_j=j!$. By \eqref{eq:smooth-test-inclusions} the sequence $(u_{m_j})$ decreases, and by \eqref{eq:smooth-Q-test-sandwich} it is bounded below by $\mathrm{Q}$, which is finite at $\vtriv$ because $\mathrm{Q}\geq\frac1m\phi_{J_m}$ and $J_m\neq0$. Thus \cref{lem:smooth-uniform-twist} and the defining decreasing Fubini--Study approximation \cite[Definition~4.1]{BJ22} give
\[
 u\coloneqq\lim_j u_{m_j}\in\PSH(X,L).
\]
Let $v\in X^{\divv}$. Then $w\coloneqq\sigma(v)$ is a divisorial valuation on $S$ (\cref{subsec:prelim-berkovich}), nontrivial since $w(t)=1$. Choose $q$ sufficiently divisible that $B_m=\tau(J_q^{m/q})$. Since $\mathrm Q\ge\frac1q\phi_{J_q}$, that is, $w(J_q)\geq-q\,\mathrm Q(v)$, \cref{lem:divisorial-trace-estimate} gives $w(B_m)\geq-m\,\mathrm Q(v)-A_S(w)$. Together with \eqref{eq:smooth-Q-test-sandwich}, this gives
\[
 \mathrm{Q}(v)\leq u_m(v)\leq\mathrm{Q}(v)+\frac{A_S(w)}m.
\]
Hence $u=\mathrm{Q}$ on $X^{\divv}$. Since $\mathrm Q\le\phi_D$ and $\phi_D$ is continuous, divisorial comparison gives $u\le\phi_D$ on $X^{\an}$, so $u\le\mathrm P_L(\phi_D)=\mathrm Q\le u$ by \eqref{eq:smooth-Q-test-sandwich}. Thus $\mathrm P_L(\phi_D)=u$ is $L$-psh, which completes the proof of \cref{thm:smooth-EP}.\qed

\section{Finite descent}
\label{sec:finite-descent}

Let $h\colon T\to X$ be a finite surjective morphism of integral projective varieties, with $X$ normal. Field norms of sections of pulled-back line bundles push Fubini--Study functions on $T$ forward to $X$; normality of $X$ makes the norm coefficients regular, while the source $T$ need not be normal. Divisorial comparison supplies the monotonicity needed for decreasing limits. This gives finite descent of the envelope property (\cref{thm:finite-normal-descent}).

Write $E=k(T)$, $K=k(X)$ and $d=[E:K]$.  For $v\in X^{\divv}$, including $\vtriv$, the support of $v$ is generic, so $K\hookrightarrow\mathcal H(v)$.  Decompose the finite algebra
\[
 E\otimes_K\mathcal H(v)=\prod_{w\mid v}B_w,
 \qquad d(w/v)\coloneqq\dim_{\mathcal H(v)}B_w,
 \qquad\sum_{w\mid v}d(w/v)=d.
\]
Here each $B_w$ is local Artinian; its residue field, finite over the complete field $\mathcal H(v)$, has a unique extending absolute value, and $E\hookrightarrow B_w$ followed by reduction is injective.  These give exactly the extensions $w$ of $v$ to $E$, with no separability assumption \cite[Ch.~VI, \S8]{AC}.

Every such extension is a point of $T^{\divv}$.  Indeed, if $v=c\ord_D$ on a normal projective model $X'\to X$, normalize $X'$ in $E$.  This is a finite normal model over $X'$ and maps birationally to $T$; the extensions of $v$ are positive rational multiples of the orders along the prime divisors over $D$.  Over $\vtriv$ the only extension is the trivial valuation of $E$.

\begin{proposition}\label{prop:weighted-coefficient-norm}
Let $L\in\Pic(X)_{\mathbf Q}$ and let
\[
 \varphi=\frac1m\max_j\{\log|s_j|+\lambda_j\}
 \in\mathcal{H}^{\mathrm{gf}}_{\mathbf{Q}}(T,h^*L),\qquad
 s_j\in H^0(T,h^*(mL)).
\]
The universal field norm of this system has coefficients $c_I\in H^0(X,mdL)$, $|I|=d$, not all zero.  They define
\begin{equation}\label{eq:weighted-norm-definition}
 \Av_h(\varphi)\coloneqq\frac1{md}\max_{|I|=d}
 \left\{\log|c_I|+\sum_j I_j\lambda_j\right\}
 \in\mathcal{H}^{\mathrm{gf}}_{\mathbf{Q}}(X,L).
\end{equation}
At every divisorial point one has
\begin{equation}\label{eq:divisorial-field-average}
 \Av_h(\varphi)(v)
   =\frac1d\sum_{w\mid v}d(w/v)\varphi(w).
\end{equation}
The construction is independent of the presentation and preserves pointwise order, including between functions presented in different line bundles.
\end{proposition}

\begin{proof}
On an affine $U=\Spec A\subset X$ trivializing $mL$ with frame $e$, write $s_j=f_jh^*e$ with $f_j\in B\coloneqq\mathcal O_T(h^{-1}U)$. Introduce independent variables $\mathbf T=(T_j)$ and form
\[
 P_U(\mathbf T)
 =\Norm_{E(\mathbf T)/K(\mathbf T)}\left(\sum_jT_jf_j\right)
 =\sum_{|I|=d}c_{I,U}\mathbf T^I.
\]
This is the determinant of multiplication on the free $K[\mathbf T]$-module $E[\mathbf T]$, so it is homogeneous of degree $d$.  The linear form is integral over $A[\mathbf T]$, and its norm is integral over that ring and lies in $K[\mathbf T]$.  Since $A[\mathbf T]$ is normal, $P_U\in A[\mathbf T]$. Replacing $e$ by $ge$ multiplies $P_U$ by $g^{-d}$; hence its coefficients glue to sections of $mdL$.  At the generic point the linear form is nonzero and has nonzero field norm, so some $c_I\ne0$. Thus \eqref{eq:weighted-norm-definition} defines a Fubini--Study function in the asserted class.

Fix $v\in X^{\divv}$, take $U$ containing its center, and extend $\mathcal H(v)$ isometrically to a complete algebraically closed field $C$.  The finite $C$-algebra $E\otimes_K C$ is a product of local Artin algebras with residue field $C$.  On a factor of length $a$, multiplication by $\sum_jT_jf_j$ has determinant equal to the $a$-th power of its residue linear form: on the filtration by powers of the maximal ideal it acts by that residue on each successive quotient.  All embeddings of a residue field of $B_w$ into $C$ are isometric, and the local factors of $E\otimes_KC$ lying over $w$ have total $C$-dimension $d(w/v)$, determinants being compatible with the base change $K[\mathbf T]\to C[\mathbf T]$.  The \emph{weighted Gauss norm} of $\sum_Iq_I\mathbf T^I\in C[\mathbf T]$ is $\max_I\{|q_I|\prod_j\mathrm{e}^{\lambda_jI_j}\}$, equivalently the supremum norm on the Berkovich polydisc of polyradius $(\mathrm e^{\lambda_j})_j$; it is multiplicative, and therefore yields
\[
 \max_{|I|=d}\left\{\log|c_I|(v)+\sum_jI_j\lambda_j\right\}
 =\sum_{w\mid v}d(w/v)
       \max_j\{\log|s_j|(w)+\lambda_j\}.
\]
Dividing by $md$ proves \eqref{eq:divisorial-field-average}. Two presentations therefore give psh functions with the same divisorial values, so divisorial comparison shows that they give the same function on $X^{\an}$.

If $\varphi\leq\psi$ on $T^{\an}$, then $\Av_h(\varphi)\leq\Av_h(\psi)$ on $X^{\divv}$ by \eqref{eq:divisorial-field-average}, hence on all of $X^{\an}$ by divisorial comparison, the left-hand side being psh and the right-hand side usc; this holds for functions presented in different numerical classes, since \cref{prop:BJ-package}(2) requires nothing more.
\end{proof}

\begin{theorem}\label{thm:finite-normal-descent}
Let $h\colon T\to X$ be a finite surjective morphism of integral projective varieties, with $X$ normal, and $\theta\in N^1(X)_{\mathbf R}$. If $\EP_T(h^*\theta)$ holds, then so does $\EP_X(\theta)$.
\end{theorem}

\begin{proof}
Let $(\varphi_i)$ be a bounded-above increasing net in $\PSH(X,\theta)$. By the envelope property upstairs and \cref{prop:EP-criterion}, there is $F\in\PSH(T,h^*\theta)$ with $F=\sup_i\varphi_i\circ h^{\an}$ on $T^{\divv}$. By \cref{prop:BJ-package}(3) and \cref{lem:finite-class-promotion}, $F$ is the limit of a decreasing sequence $F_j\in\mathcal{H}^{\mathrm{gf}}_{\mathbf{Q}}(T,h^*L_j)$ with $c_1(L_j)\to\theta$. By \cref{prop:weighted-coefficient-norm}, the functions $\Av_h(F_j)\in\mathcal{H}^{\mathrm{gf}}_{\mathbf{Q}}(X,L_j)$ decrease on $X^{\an}$. The only extension of $\vtriv$ is the trivial valuation of $E$, so $\Av_h(F_j)(\vtriv)=F_j(\vtriv)\to F(\vtriv)\in\mathbf R$, and \cite[Definition~4.1]{BJ22} gives $G\coloneqq\lim_j\Av_h(F_j)\in\PSH(X,\theta)$ directly, since these approximants are themselves Fubini--Study functions with classes tending to $\theta$. For $v\in X^{\divv}$, every extension $w$ of $v$ lies in $T^{\divv}$ and satisfies $h(w)=v$, so passing to the limit in \eqref{eq:divisorial-field-average} gives
\[
 G(v)=\frac1d\sum_{w\mid v}d(w/v)F(w)=\frac1d\sum_{w\mid v}d(w/v)\sup_i\varphi_i(v)=\sup_i\varphi_i(v).
\]
\Cref{prop:EP-criterion} proves $\EP_X(\theta)$.
\end{proof}

\section{Local regularization on normal varieties}\label{sec:normal-local-regularization}

We prove that, over a countable field, a function on a normal projective variety is psh in an ample class whenever its pullback to a normal projective modification is psh (\cref{cor:normal-weak-global}). This is a non-Archimedean analogue of the theorem of Forn\ae ss--Narasimhan \cite{FN80}. Apart from this corollary, only \cref{lem:local-basic,lem:global-to-local} are used later.

The argument follows Wu \cite{Wu26}, with the modifications needed in arbitrary characteristic. On an affine chart, we consider \emph{regularizable} functions, defined as decreasing limits of affine Fubini--Study functions. Polynomial hulls and norms along generically \'etale finite projections extend regularizability from the complement of a closed subset (\cref{thm:normal-local-regularization}). For the global application, regularizability away from the exceptional locus follows directly from the psh pullback. We cite the local extension and hull arguments of \cite{Wu26}, recording minor corrections in footnotes. We give proofs where the passage to arbitrary characteristic or the form needed here requires a further argument.

\emph{Conventions.} Throughout this section $k$ is countable, algebraically closed and trivially valued. Closed subsets carry their reduced induced structure. Countability enters only through the countability of $k[Y]$ for an affine variety $Y$, in \cref{lem:polyhedral-neighborhoods}(2) and \cref{lem:strict-exhaustion}; in particular $Y^{\an}$ is second countable. For a variety $Y$, the set $Y^{\divv}\subseteq Y^{\an}$ of divisorial points consists of $\vtriv$ and of the valuations $c\ord_E$ of $k(Y)$ described in \cref{subsec:prelim-berkovich}; thus $Y^{\divv}=P^{\divv}$ if $Y$ is an open subset of a projective variety $P$. For a regular function $f$ we put $\log|f|(y)=-y(f)$. If $M$ is a line bundle on a projective variety $P$ and $0\ne s\in H^0(P,M)$, we write $P_s\subseteq P$ for the open subset where $s$ does not vanish. Thus $P_s^{\an}=\{\log|s|>-\infty\}$, and $\log|\sigma|-\log|s|=\log|\sigma/s|$ on $P_s^{\an}$ for every $\sigma\in H^0(P,M)$.

\subsection{Local psh functions}\label{subsec:local-definitions}

In this subsection $F$ is a complete non-Archimedean field, which will be $k$ or a complete extension of $k$, and $V$ is a scheme of finite type over $F$. Complete extensions of $k$ are used in \cref{lem:local-PL-positivity} and \cref{subsec:norms}.

\begin{definition}\label{def:local-psh}
Let $U\subseteq V^{\an}$ be open.
\begin{enumerate}
\item A function $u\colon U\to\mathbf R$ is \emph{convex PL} if every point of $U$ has an open neighbourhood $W\subseteq U$ with the following property: there are a Zariski open subset $V_0\subseteq V$ with $W\subseteq V_0^{\an}$, invertible functions $f_1,\ldots,f_r\in\mathcal O(V_0)^\times$ and rational numbers $a_{ij},b_j$ such that
\[
u=\max_{1\le j\le N}\Bigl(\sum_{i=1}^ra_{ij}\log|f_i|+b_j\Bigr)\qquad\text{on }W.
\]
\item A function on $U$ is \emph{locally PL} if it is locally the difference of two convex PL functions.
\item A function $u\colon U\to[-\infty,\infty)$ is \emph{locally psh} if every point of $U$ has an open neighbourhood $W\subseteq U$ on which $u$ is the pointwise limit of a decreasing sequence of convex PL functions on $W$.
\end{enumerate}
\end{definition}

Convex PL functions are continuous, and locally psh functions are upper semicontinuous; the latter may be identically $-\infty$ on some open subsets. \Cref{def:local-psh} is a variant of \cite[Definitions~3.7 and~3.8]{Wu26}, using invertible functions and explicit finite maxima.\footnote{In \cite[Definition~3.1(3)]{Wu26}, the supporting affine function must touch the graph: replace the displayed condition by $u(x)=l(x)$ and $u\ge l$ near $x$. Our formulation also avoids conventions at infinite tropical coordinates.}

\begin{lemma}\label{lem:local-basic}
Let $U\subseteq V^{\an}$ be open.
\begin{enumerate}
\item Finite maxima, finite sums, positive rational multiples and rational translates of convex PL \textup(resp.\ locally psh\textup) functions on $U$ are convex PL \textup(resp.\ locally psh\textup). If $f$ is invertible on a Zariski open subset $V_0\subseteq V$, then adding a rational multiple of $\log|f|$ preserves both classes on $U\cap V_0^{\an}$.
\item Let $g\colon V'^{\an}\to V^{\an}$ be the analytification of a morphism $V'\to V$ of schemes of finite type over $F$, or the natural map $V_{F'}^{\an}\to V^{\an}$ for a complete extension $F'/F$. If $u$ is convex PL \textup(resp.\ locally psh\textup) on $U$, then so is $u\circ g$ on $g^{-1}(U)$.
\item Let $U'\subseteq U$ be open, $u$ locally psh on $U$ and $v$ locally psh on $U'$, and assume that the closure in $U$ of $\{x\in U':v(x)>u(x)\}$ is contained in $U'$. Then the function equal to $\max\{u,v\}$ on $U'$ and to $u$ on $U\setminus U'$ is locally psh on $U$, and it is convex PL if $u$ and $v$ are.
\item Every locally psh function is classically psh in the sense of \cite[Definition~4.1]{GublerRabinoff25}. In particular its pullback under a morphism from an analytic curve is subharmonic, and it satisfies the local maximum principle \cite[Theorem~6.1]{GublerRabinoff25}.
\end{enumerate}
\end{lemma}

\begin{proof}
(1) A sum of two finite maxima is the maximum of the pairwise sums, and maxima and sums of decreasing sequences decrease to the maxima and sums of the limits. Moreover, if $u_j$ decrease to $u$ and rational $q_j$ increase to $c\in\mathbf R$, then $u_j-q_j$ decrease to $u-c$; so locally psh functions are also stable under real translates. (2) Invertible functions pull back to invertible functions, and $\log|f|\circ g=\log|f\circ g|$. (3) On $U'$ this is (1). A point of $U\setminus U'$ has a neighbourhood disjoint from the closure in question, on which the function equals $u$. (4) For $f$ invertible and analytic, $\log|f|$ and $-\log|f|=\log|f^{-1}|$ are classically psh \cite[Proposition~4.2(2)]{GublerRabinoff25}, and so are constants. Hence convex PL functions are classically psh by the cone property, and so are locally psh functions, by stability under decreasing limits and the sheaf property \cite[Proposition~4.2(1),(3),(4)]{GublerRabinoff25}. On curves, classically psh functions are the subharmonic ones \cite[\S4.1.2]{GublerRabinoff25}.
\end{proof}

We now compare these notions with global psh functions. Let $P$ be a projective variety over $k$ and $L\in\Pic(P)_{\mathbf Q}$ semiample. We call $(m,s)$ an \emph{admissible pair} for $L$ if $m\ge1$, $mL$ is a globally generated line bundle and $0\ne s\in H^0(P,mL)$. Choose generators $\sigma_j$ of $mL$ with $s$ among them. Then $\max_j\log|\sigma_j|=0$, so
\begin{equation}\label{eq:reference-potential}
-\frac1m\log|s|=\frac1m\max_j\log|\sigma_j/s|\qquad\text{on }P_s^{\an}.
\end{equation}
Near a point $x\in P_s^{\an}$ the terms with $\sigma_j/s$ vanishing at the support of $x$ may be omitted, and the remaining $\sigma_j/s$ are invertible near that support. Hence $-m^{-1}\log|s|$ is convex PL on $P_s^{\an}$. If $(m',s')$ is another admissible pair, then $m^{-1}\log|s|-m'^{-1}\log|s'|=(mm')^{-1}\log|s^{m'}/s'^m|$, where $s^{m'}/s'^m$ is invertible on $P_s\cap P_{s'}$.

\begin{definition}\label{def:locally-L-psh}
Let $U\subseteq P^{\an}$ be open. A function $\varphi\colon U\to[-\infty,\infty)$ is \emph{locally $L$-psh} \textup(resp.\ \emph{locally $L$-convex}\textup) if $\varphi-m^{-1}\log|s|$ is locally psh \textup(resp.\ convex PL\textup) on $U\cap P_s^{\an}$ for every admissible pair $(m,s)$ for $L$.
\end{definition}

By the last observation and \cref{lem:local-basic}(1), it suffices that each point of $U$ lies in some $P_s^{\an}$ near which the condition holds for $(m,s)$.

\begin{lemma}\label{lem:global-to-local}
Let $P$ be a projective variety and $L\in\Pic(P)_{\mathbf Q}$ semiample.
\begin{enumerate}
\item Every $\varphi\in\mathcal{H}^{\mathrm{gf}}_{\mathbf{Q}}(P,L)$ defined by sections without common zero is locally $L$-convex.
\item Every $\Phi\in\PSH(P,L)$ is locally $L$-psh.
\end{enumerate}
\end{lemma}

\begin{proof}
(1) Write $\varphi=N^{-1}\max_\ell\{\log|\tau_\ell|+\lambda_\ell\}$ with $\tau_\ell\in H^0(P,NL)$ without common zero. On $P_s^{\an}$ we have $\varphi-m^{-1}\log|s|=(Nm)^{-1}\max_\ell\{\log|\tau_\ell^m/s^N|+m\lambda_\ell\}$, a maximum of logarithms of regular functions on $P_s$ plus constants. Near a point we may omit the terms that vanish at its support, and at least one term remains since the $\tau_\ell$ have no common zero.

(2) Since $L$ is nef, \cite[Corollary~12.18(ii)]{BJ22} gives $M_j\in\Pic(P)_{\mathbf Q}$ with $c_1(M_j)\to c_1(L)$ and $\Phi_j\in\mathcal{H}^{\mathrm{gf}}_{\mathbf{Q}}(P,M_j)$, defined by sections without common zero, decreasing to $\Phi$. Fix a globally generated ample line bundle $H$. Passing to a subsequence, choose rational $\epsilon_j\downarrow0$ such that $L+\epsilon_jH-M_j$ is ample. Tensoring the sections defining $\Phi_j$ with a generating family of a multiple of $L+\epsilon_jH-M_j$, whose homogeneous Fubini--Study function is zero, shows that $\Phi_j\in\mathcal{H}^{\mathrm{gf}}_{\mathbf{Q}}(P,L+\epsilon_jH)$, still defined by sections without common zero. Given $x\in P_s^{\an}$, choose $t\in H^0(P,H)$ not vanishing at the support of $x$ and put $W\coloneqq(P_s\cap P_t)^{\an}$. By (1), applied to $L+\epsilon_jH$ and the admissible pair given by a common power of $s$ and $t$, the functions $\Phi_j-m^{-1}\log|s|-\epsilon_j\log|t|$ are convex PL on $W$. They decrease to $\Phi-m^{-1}\log|s|$, because $-\log|t|\ge0$ and $\epsilon_j\downarrow0$.
\end{proof}

We write $\PL(P^{\an})$ for the $\mathbf Q$-vector space of rational PL functions on $P^{\an}$ \cite[\S2]{BJ22}; it is stable under finite maxima and minima, and it is dense in $C^0(P^{\an})$ \cite[Theorem~2.2]{BJ22}.

The following is a special case of \cite[Theorem~3.23]{Wu26}, which relies on the positivity criterion of Chambert-Loir and Ducros \cite{CLD12}. Since that criterion is formulated for nontrivially valued fields, we explain the passage from the trivially valued setting.

\begin{lemma}\label{lem:local-PL-positivity}
Let $P$ be a normal projective variety and $L\in\Pic(P)_{\mathbf Q}$ ample. Every locally $L$-convex function $\varphi\colon P^{\an}\to\mathbf R$ belongs to $\PSH(P,L)$.
\end{lemma}

\begin{proof}
By \eqref{eq:reference-potential}, $\varphi$ is locally PL. We first show that $\varphi$ agrees near each point $x$ with a function in $\PL(P^{\an})$. By the definition of locally PL functions, it suffices to treat $\log|f|$, where $f$ is a rational function that is regular and invertible at the support $\xi$ of $x$. Fix an ample line bundle $A$ and let $\mathfrak d$ be the ideal of denominators of $f$; it is trivial at $\xi$. For $m$ large choose $s_0\in H^0(P,mA\otimes\mathfrak d)$ with $s_0(\xi)\ne0$; then $s_1\coloneqq fs_0\in H^0(P,mA)$ and $s_1(\xi)\ne0$. Complete $s_0,s_1$ to a generating family $(s_j)$ of $mA$. Choose one rational $\lambda>\max_{i=0,1}\{-\log|s_i|(x)\}$. For $i=0,1$, the Fubini--Study function $\psi_i\coloneqq m^{-1}\max\{\log|s_i|+\lambda,\max_{j\ne i}\log|s_j|\}$ equals $m^{-1}(\log|s_i|+\lambda)$ near $x$, since $\log|s_j|\le0$ for all $j$. Hence $\log|f|=m(\psi_1-\psi_0)$ near $x$.

By compactness, $P^{\an}$ is covered by the interiors of finitely many compact sets $Q_i$ such that $\varphi=G_i$ on an open neighbourhood $V_i$ of $Q_i$, with $G_i\in\PL(P^{\an})$. Choose continuous $\chi_i\colon P^{\an}\to[0,1]$ equal to $1$ on $Q_i$ with support in $V_i$, and $\tau_i\in\PL(P^{\an})$ with $|\tau_i-\chi_i|\le1/8$. Then $b_i\coloneqq\max\{0,\min\{1,2\tau_i-1/2\}\}\in\PL(P^{\an})$ takes values in $[0,1]$, equals $1$ on $Q_i$, and vanishes outside $\{\chi_i\ge1/8\}\subseteq V_i$. Let $M\in\mathbf Q$ bound all $|G_i|$ and $|\varphi|$, which are continuous, and put $G\coloneqq\min_i\{G_i+4M(1-b_i)\}\in\PL(P^{\an})$. Let $y\in\operatorname{int}Q_{i_0}$. The term of index $i_0$ equals $\varphi(y)$. If $b_i(y)>0$, then $y\in V_i$ and the term of index $i$ is at least $\varphi(y)$; if $b_i(y)=0$, it is at least $-M+4M\ge\varphi(y)$. Thus $\varphi=G\in\PL(P^{\an})$.

Replacing $L$ and $\varphi$ by a common positive multiple, we may write $\varphi=\phi_{\mathcal L}$ for a line bundle $\mathcal L$ on a normal test configuration $\mathcal P$ of $(P,L)$ dominating $P\times\mathbf A^1$ \cite[Theorem~2.31]{BJ22}. By \cite[Definition~3.2 and Lemma~4.3]{BJ22} it suffices to show that $\mathcal L$ has nonnegative degree on every integral curve $C$ contained in the central fibre $\mathcal P_0$.

Let $K=k((\varpi))$ with $|\varpi|=\mathrm e^{-1}$, let $\mathfrak P$ be the formal completion of $\mathcal P\times_{\mathbf A^1}\operatorname{Spec}k[[\varpi]]$ along its special fibre $\mathcal P_0$, and let $\mathfrak L$ be the induced formal line bundle. This is a formal model of $(P_K,L_K)$ and defines an integral PL metric on $L_K$. Let $r\colon P_K^{\an}\to P^{\an}$ be the restriction map. Then $\varphi\circ r$ is the model function of $\mathfrak L$ relative to the trivial model of $L$. Indeed, both sides are additive in $\mathcal L$, so it suffices to treat $\log|\mathfrak a|$ for a $\mathbf G_m$-invariant ideal $\mathfrak a=\sum_\lambda\mathfrak a_\lambda t^\lambda$ on $P\times\mathbf A^1$, which is generated by elements $a\,t^\lambda$ with $a\in\mathfrak a_\lambda$; at $w\in P_K^{\an}$ both sides then equal $\max_\lambda\{\log|\mathfrak a_\lambda|(r(w))-\lambda\}$. By hypothesis and \cref{lem:local-basic}(2), near every point of $P_K^{\an}$ the metric has a local potential of the form $h\circ(\log|f_1|,\ldots,\log|f_r|)$ with $f_i$ invertible and $h$ convex. Such a function is psh in the sense of \cite{CLD12} by \cite[Corollaire~18.8.9]{CLD12}, so the curvature current of the metric is positive \cite[Proposition~8.3.1]{CLD12}. By \cite[Th\'eor\`eme~18.8.4 and Corollaire~18.8.8]{CLD12}, the residual line bundle of $L_K$ is numerically effective at every point of $P_K^{\an}$. Choose $z\in P_K^{\an}$ whose reduction is the generic point of $C$. By \cite[\S17.7.9 and Lemme~17.7.10]{CLD12}, the pair consisting of $C$ and $\mathcal L|_C$ is a premodel of the residual line bundle at $z$; hence $\deg(\mathcal L|_C)\ge0$.
\end{proof}

Finally we recall the notion of weakly psh functions \cite[Definition~6.1]{Wu26}.

\begin{definition}\label{def:weakly-psh}
Let $Y$ be a variety over $k$ and $U\subseteq Y^{\an}$ open. A function $u\colon U\to[-\infty,\infty)$ is \emph{weakly psh} if there is a proper birational morphism $\mu\colon Y'\to Y$ such that $u\circ\mu^{\an}$ is locally psh on $(\mu^{\an})^{-1}(U)$. Such a $\mu$ is called a \emph{witness}.
\end{definition}

By Chow's lemma and normalization we may take $Y'$ normal and projective over $Y$. Since $\mu^{\an}$ is proper and surjective, weakly psh functions are upper semicontinuous. Locally psh functions are weakly psh, and restrictions of weakly psh functions to open subsets are weakly psh.

\subsection{Regularizable functions on affine varieties}\label{subsec:affine-regularizable}

In this subsection $Y$ is a normal affine variety. Regularizable functions on $Y$ are to psh functions on $Y^{\an}$ what psh functions of logarithmic growth are to psh functions on $\mathbf C^n$: they are the functions that extend, after twisting by a multiple of the boundary, to global psh functions on a projective compactification (\cref{lem:compactification}).

\begin{definition}\label{def:regularizable}
\begin{enumerate}
\item An \emph{affine Fubini--Study function} on $Y$ is a function of the form
\[
g=\max\{a_1\log|f_1|+b_1,\ldots,a_r\log|f_r|+b_r,\,b_0\},
\]
where $f_1,\ldots,f_r\in k[Y]$ include a set of generators of the $k$-algebra $k[Y]$, $a_i\in\mathbf Q_{>0}$ and $b_i\in\mathbf Q$. We write $\FS(Y)$ for the set of these functions.
\item A function $u$ on $Y^{\an}$ has \emph{logarithmic growth} if $u\le g+C$ for some $g\in\FS(Y)$ and $C\in\mathbf R$.
\item A function $u\colon Y^{\an}\to[-\infty,\infty)$ is \emph{regularizable} if it is not identically $-\infty$ and is the pointwise limit of a decreasing sequence in $\FS(Y)$. We write $\mathcal{R}(Y)$ for the set of regularizable functions.
\end{enumerate}
\end{definition}

The class $\mathcal{R}(Y)$ is denoted $\PSH^{\log}_{\mathrm{reg}}(Y^{\an})$ in \cite[Definition~3.15]{Wu26}.\footnote{The growth condition in \cite[Definition~3.15(2)]{Wu26} should read $u\le g+O(1)$, with $g\in\FS(Y)$, as used in its subsequent proofs.} Functions in $\FS(Y)$ are continuous, bounded below, and convex PL by the argument following \eqref{eq:reference-potential}, and their sublevel sets are compact, since they lie in closed polydiscs. Hence regularizable functions are locally psh and have logarithmic growth.

\begin{lemma}\label{lem:growth}
\begin{enumerate}
\item Let $\pi\colon Y'\to Y$ be a morphism from a normal affine variety \textup(for instance the inclusion of an affine open subset, or the identity\textup). For $g\in\FS(Y)$ and $g'\in\FS(Y')$ there are constants $C,C'$ with $g\circ\pi\le Cg'+C'$.
\item Let $\pi\colon Y\to S$ be a finite surjective morphism to a normal affine variety $S$. For $g\in\FS(Y)$ and $h\in\FS(S)$ there are constants $C,C'$ with $g\le C\,h\circ\pi+C'$.
\end{enumerate}
\end{lemma}

\begin{proof}
These are the estimates of \cite[Lemma~3.16]{Wu26}.\footnote{In its proof, for $f=\sum_\alpha a_\alpha z^\alpha$, use $\log|f|\le\max_\alpha\log|a_\alpha z^\alpha|$. For a monic relation $f^m+\sum_{i<m}b_if^i=0$, the root bound is $\log|f|\le\max_{i<m}(m-i)^{-1}\log|b_i|$, with the coefficient functions pulled back to $Y$.} The polynomial estimate in (1) applies to the finitely many pullbacks of functions defining $g$, expressed in the generators occurring in $g'$. For (2), apply the integral-equation estimate to each function defining $g$, and then apply (1) to its coefficients on $S$. Both arguments are independent of the characteristic.
\end{proof}

By (1), the restriction of a function of logarithmic growth to an affine open subset has logarithmic growth; we use this without further comment. Sums and finite maxima of functions of logarithmic growth have logarithmic growth, since $\FS(Y)$ is stable under sums and finite maxima.

\begin{lemma}\label{lem:boundary-presentation}
For every $g\in\FS(Y)$ there are a normal projective variety $P$ containing $Y$ as an open subset, an ample line bundle $H$ on $P$, a section $s\in H^0(P,H)$ with $P_s=Y$, a function $a\in\mathcal{H}^{\mathrm{gf}}_{\mathbf{Q}}(P,H)$ defined by sections without common zero, and $\kappa\in\mathbf Q_{>0}$, such that $g=\kappa(a-\log|s|)$ on $Y^{\an}$.
\end{lemma}

\begin{proof}
Apply \cite[Proposition~3.17(2)]{Wu26}, absorbing its additive constant into the weights of $a$.\footnote{The displayed sections in that proof are a generating family, not necessarily a basis; this suffices to define the Fubini--Study metric.} Normalizing the compactification preserves $Y$ and ampleness, and gives the stated normal $P$.
\end{proof}

The next lemma relates regularizable functions to global psh functions. In (2), decreasing convex PL presentations near compact subsets suffice; this is the form needed in \cref{lem:hartogs}.

\begin{lemma}\label{lem:compactification}
Let $P$ be a normal projective variety, $H$ an ample line bundle on $P$ and $s\in H^0(P,H)$ nonzero, and put $Y\coloneqq P_s$.
\begin{enumerate}
\item Let $c\in\mathbf Q_{>0}$ and $\Phi\in\PSH(P,cH)$. There are $u_j\in\FS(Y)$ decreasing to $\Phi-c\log|s|$ on $Y^{\an}$ such that each $u_j+c\log|s|$ is bounded. In particular $(\Phi-c\log|s|)|_{Y^{\an}}\in\mathcal{R}(Y)$.
\item Let $u\colon Y^{\an}\to[-\infty,\infty)$ be a function, not identically $-\infty$, such that $u\le-\kappa\log|s|+C$ for some $\kappa\in\mathbf Q_{>0}$ and $C\in\mathbf R$. Assume that every compact subset of $Y^{\an}$ has an open neighbourhood on which $u$ is the limit of a decreasing sequence of convex PL functions. Then for every rational $c>\kappa$, the function equal to $u+c\log|s|$ on $Y^{\an}$ and to $-\infty$ on $P^{\an}\setminus Y^{\an}$ belongs to $\PSH(P,cH)$.
\end{enumerate}
\end{lemma}

\begin{proof}
(1) Apply \cite[Corollary~3.18 and Proposition~3.17(1)]{Wu26}.\footnote{The limit in the proof of Corollary~3.18 is $g+\varphi$. In Proposition~3.17(1), replace the constant $0$ in the maximum by a lower bound $-C_j$. To include generators $z=\sigma/s^N$ of $k[Y]$, add the dominated terms $\epsilon\log|z|-C_j$ with $0<\epsilon\le c/N$; here $u_j\ge-c\log|s|-C_j$ and $\log|z|\le-N\log|s|$.} The fixed-class regularization used there is available by \cite[Corollary~12.18(iii) and Example~3.9]{BJ22}: take $\Phi_j\in\mathcal{H}^{\mathrm{gf}}_{\mathbf{Q}}(P,cH)$ defined by sections without common zero, decreasing to $\Phi$, and put $u_j=\Phi_j-c\log|s|$. Then $u_j+c\log|s|=\Phi_j$ is bounded.

(2) Use the truncation and gluing argument of \cite[Lemmas~3.25(1) and~3.26]{Wu26}, with \cref{lem:local-PL-positivity} as the positivity input. Fix $C_1\in\mathbf Q$. Since $c>\kappa$, choose $t_0$ so that $u+c\log|s|\le C_1-2$ where $\log|s|\le t_0$. Take convex PL approximants $u_j\downarrow u$ near the compact set $\{\log|s|\ge t_0-1\}$. On the compact collar $t_0-1\le\log|s|\le t_0$, one has $u_j+c\log|s|<C_1-1$ for all large $j$.\footnote{This step uses compactness: the open sets where the strict inequality holds increase and cover the collar. Continuity of the limit, as required by the cited use of Dini's theorem in \cite[Lemma~3.25(1)]{Wu26}, is unnecessary.} Thus $\max\{u_j+c\log|s|,C_1\}$ glues across the collar to the constant $C_1$ near the boundary of $Y$. The glued function is locally $cH$-convex, hence $cH$-psh. First let $j\to\infty$, then $C_1\to-\infty$, using PSH2 and finiteness of $u$ somewhere.
\end{proof}

We use the following closure properties.

\begin{lemma}\label{lem:regularizable-closure}
\begin{enumerate}
\item Every $u\in\mathcal{R}(Y)$ is finite on $Y^{\divv}$.
\item If $u_i\in\mathcal{R}(Y)$ decrease to a function $u\not\equiv-\infty$, then $u\in\mathcal{R}(Y)$.
\item Let $f_1,\ldots,f_r\in k[Y]$, $a_i\in\mathbf Q_{>0}$ and $b_i\in\mathbf Q$. If $q\coloneqq\max_i(a_i\log|f_i|+b_i)$ is not identically $-\infty$, then $q\in\mathcal{R}(Y)$. The class $\mathcal{R}(Y)$ is stable under finite maxima, finite sums, multiplication by positive rational numbers and addition of real constants.
\item Let $\pi\colon Y'\to Y$ be a morphism from a normal affine variety. For $u\in\mathcal{R}(Y)$, the function $u\circ\pi$ is either identically $-\infty$ or in $\mathcal{R}(Y')$. In particular, restriction to affine open subsets preserves regularizability.
\item Let $u\in\mathcal{R}(Y)$ and $g\in\FS(Y)$ with $u\le g+C$. For every rational $c>1$ there are $u_j\in\FS(Y)$ decreasing to $u$ such that each $u_j-cg$ is bounded.
\item Let $v\colon Y^{\an}\to[-\infty,\infty)$ be a function of logarithmic growth, not identically $-\infty$, such that every compact subset of $Y^{\an}$ has an open neighbourhood on which $v$ is the limit of a decreasing sequence of convex PL functions. Then $v\in\mathcal{R}(Y)$.
\end{enumerate}
\end{lemma}

\begin{proof}
Let $v$ be as in (6), or $v\in\mathcal{R}(Y)$. Fix $g\in\FS(Y)$ with $v\le g+C$ and a presentation $g=\kappa(a-\log|s|)$ as in \cref{lem:boundary-presentation}. Since $a$ is bounded, $v\le-\kappa\log|s|+C'$, so \cref{lem:compactification}(2) shows that $v+c\log|s|\in\PSH(P,cH)$ for rational $c>\kappa$. This and \cref{lem:compactification}(1) prove (6). The positivity of the glued truncations in this argument is supplied by \cref{lem:local-PL-positivity}. Thus (6) applies to the local convex PL presentations in its statement; \cite[Lemma~3.26]{Wu26}, which assumes a single tropicalization and logarithmic growth, does not by itself give this form.

(1) Psh functions are finite on $P^{\divv}=Y^{\divv}$, and so is $\log|s|$. (2) and (5) are the compactification arguments of \cite[Lemma~3.25(3),(2)]{Wu26}, using \cref{lem:compactification} and PSH2. In (2), a growth bound for $u_1$ controls the whole sequence; in (5), use the class $c\kappa H$ in the presentation $g=\kappa(a-\log|s|)$. (3) Let $b\coloneqq\max\{0,\log|z_1|,\ldots,\log|z_n|\}$ for generators $z_l$ of $k[Y]$. The functions $\max\{q,b/j-j\}$ lie in $\FS(Y)$ and decrease to $q$ because $b\ge0$. Maxima and sums of approximants are approximants of maxima and sums, since $a\log|f|+a'\log|f'|=N^{-1}\log|f^{Na}f'^{Na'}|$; the limits are not identically $-\infty$ by (1). For a real translate $u+c$, use approximants $u_j+c_j$ with rational $c_j\downarrow c$. (4) If $u_j\in\FS(Y)$ decrease to $u$ and $b'\ge0$ is as in (3) for $Y'$, then $\max\{u_j\circ\pi,b'/j-j\}\in\FS(Y')$ decrease to $u\circ\pi$.
\end{proof}

\begin{lemma}\label{lem:principal-extension}
Let $0\ne f\in k[Y]$ and $u\in\mathcal{R}(Y_f)$, and assume $u\le g+C$ on $Y_f^{\an}$ for some $g\in\FS(Y)$. For every $c\in\mathbf Q_{>0}$, the function equal to $u+c\log|f|$ on $Y_f^{\an}$ and to $-\infty$ on $V(f)^{\an}$ belongs to $\mathcal{R}(Y)$.
\end{lemma}

\begin{proof}
The truncation proof of \cite[Proposition~3.30]{Wu26} applies with \cref{lem:regularizable-closure}(5),(6),(2).\footnote{In that proof, read $u_\epsilon$ as $g_\epsilon=\max\{g,-\epsilon\log|f|\}$, with $\epsilon=c/3$, and retain the constant in $u_j\le2g_\epsilon+C_j$. If $h=\max\{2g+c\log|f|,(c/3)\log|f|\}$, the neighbourhood on which the truncation at $-j$ is constant is $\{h+C_j<-j\}$; it may depend on $j$.} Its truncations extend as convex PL functions of logarithmic growth on $Y$, and decrease to the asserted extension. The limit is finite at $\vtriv$ by \cref{lem:regularizable-closure}(1). No characteristic assumption is used.
\end{proof}

The following lemma replaces \cite[Proposition~6.11]{Wu26}. Its proof there uses norms of metrics \cite{BE21}, which are established under geometric reducedness assumptions, and a flattening argument. We use instead the field-norm construction of \cref{sec:finite-descent}.

\begin{lemma}\label{lem:finite-flat-norm}
Let $\pi\colon Y\to S$ be a finite flat morphism of normal affine varieties, of degree $d$. For $s\in S^{\an}$, write $\mathcal H(s)\otimes_{k[S]}k[Y]=\prod_yA_y$ as a product of local Artinian $\mathcal H(s)$-algebras, indexed by the points $y\in Y^{\an}$ with $\pi(y)=s$, and put $m_y(\pi)\coloneqq\dim_{\mathcal H(s)}A_y$, so that $\sum_{\pi(y)=s}m_y(\pi)=d$. For a function $u$ on $Y^{\an}$ let
\[
\Norm_\pi u(s)\coloneqq\sum_{\pi(y)=s}m_y(\pi)\,u(y).
\]
Then $\Norm_\pi$ maps $\mathcal{R}(Y)$ to $\mathcal{R}(S)$. If $F/k$ is a complete extension, $r$ denotes the natural maps from $F$-analytic spaces to $k$-analytic spaces, and $\Norm_{\pi_F}$ is defined by the same formula, then $(\Norm_\pi u)\circ r=\Norm_{\pi_F}(u\circ r)$ for $u\in\mathcal{R}(Y)$.
\end{lemma}

\begin{proof}
Let $u=N^{-1}\max_j\{\log|b_j|+\lambda_j\}\in\FS(Y)$, with $b_j\in k[Y]$ (the constant term corresponds to $b_j=1$). Since $k[Y]$ is a projective $k[S]$-module of rank $d$, the determinant of multiplication by $\sum_jT_jb_j$ on $k[Y][T]$, where $T=(T_j)$ are indeterminates, is a form $\sum_{|\alpha|=d}c_\alpha T^\alpha$ with $c_\alpha\in k[S]$. We claim that
\[
N\Norm_\pi u=\max_{|\alpha|=d}\Bigl\{\log|c_\alpha|+\sum_j\alpha_j\lambda_j\Bigr\}\qquad\text{on }S^{\an}.
\]
At a point $s$, this is the computation in the proof of \cref{prop:weighted-coefficient-norm}, with the finite $\mathcal H(s)$-algebra $\prod_yA_y$ in place of $E\otimes_K\mathcal H(v)$: after base change to a complete algebraically closed field, the determinant on a local Artinian factor of length $\ell$ is the $\ell$-th power of the residue of the linear form, and the weighted Gauss norm is multiplicative. No separability or reducedness of the fibres is needed. The right-hand side of the claim lies in $\mathcal{R}(S)$ by \cref{lem:regularizable-closure}(3).

For general $u\in\mathcal{R}(Y)$ with approximants $u_i\in\FS(Y)$, the functions $\Norm_\pi u_i$ decrease to $\Norm_\pi u$. The only point of $Y^{\an}$ over the trivial valuation of $k(S)$ is the trivial valuation of $k(Y)$, with multiplicity $d$; hence $\Norm_\pi u(\vtriv)=du(\vtriv)$ is finite by \cref{lem:regularizable-closure}(1), and \cref{lem:regularizable-closure}(2) gives $\Norm_\pi u\in\mathcal{R}(S)$. Finally, the claim for $\pi_F$ involves the same coefficients $c_\alpha$, which gives the compatibility with base change for $u\in\FS(Y)$, and hence for $u\in\mathcal{R}(Y)$ by passing to decreasing limits.
\end{proof}

\subsection{Polynomial convexity}\label{subsec:polynomial-convexity}

In this subsection $Y$ is again a normal affine variety. Polynomial hulls replace the pseudoconvexity arguments of complex analysis: the hull inequality for regularizable functions below is the analogue of the maximum principle, and the Runge property of sublevel sets is what allows us to recover a function from its Hartogs radius in \cref{lem:hartogs}. For a compact subset $K\subseteq Y^{\an}$, its \emph{polynomial hull} is
\[
\widehat K\coloneqq\{y\in Y^{\an}:|f(y)|\le\sup_K|f|\text{ for all }f\in k[Y]\}.
\]
We call $K$ \emph{polynomially convex} if $\widehat K=K$, and an open subset $D\subseteq Y^{\an}$ is \emph{Runge} if $\widehat K\subseteq D$ for every compact $K\subseteq D$. The hull $\widehat K$ is compact, being closed and contained in a closed polydisc. We say that a function $u$ on $Y^{\an}$ satisfies the \emph{hull inequality} if $u(y)\le\sup_Ku$ for every compact $K\subseteq Y^{\an}$ and every $y\in\widehat K$. Every $u\in\FS(Y)$ satisfies it, since each of its terms does. So does every $u\in\mathcal{R}(Y)$: if $u_j\in\FS(Y)$ decrease to $u$, then $u(y)\le u_j(y)\le\sup_Ku_j$, and $\sup_Ku_j$ decreases to $\sup_Ku$ because $K$ is compact and the $u_j$ are upper semicontinuous. An increasing union of Runge open subsets is Runge, because a compact subset of the union lies in one of them.

\begin{remark}
    For an affine scheme $Y$, let $A = \Gamma(Y^{\rm an}, \mathcal{O}_{Y^{\rm an}})$. For a compact subset $K\subset Y^{\rm an}$, we can also define a convex hull
    \begin{align*}
        \widehat{K}_{\rm an}:= \{ y\in Y^{\rm an}: |f(y)|\leq \sup_{K} |f| \text{ for all } f\in A\}.
    \end{align*}
    The convex hull $\widehat K_{\rm an}$ is the counterpart of the holomorphic convex hull used in \cite{FN80}. Essentially due to the fact that any element in $A$ can by approximated by elements in $k[Y]$, we have
    \begin{align*}
        \widehat K = \widehat K_{\rm an}.
    \end{align*}
    Indeed, the argument we use in this section is parallel to the one used in \cite{FN80}.
\end{remark}

We use the following strict-sublevel form of \cite[Proposition~5.4]{Wu26}.\footnote{Throughout the arguments from \cite[\S5]{Wu26}, use strict polynomial inequalities rather than identifying $\operatorname{int}\{g\le t\}$ with $\{g<t\}$. These sets can differ: $\{|T-1|<1\}\subseteq\{|T|\le1\}$ contains the rigid point $T=1$. The separation argument below gives the required strict inequalities directly.}

\begin{lemma}\label{lem:polyhedral-neighborhoods}
\begin{enumerate}
\item Let $H\subseteq Y^{\an}$ be a polynomially convex compact subset and $A$ an open neighbourhood of $H$. There is $g\in\FS(Y)$ with $g\ge0$, $H\subseteq\operatorname{int}\{g=0\}$ and $\{g\le1\}\subseteq A$.
\item Every Runge open subset $D\subseteq Y^{\an}$ is the union of an increasing sequence of open subsets $\{g_j<1\}$ with $g_j\in\FS(Y)$.
\end{enumerate}
\end{lemma}

\begin{proof}
(1) Let $z_1,\ldots,z_n$ generate $k[Y]$ and choose $R\in\mathbf Q$ with $H\subseteq\{\max_l\log|z_l|<R\}$. The set $E\coloneqq\{\max_l\log|z_l|\le R+1\}\setminus A$ is compact. For $y\in E$ we have $y\notin\widehat H$, so, replacing a function by a power, there are $f_y\in k[Y]$ and $a_y\in\mathbf Q$ with $\sup_H\log|f_y|<a_y<\log|f_y(y)|-1$. Finitely many of the open sets $\{\log|f_y|>a_y+1\}$, say for $y_1,\ldots,y_p$, cover $E$. Put
\[
g\coloneqq\max\Bigl\{0,\ \max_l(\log|z_l|-R),\ \max_{\ell}\bigl(\log|f_{y_\ell}|-a_{y_\ell}\bigr)\Bigr\}.
\]
All terms except $0$ are negative on the open set $\{\max_l\log|z_l|<R\}\cap\bigcap_\ell\{\log|f_{y_\ell}|<a_{y_\ell}\}$, which contains $H$, so $g=0$ there. If $y\notin A$, then either $\max_l\log|z_l|(y)>R+1$ or $y\in E$, and in both cases $g(y)>1$.

(2) Since $Y^{\an}$ is second countable and locally compact, $D$ is the union of an increasing sequence of compact sets $K_j$. Put $Q_0\coloneqq\varnothing$ and inductively let $g_j$ be given by (1) for the polynomially convex compact set $H_j\coloneqq\widehat{K_j\cup Q_{j-1}}$, which lies in $D$ because $D$ is Runge, and for $A=D$; put $Q_j\coloneqq\{g_j\le1\}$. Then $K_j\cup Q_{j-1}\subseteq H_j\subseteq\{g_j<1\}\subseteq Q_j\subseteq D$, which gives the assertion.
\end{proof}

We next use Rossi's local maximum principle.

\begin{lemma}\label{lem:rossi}
Let $K\subseteq Y^{\an}$ be compact, $x\in\widehat K$, and $V$ a relatively compact open neighbourhood of $x$ with $\overline V\cap K=\varnothing$. Then $x$ lies in the hull of $\partial V\cap\widehat K$.
\end{lemma}

\begin{proof}
Apply the proof of \cite[Theorem~5.6]{Wu26}, using \cref{lem:polyhedral-neighborhoods}(1) for its neighbourhood construction and choosing the scaling constant rational. The pasted maximum in that proof is convex PL and agrees with a function in $\FS(Y)$ outside a compact set. Hence it belongs to $\mathcal{R}(Y)$ by \cref{lem:regularizable-closure}(6), and the hull inequality gives the contradiction directly. These inputs apply in every characteristic.
\end{proof}

Let $\widetilde Y\coloneqq Y\times\mathbf A^1$, with fibre coordinate $w$ and projection $p\colon\widetilde Y\to Y$. The Hartogs radius will recover regularizability from the Runge property of $\{u\circ p+\log|w|<0\}$. For $y\in Y^{\an}$ the fibre $p^{-1}(y)\subseteq\widetilde Y^{\an}$ is identified with the Berkovich affine line over $\mathcal H(y)$, and $(y,0)$ denotes its origin. For an open subset $D\subseteq\widetilde Y^{\an}$ containing the zero section $\{(y,0):y\in Y^{\an}\}$, let $d_D(y)\in(0,+\infty]$ be the supremum of the radii $r$ such that the closed disc $\{|w|\le r\}$ of $p^{-1}(y)$ is contained in $D$.

\begin{lemma}\label{lem:hartogs}
Let $u\colon Y^{\an}\to[-\infty,\infty)$ be upper semicontinuous, of logarithmic growth and not identically $-\infty$. If $\{u\circ p+\log|w|<0\}\subseteq\widetilde Y^{\an}$ is Runge, then $u\in\mathcal{R}(Y)$.
\end{lemma}

\begin{proof}
Apply the proof of \cite[Proposition~5.7]{Wu26} to $D=\{u\circ p+\log|w|<0\}$, whose Hartogs radius is $d_D=\mathrm e^{-u}$. By \cref{lem:polyhedral-neighborhoods}(2), exhaust $D$ by increasing strict polyhedra $W_j=\{g_j<1\}\cap\{\log|w|<j\}$, and put $\Omega_j=\{y:(y,0)\in W_j\}$. The fibrewise calculation of \cite[Lemma~5.8]{Wu26}, made with these strict inequalities, shows that $-\log d_{W_j}$ is convex PL on $\Omega_j$.\footnote{If the defining inequalities are $\log|\sum_\ell a_{\alpha\ell}(y)w^\ell|<c_\alpha$, on $\Omega_j$ the Gauss norm gives $-\log d_{W_j}=\max\{-j,M_j\}$, where $M_j=\max_{\alpha,\ell\ge1}(\log|a_{\alpha\ell}|-c_\alpha)/\ell$. The constant term inequalities hold there by definition, and coefficients vanishing at a support are locally dominated by $-j$.} Compact fibre discs show that these radii increase to $d_D$. Since $\Omega_j\uparrow Y^{\an}$, the resulting decreasing approximation exists on a neighbourhood of every compact subset of $Y^{\an}$. Now apply \cref{lem:regularizable-closure}(6), as in Wu's argument.
\end{proof}

The following lemma is \cite[Proposition~5.10]{Wu26}. The proof given there rests on an incorrect ideal-theoretic claim \cite[Lemma~5.11]{Wu26}; we give a direct argument.

\begin{lemma}\label{lem:principal-hull}
Let $0\ne f\in k[Y]$, let $K\subseteq Y^{\an}$ be compact, and let $x\in\widehat K\cap Y_f^{\an}$, where $\widehat K$ is the hull in $Y$. If $v\in\mathcal{R}(Y_f)$ satisfies $v\le g+C$ for some $g\in\FS(Y)$, then $v(x)\le\sup_{K\cap Y_f^{\an}}v$.
\end{lemma}

\begin{proof}
Choose $M\ge0$ with $\log|f|\le M$ on $K$ and $|\log|f(x)||\le M$. The set $K\cap Y_f^{\an}$ is nonempty, since $x\in\widehat K$ and $|f(x)|>0$. For $c\in\mathbf Q_{>0}$, let $v_c$ be the extension of $v+c\log|f|$ by $-\infty$. By \cref{lem:principal-extension}, $v_c\in\mathcal{R}(Y)$, so the hull inequality gives
\[
v(x)-cM\le v_c(x)\le\sup_Kv_c\le\sup_{K\cap Y_f^{\an}}v+cM .
\]
The last supremum is bounded above by $\sup_Kg+C$. If it is finite, let $c\downarrow0$; if it is $-\infty$, the same inequality gives $v(x)=-\infty$.
\end{proof}

\subsection{A strictly psh exhaustion and gluing}\label{subsec:gluing}

In this subsection $Y$ is a normal affine variety. We need a regularizable function that admits supporting functions with a strict drop outside prescribed neighbourhoods, as in \cite[Lemmas~6.6 and~6.7]{Wu26}. We construct it using the countability of $k$ and a series of the functions $\log(1+|f|)$, $f\in k[Y]$, whose supporting lines are explicit.

Let $\mathbf T\coloneqq\mathbf R\cup\{+\infty\}$ with its natural topology, and let $\lambda\colon\mathbf T\to[0,\infty)$ be $\lambda(\tau)\coloneqq\log(1+\mathrm e^{-\tau})$, with $\lambda(+\infty)=0$. Thus $\lambda(-\log|f|)=\log(1+|f|)$ on $Y^{\an}$ for $f\in k[Y]$.

\begin{lemma}\label{lem:strict-exhaustion}
There is $\rho\in\mathcal{R}(Y)$ with $\rho\ge0$ such that for every $x\in Y^{\an}$ and every open neighbourhood $V$ of $x$ there are $\rho'\in\mathcal{R}(Y)$ and $\delta>0$ with
\[
\rho'\le\rho,\qquad\rho'(x)=\rho(x),\qquad\rho'\le\rho-\delta\ \text{on }Y^{\an}\setminus V.
\]
\end{lemma}

\begin{proof}
\emph{Supporting lines.} For $\tau^0\in\mathbf R$ let $p\coloneqq\mathrm e^{-\tau^0}/(1+\mathrm e^{-\tau^0})\in(0,1)$ and $\ell_{\tau^0}(\tau)\coloneqq\lambda(\tau^0)-p(\tau-\tau^0)\in[-\infty,\infty)$; for $\tau^0=+\infty$ let $\ell_{\tau^0}\coloneqq0$. We claim that $\lambda-\ell_{\tau^0}\ge0$ on $\mathbf T$, with equality only at $\tau^0$, and that it is bounded below by a positive constant outside every neighbourhood of $\tau^0$. For $\tau^0=+\infty$ this holds because $\lambda$ is positive and decreasing on $\mathbf R$. For $\tau^0\in\mathbf R$, the function $\lambda-\ell_{\tau^0}$ is strictly convex on $\mathbf R$ with vanishing derivative at $\tau^0$, it tends to $+\infty$ at both ends of $\mathbf R$, since $\lambda(\tau)\ge\max\{0,-\tau\}$ and $p<1$, and it equals $+\infty$ at $+\infty$.

\emph{Regularizability.} Let $f\in k[Y]$ and $x\in Y^{\an}$. We claim that $\log(1+|f|)\in\mathcal{R}(Y)$, and that $\ell_{\tau^0}(-\log|f|)\in\mathcal{R}(Y)$ for $\tau^0\coloneqq-\log|f|(x)$. For the first claim, the binomial formula applied to the nonnegative real number $|f|$ gives $P_m\le\log(1+|f|)\le P_m+m^{-1}\log(m+1)$, where
\[
P_m\coloneqq\frac1m\max\left\{0,\ \max_{1\le\alpha\le m}\left(\log\binom m\alpha+\log|f^\alpha|\right)\right\}.
\]
Rounding the constants of $P_m+m^{-1}\log(m+1)$ up to rational numbers, for $m$ large, gives functions $Q_j$ of the form in \cref{lem:regularizable-closure}(3) with $\log(1+|f|)\le Q_j\le\log(1+|f|)+2^{-j-2}$. The functions $Q_j+2^{-j}$ are regularizable and decrease to $\log(1+|f|)\ge0$, so \cref{lem:regularizable-closure}(2) applies. For the second claim we may assume $\tau^0$ finite; then $f\ne0$ and $\ell_{\tau^0}(-\log|f|)=C+p\log|f|$ with $C\coloneqq\lambda(\tau^0)+p\tau^0$. Fix generators $z_1,\ldots,z_n$ of $k[Y]$ and put $b\coloneqq1+\max\{0,\log|z_1|,\ldots,\log|z_n|\}\in\FS(Y)$. Choose rational numbers $a_m\uparrow p$, $a'_m\downarrow p$ and $C_m\downarrow C$ with $a_m>0$. The functions
\[
q_m\coloneqq\max\bigl\{C_m+a_m\log|f|,\ C_m+a'_m\log|f|,\ b/m-m\bigr\}
\]
lie in $\FS(Y)$. They decrease, because $\max\{a\tau,a'\tau\}$ decreases as the interval $[a,a']$ shrinks, whatever the sign of $\tau$, and their limit is $C+p\log|f|$, which is finite at $\vtriv$.

\emph{Construction of $\rho$.} Enumerate the countable set $k[Y]$ as $f_1,f_2,\ldots$, and put $h_j\coloneqq\log(1+|f_j|)$. By the degree estimate of \cref{lem:growth}(1) there are integers $C_j\ge1$ with $0\le h_j\le C_jb$. Put $\gamma_j\coloneqq2^{-j}/C_j$ and $\rho\coloneqq\sum_j\gamma_jh_j$, so that $0\le\rho\le b$. The regularizable functions $R_N\coloneqq\sum_{j\le N}\gamma_jh_j+2^{-N}b$ decrease to $\rho$, because $R_{N+1}-R_N=\gamma_{N+1}h_{N+1}-2^{-N-1}b\le0$; hence $\rho\in\mathcal{R}(Y)$.

\emph{Strictness.} Let $x\in Y^{\an}$ and $V$ an open neighbourhood of $x$. Since the topology of $Y^{\an}$ is the coarsest one making all $|f|$, $f\in k[Y]$, continuous, there are a nonempty finite set $J$ of indices and open neighbourhoods $U_j$ of $\tau_j\coloneqq-\log|f_j|(x)$ in $\mathbf T$ such that $\{y:-\log|f_j|(y)\in U_j\text{ for all }j\in J\}\subseteq V$. Let $\delta_j>0$ be a lower bound of $\lambda-\ell_{\tau_j}$ outside $U_j$, and put
\[
\rho'\coloneqq\rho-\sum_{j\in J}\gamma_j\,(\lambda-\ell_{\tau_j})(-\log|f_j|).
\]
Then $\rho-\rho'$ is nonnegative, vanishes at $x$, and is at least $\delta\coloneqq\min_{j\in J}\gamma_j\delta_j$ outside $V$. Finally, for $N\ge\max J$ the regularizable functions $R_N-\sum_{j\in J}\gamma_j\bigl(h_j-\ell_{\tau_j}(-\log|f_j|)\bigr)$ decrease to $\rho'$, which is finite at $x$; so $\rho'\in\mathcal{R}(Y)$.
\end{proof}

The Runge property needed in \cref{lem:hartogs} will be obtained from hull inequalities on $\widetilde Y=Y\times\mathbf A^1$, as follows. Recall that $p\colon\widetilde Y\to Y$ is the projection and $w$ the coordinate on $\mathbf A^1$.

\begin{lemma}\label{lem:hull-criterion}
Let $u\colon Y^{\an}\to[-\infty,\infty)$ be upper semicontinuous, of logarithmic growth and not identically $-\infty$. Put $\widetilde u\coloneqq u\circ p+\log|w|$, and let $\widetilde\rho$ be a function as in \cref{lem:strict-exhaustion} for $\widetilde Y$. If $\widetilde u+c\widetilde\rho$ satisfies the hull inequality for every rational $c>0$, then $u\in\mathcal{R}(Y)$.
\end{lemma}

\begin{proof}
Let $c>0$ be rational. If $K\subseteq\{\widetilde u+c\widetilde\rho<0\}$ is compact, then $\sup_K(\widetilde u+c\widetilde\rho)<0$ by upper semicontinuity, so $\widehat K\subseteq\{\widetilde u+c\widetilde\rho<0\}$ by the hull inequality. Thus these open sets are Runge. As $c\downarrow0$ they increase to $\{\widetilde u<0\}$, since $\widetilde\rho$ is finite and nonnegative. Hence $\{\widetilde u<0\}$ is Runge, and \cref{lem:hartogs} gives $u\in\mathcal{R}(Y)$.
\end{proof}

The following gluing statement is \cite[Lemma~6.8 and Corollary~6.9]{Wu26}, with the local regularizability hypothesis used in their proofs.\footnote{In Lemma~6.8, assume $u|_{U_\alpha^{\an}}\in\mathcal{R}(U_\alpha)$, as in the statement below; regularizability of $u+c\rho$ alone does not justify the auxiliary function $u+c\rho'+\epsilon v$. One may take $v=\max_i\log|f_i|$ for generators of the ideal of $Z$, and choose $\epsilon$ rational.}

\begin{proposition}\label{prop:gluing}
Let $Z\subsetneq Y$ be closed \textup(possibly empty\textup) and $u$ a weakly psh function on $Y^{\an}$ of logarithmic growth. Assume that $Y\setminus Z$ is covered by affine open subsets $U_\alpha$ with $u|_{U_\alpha^{\an}}\in\mathcal{R}(U_\alpha)$, and let $\rho$ be as in \cref{lem:strict-exhaustion}.
\begin{enumerate}
\item For every rational $c>0$, every compact $K\subseteq Y^{\an}$ and every $x\in\widehat K\setminus Z^{\an}$, we have $(u+c\rho)(x)\le\sup_K(u+c\rho)$.
\item If $Z=\varnothing$, then $u\in\mathcal{R}(Y)$.
\end{enumerate}
\end{proposition}

\begin{proof}
(1) Use the proof of \cite[Lemma~6.8]{Wu26}, with \cref{lem:strict-exhaustion,lem:rossi,lem:principal-hull} as its inputs. In the notation of that proof, the auxiliary function $u+c\rho'+\epsilon v$ is regularizable on the principal neighbourhood $Y_g\subseteq U_\alpha$ by \cref{lem:regularizable-closure}(3),(4), and has logarithmic growth on $Y$ since $\rho'\le\rho$. Thus its strict boundary drop contradicts the principal-open hull inequality exactly as there. The proof uses only these inputs, which hold over $k$.

(2) Apply (1) with $Z=\varnothing$ to $\widetilde Y=Y\times\mathbf A^1$ and $\widetilde u=u\circ p+\log|w|$, as in \cite[Corollary~6.9]{Wu26}. A witness for $u$ pulls back to a witness on $\widetilde Y$, and \cref{lem:regularizable-closure}(3),(4) gives regularizability on each $U_\alpha\times\mathbf A^1$. Then \cref{lem:hull-criterion} applies; $u$ is not identically $-\infty$, being regularizable on some $U_\alpha$.
\end{proof}

\subsection{Norms along generically \'etale projections}\label{subsec:norms}

To extend regularizability across a closed subset, we project a neighbourhood of one of its points finitely onto an open subset of $\mathbf A^n$ and push functions forward by a weighted sum over the fibres. Outside a dense open subset over which the projection is \'etale, the weights are defined by restricting to general lines, and regularity of the pushforward is checked on curves. We use subharmonic functions on analytic curves over an arbitrary complete non-Archimedean field in the sense of \cite[Definition~3.4]{GublerRabinoff25}, which may be identically $-\infty$ on a component. By \cref{lem:local-basic}(4), the pullback of a locally psh function under a morphism from a curve is subharmonic.

\begin{lemma}\label{lem:curves}
Let $F$ be a complete non-Archimedean field.
\begin{enumerate}
\item Let $q\colon C'\to C$ be a finite flat surjective morphism of normal algebraic $F$-curves, with $C$ integral. For $y\in C^{\an}$ write $\mathcal H(y)\otimes_{\mathcal O_C}q_*\mathcal O_{C'}=\prod_xA_x$, the product of local Artinian algebras indexed by the points $x\in C'^{\an}$ over $y$, and put $m_x(q)\coloneqq\dim_{\mathcal H(y)}A_x$. If $u$ is subharmonic on $C'^{\an}$, then $\Norm_qu(y)\coloneqq\sum_{q(x)=y}m_x(q)u(x)$ is subharmonic on $C^{\an}$.
\item Two subharmonic functions on the analytification of an algebraic $F$-curve that agree outside finitely many rigid points are equal.
\item Assume that $F$ is an algebraically closed complete extension of $k$, let $r$ denote the natural maps from $F$-analytic spaces to $k$-analytic spaces, and let $u$ be a weakly psh function on $Y^{\an}$ for a variety $Y$ over $k$. For every morphism $p\colon C\to Y_F$ from a smooth $F$-curve, $u\circ r\circ p^{\an}$ is subharmonic on $C^{\an}$.
\end{enumerate}
\end{lemma}

\begin{proof}
(1) Let $K=F(C)$, let $E_i$ be the function fields of the components $C'_i$ of $C'$, and let $\nu_i$ be the inseparable degree of $E_i/K$. Choose a finite normal extension $L/K$ containing a conjugate of each $E_i$, and let $r_L\colon D\to C$ be the normalization of $C$ in $L$, which is finite. Every $K$-embedding $\sigma\colon E_i\to L$ induces a $C$-morphism $p_{i,\sigma}\colon D\to C'_i$, since $C'_i$ is the normalization of $C$ in $E_i$. We claim that
\begin{equation}\label{eq:curve-trace-splitting}
(\Norm_qu)\circ r_L^{\an}=\sum_i\nu_i\sum_{\sigma\colon E_i\hookrightarrow L}u\circ p_{i,\sigma}^{\an}
\end{equation}
for every function $u$. Work over an affine open $\operatorname{Spec}R\subseteq C$ over which $q_*\mathcal O_{C'}$ is free with basis $b_1,\ldots,b_d$, where $d$ is the degree of $q$, and let $S$ be the integral closure of $R$ in $L$. Since the generic fibre of $q$ decomposes into the fields $E_i$, and the characteristic polynomial of an element of $E_i$ over $K$ is the product of its conjugates, each with multiplicity $\nu_i$, we have
\[
\det\bigl(Z-m_{b(T)}\bigr)=\prod_i\prod_{\sigma\colon E_i\hookrightarrow L}\bigl(Z-\sigma(b(T))\bigr)^{\nu_i}\qquad\text{in }S[T_1,\ldots,T_d,Z],
\]
where $b(T)=\sum_jT_jb_j$ and $m_{b(T)}$ is multiplication by $b(T)$. Fix $z\in D^{\an}$ and put $y=r_L(z)$. Specialize to $\mathcal H(z)$ and substitute for $T$ the coordinates of the idempotent of $A_x$. The left side becomes $Z^{d-m_x(q)}(Z-1)^{m_x(q)}$, and a factor on the right becomes $Z-1$ exactly when $p_{i,\sigma}(z)=x$. Hence $m_x(q)=\sum_{(i,\sigma):\,p_{i,\sigma}(z)=x}\nu_i$, which proves \eqref{eq:curve-trace-splitting}. Its right side is subharmonic by \cite[Proposition~3.5(4),(8)]{GublerRabinoff25}. Since $r_L^{\an}$ is closed and surjective, $\Norm_qu$ is upper semicontinuous, and it is subharmonic by finite surjective descent \cite[Proposition~3.5(8)]{GublerRabinoff25}.

(2) Let $h$ be a subharmonic function and $z$ a rigid point. By upper semicontinuity $h(z)\ge\limsup_{x\to z,\,x\ne z}h(x)$. If the inequality were strict, $h$ would have a strict local maximum at $z$, contradicting the maximum principle \cite[Proposition~3.5(6)]{GublerRabinoff25}: the analytification of an algebraic curve has no boundary and no isolated points. Hence $h(z)$ is determined by the values of $h$ near $z$, also when $h(z)=-\infty$.

(3) Subharmonicity being local \cite[Proposition~3.5(1)]{GublerRabinoff25}, we may assume $C$ connected, hence irreducible. Let $\mu\colon Y'\to Y$ be a witness. Then $u\circ r\circ\mu_F^{\an}$ is locally psh on $Y_F'^{\an}$ by \cref{lem:local-basic}(2). The fibre product $C\times_{Y_F}Y'_F$ is proper and surjective over $C$; choose an irreducible component dominating $C$, a closed point of its generic fibre over $C$, and let $C''$ be the normalization of its closure. Then $C''$ is a smooth curve, finite and surjective over $C$, with a morphism $C''\to Y'_F$ lifting $p$. The pullback of $u\circ r\circ p^{\an}$ to $C''^{\an}$ is the pullback of $u\circ r\circ\mu_F^{\an}$, hence subharmonic, and finite surjective descent \cite[Proposition~3.5(8)]{GublerRabinoff25} shows that $u\circ r\circ p^{\an}$ is subharmonic. This works even if $p(C)$ lies in the locus over which $\mu$ is not an isomorphism.
\end{proof}

The next proposition is \cite[Proposition~6.12]{Wu26}, extended to positive characteristic. The proof in \cite{Wu26} applies a regularization theorem for smooth affine curves over the completed residue field $\mathcal H(s)$ of a point $s$, which need not be discretely or trivially valued; we only use subharmonicity over an algebraically closed complete extension, via \cref{lem:curves}. Instead of \cite[Proposition~6.2(2)]{Wu26}, applied there on the \'etale locus, we assume regularizability on an affine open cover of $Y\setminus Z$ and take the \'etale locus off $\pi(Z)$.

\begin{proposition}\label{prop:nonflat-norm}
Let $Y$ be a normal affine variety of dimension $n\ge1$, $S\subseteq\mathbf A^n$ a principal open subset, $\pi\colon Y\to S$ a finite surjective generically \'etale morphism of degree $d$, and $Z\subsetneq Y$ closed. There are positive integers $m_y(\pi)$, $y\in Y^{\an}$, with $\sum_{\pi(y)=s}m_y(\pi)=d$ for every $s\in S^{\an}$, such that the following holds. If $u$ is a weakly psh function on $Y^{\an}$ of logarithmic growth, and $Y\setminus Z$ is covered by affine open subsets on which $u$ is regularizable, then
\[
\Norm_\pi u(s)\coloneqq\sum_{\pi(y)=s}m_y(\pi)u(y)
\]
defines an element of $\mathcal{R}(S)$.
\end{proposition}

\begin{proof}
\emph{Choice of weights.} Choose $0\ne h\in k[S]$ such that $S_0\coloneqq S_h\subseteq S\setminus\pi(Z)$ and $\pi$ is \'etale over $S_0$, and put $Y_0\coloneqq\pi^{-1}(S_0)$, an affine open subset of $Y\setminus Z$. Over $S_0^{\an}$ use the weights of \cref{lem:finite-flat-norm}. Let $s\in S^{\an}\setminus S_0^{\an}$, let $F$ be the completion of an algebraic closure of $\mathcal H(s)$, and let $s_F\in S_F(F)$ be the rigid point induced by $s$. Choose a line $\Lambda\subseteq\mathbf A^n_F$ through $s_F$ that is not contained in the closure of $(S\setminus S_0)_F$ in $\mathbf A^n_F$; this is possible because the lines through $s_F$ contained in a closed subset of dimension less than $n$ form a proper closed subset of the space of lines through $s_F$. Let $C_s\coloneqq\Lambda\cap S_F$; it is a smooth affine curve, and $C_s\cap(S_0)_F$ is dense in $C_s$. Let $C'_s$ be the normalization of the reduction of $Y_F\times_{S_F}C_s$. Since $k$ is algebraically closed, $Y_F$ is integral of dimension $n$. By Krull's height theorem every irreducible component of $Y_F\times_{S_F}C_s$, which is cut out in $Y_F$ by $n-1$ equations, has dimension at least one, and finiteness over $C_s$ gives dimension at most one; so every component is a curve dominating $C_s$. Over $C_s\cap(S_0)_F$ the fibre product is \'etale over $C_s$, hence reduced and normal, so it is not changed there by reduction and normalization. Thus $C'_s$ is a smooth curve, and $q_s\colon C'_s\to C_s$ is finite, flat (being torsion free over a smooth curve) and of degree $d$. Every point of $Y^{\an}$ over $s$ is the image of a point of $Y_F^{\an}$ over $s_F$, which in turn is the image of a rigid point of $C'^{\an}_s$ over $s_F$. For $y\in\pi^{-1}(s)$, let $m_y(\pi)$ be the sum of the multiplicities $m_{x}(q_s)$ over the points $x\in q_s^{-1}(s_F)$ mapping to $y$. Then $m_y(\pi)\ge1$ and $\sum_ym_y(\pi)=d$. These choices do not depend on $u$.

\emph{Regularity.} Let $u$ be as in the statement. The affine open subset $Y_0$ of $Y\setminus Z$ is covered by the intersections of $Y_0$ with the given affine open subsets, which are affine and on which $u$ is regularizable by \cref{lem:regularizable-closure}(4); so $u|_{Y_0^{\an}}\in\mathcal{R}(Y_0)$ by \cref{prop:gluing}(2), and $v_0\coloneqq\Norm_{\pi|_{Y_0}}u\in\mathcal{R}(S_0)$ by \cref{lem:finite-flat-norm}. By \cref{lem:growth}(2), and since functions in $\FS(S)$ are bounded below, $u\le B'g\circ\pi+B''$ for some $g\in\FS(S)$ and constants $B'\in\mathbf Q_{>0}$, $B''\in\mathbf R$, so $v_0\le dB'g+dB''$ on $S_0^{\an}$, where $dB'g\in\FS(S)$.

\emph{Extension to $S$.} We extend $v_0$ to $S$ using the envelope property of $\mathbf P^n$; this replaces \cite[Theorem~3.31]{Wu26}, whose proof uses a smooth compactification with ample boundary. Write $S=\{h_1\ne0\}$ with $h_1\in k[x_1,\ldots,x_n]$ of degree $e$, regard $\mathbf A^n$ as the chart $\{x_0\ne0\}$ of $P\coloneqq\mathbf P^n$, and let $\sigma\in H^0(P,H)$, $H\coloneqq\mathcal O_P(e+1)$, be the product of $x_0$ and the homogenization of $h_1$, so that $P_\sigma=S$. Every $f\in k[S]$ has the form $\tau/\sigma^N$ with $\tau\in H^0(P,NH)$, so $\log|f|\le-N\log|\sigma|$. Fix rational numbers $c_i\downarrow0$. By \cref{lem:principal-extension}, the extension $w_i$ of $v_0+c_i\log|h|$ by $-\infty$ lies in $\mathcal{R}(S)$, and by the previous remark there are $\kappa\in\mathbf Q_{>0}$ and $B_0\in\mathbf R$, independent of $i$, with $w_i\le-\kappa\log|\sigma|+B_0$. Fix a rational $a>\kappa$. By \cref{lem:compactification}(2), the extension $\phi_i$ of $w_i+a\log|\sigma|$ by $-\infty$ lies in $\PSH(P,aH)$, and $\phi_i\le B_0$. Since $P$ is smooth, $\EP_P(c_1(aH))$ holds by \cref{thm:smooth-EP}, so $\psi_j\coloneqq\sups_{i\ge j}\phi_i$ lies in $\PSH(P,aH)$. On the open set $S_0^{\an}$ we have $\sup_{i\ge j}\phi_i=v_0+a\log|\sigma|+\max\{c_j\log|h|,0\}$, which is upper semicontinuous there; hence $\psi_j$ is given by the same formula on $S_0^{\an}$. The $\psi_j$ decrease to a function $\psi$ equal to $v_0+a\log|\sigma|$ on $S_0^{\an}$, which is finite at $\vtriv\in S_0^{\an}$ by \cref{lem:regularizable-closure}(1); so $\psi$ is $aH$-psh by PSH2. By \cref{lem:compactification}(1), $v\coloneqq(\psi-a\log|\sigma|)|_{S^{\an}}$ lies in $\mathcal{R}(S)$, and $v=v_0$ on $S_0^{\an}$.

\emph{Identification off $S_0$.} It remains to show $v(s)=\Norm_\pi u(s)$ for $s\notin S_0^{\an}$. With the notation above, $v\circ r$ restricts to a locally psh, hence subharmonic, function on $C_s^{\an}$. The pullback of $u\circ r$ to $C'^{\an}_s$ is subharmonic by \cref{lem:curves}(3), so its norm under $q_s$ is subharmonic on $C_s^{\an}$ by \cref{lem:curves}(1). Over $C_s\cap(S_0)_F$, the map $C'_s\to Y_F\times_{S_F}C_s$ is an isomorphism, so the fibre algebras of $q_s$ and of $\pi_F$ agree there; by the base change compatibility in \cref{lem:finite-flat-norm}, the two subharmonic functions agree on $(C_s\cap(S_0)_F)^{\an}$. They differ at most at the finitely many rigid points of $C_s\setminus(S_0)_F$, hence are equal by \cref{lem:curves}(2). Evaluating at $s_F$ gives $v(s)=\sum_ym_y(\pi)u(y)$.
\end{proof}

\subsection{Extension across a closed subset}\label{subsec:singular-locus}

We now show that regularizability extends across proper closed subsets. We first recall how to project a neighbourhood of a smooth point of a closed subset finitely onto an open subset of $\mathbf A^n$. The following is \cite[Lemma~6.16]{Wu26}, supplemented by generic \'etaleness, which is needed in positive characteristic.

\begin{lemma}\label{lem:projection-existence}
Let $Y$ be an affine variety of dimension $n\ge1$, $Z\subsetneq Y$ a closed subset, and $z\in Z$ a closed point at which $Z$ is smooth. There are an affine open neighbourhood $U$ of $z$, a principal open subset $S\subseteq\mathbf A^n$ and a finite surjective generically \'etale morphism $\pi\colon U\to S$ whose restriction to $Z\cap U$ is a closed immersion.
\end{lemma}

\begin{proof}
Use the successive general projections in the proof of \cite[Lemma~6.16]{Wu26}, and shrink the target to a principal affine neighbourhood of the image of $z$.\footnote{At the end of that proof, take $U=r^{-1}(S)$ rather than $r^{-1}(W)$. The projection argument only uses smoothness of $Z$ at $z$: apply \cite[\href{https://stacks.math.columbia.edu/tag/0B1Q}{Tag 0B1Q}]{Stacks} to the component through $z$, and avoid the cones joining $z$ to the other components and to the boundary of the compactification.} To ensure generic \'etaleness in positive characteristic, fix also a smooth closed point $w$ of the projective closure of $Y$. At each stage choose the projection centre so that the projection is a local isomorphism at the image of $w$, or unramified there at the last stage. This excludes a further proper closed subset of centres by \cite[\href{https://stacks.math.columbia.edu/tag/0B1Q}{Tag 0B1Q}]{Stacks}. The resulting finite map is unramified at $w$, hence generically \'etale.
\end{proof}

The main step is the following variant of \cite[Lemma~6.15]{Wu26}. We follow the strategy there: at a point where the hull inequality fails, a norm along the projection produces a regularizable function on $S$ that violates it. The barrier is built from a partition $f_1+f_2=1$ separating $Z$ from the rest of its preimage, which avoids the polyhedral neighbourhoods of \cite[Lemma~6.15, Step~2]{Wu26} and the constancy of a local degree asserted in \cite[Lemma~6.15, Step~3]{Wu26}, which is not justified as stated.

\begin{proposition}\label{prop:projection-step}
Let $Y$ be a normal affine variety of dimension $n\ge1$, $S\subseteq\mathbf A^n$ a principal open subset, $\pi\colon Y\to S$ a finite surjective generically \'etale morphism, and $Z\subsetneq Y$ a closed subset such that $\pi|_Z$ is a closed immersion and $Z'\coloneqq\pi^{-1}(\pi(Z))\setminus Z$ is closed in $Y$. Let $u$ be a weakly psh function on $Y^{\an}$ of logarithmic growth, and assume that $Y\setminus Z$ is covered by affine open subsets on which $u$ is regularizable. Then $u\in\mathcal{R}(Y)$.
\end{proposition}

\begin{proof}
Let $\rho$ be as in \cref{lem:strict-exhaustion}, let $d$ be the degree of $\pi$, and let $m_y\coloneqq m_y(\pi)$ be the weights of \cref{prop:nonflat-norm} for $\pi$ and $Z$. We show that $u_c\coloneqq u+c\rho$ satisfies the hull inequality for every rational $c>0$. Since $\widetilde Y=Y\times\mathbf A^1$, $Z\times\mathbf A^1$, $\pi\times\mathrm{id}$ and $u\circ p+\log|w|$ satisfy the same hypotheses (as in the proof of \cref{prop:gluing}(2)), \cref{lem:hull-criterion} then gives $u\in\mathcal{R}(Y)$; note that $u$ is not identically $-\infty$, being regularizable on a nonempty open subset.

Suppose that $u_c(x_0)>\sup_Lu_c$ for some compact $L\subseteq Y^{\an}$ and $x_0\in\widehat L$. The upper semicontinuous function $u_c$ attains its maximum over the compact set $\widehat L$ at some point $x$. Subtracting the finite constant $u_c(x)$ from $u$, we may assume that $u_c(x)=0=\max_{\widehat L}u_c$. Choose a rational $\delta>0$ with $\sup_Lu_c\le-\delta$; in particular $x\notin L$. By \cref{prop:gluing}(1), $x\in Z^{\an}$.

Since $Z$ and $Z'$ are disjoint closed subsets of the affine variety $Y$, there are $f_1,f_2\in k[Y]$ with $f_1+f_2=1$, $f_1=0$ on $Z'$ and $f_2=0$ on $Z$; then $f_1\ne0$, as $Z\ne\varnothing$. By the ultrametric inequality, the open sets $O_1\coloneqq\{|f_2|<1\}\supseteq Z^{\an}$ and $O_2\coloneqq\{|f_1|<1\}\supseteq Z'^{\an}$ are disjoint, $|f_1|=1$ on $O_1$ and $|f_2|=1$ on $O_2$. The map $\pi^{\an}$ is proper, hence closed, and it is injective on $Z^{\an}$ since $\pi|_Z$ is a closed immersion. Hence the fibre $\pi^{-1}(\pi(x))$ consists of $x$ and of a compact subset $\Gamma$ of $Z'^{\an}$. We choose, in this order:
\begin{enumerate}
\item a relatively compact open neighbourhood $A$ of $\Gamma$ with $\overline A\subseteq\{\log|f_1|<-1\}$, and a rational $M\ge0$ with $u_c\le M$ on $\overline A$;
\item an open neighbourhood $V_x\subseteq\{\log|f_2|<-1\}$ of $x$ on which $u_c<\epsilon\coloneqq\delta/(2d)$;
\item the open neighbourhood $V\coloneqq S^{\an}\setminus\pi(Y^{\an}\setminus(V_x\cup A))$ of $\pi(x)$, which satisfies $\pi^{-1}(V)\subseteq V_x\cup A$;
\item a relatively compact open neighbourhood $B$ of $x$ with $\overline B\subseteq V_x\cap\pi^{-1}(V)$ and $\overline B\cap L=\varnothing$;
\item $\rho'\in\mathcal{R}(Y)$ and $\delta'>0$ as in \cref{lem:strict-exhaustion} for $x$ and $B$.
\end{enumerate}
Put $u'\coloneqq u+c\rho'\le u_c$ and $\eta\coloneqq\max\{\delta\log|f_2|,-\delta\}\in\mathcal{R}(Y)$, and
\[
u^\sharp\coloneqq\max\bigl\{u'+(M+1)\log|f_1|,\ \eta\bigr\}.
\]
This function is weakly psh, with the witness of $u$, and of logarithmic growth, and it is regularizable on every affine open subset on which $u$ is, by \cref{lem:regularizable-closure}(3),(4). On $V_x\subseteq O_1$ we have $\log|f_1|=0$ and $\eta=-\delta$, so $u^\sharp=\max\{u',-\delta\}<\epsilon$. On $A\subseteq O_2$ we have $\eta=0$ and $u'+(M+1)\log|f_1|<M-(M+1)$, so $u^\sharp=0$. By \cref{prop:nonflat-norm}, $h\coloneqq\Norm_\pi u^\sharp\in\mathcal{R}(S)$. Since $u^\sharp=0$ on $\Gamma\subseteq A$ and $u'(x)=u_c(x)=0$, we get $h(\pi(x))=m_x u^\sharp(x)=0$.

By \cref{lem:rossi}, $x$ lies in the hull of the compact set $K\coloneqq\partial B\cap\widehat L$; hence $\pi(x)$ lies in the hull of $\pi(K)$ in $S$, and the hull inequality for $h$ gives $0\le\sup_Kh\circ\pi$. Since $h\circ\pi$ is upper semicontinuous, it suffices to show $h(\pi(y))<0$ for every $y\in K$ to reach a contradiction. We have $y\in V_x$, and $\pi^{-1}(\pi(y))\subseteq V_x\cup A$ since $\pi(y)\in V$. As $y\in\widehat L\setminus B$, we have $u'(y)\le u_c(y)-c\delta'\le-c\delta'$. If $y\notin Z^{\an}$, then $u_c(y)\le\sup_Lu_c\le-\delta$ by \cref{prop:gluing}(1), so $u^\sharp(y)=-\delta$; since $u^\sharp\le\epsilon$ at the other points of the fibre and $m_y\ge1$,
\[
h(\pi(y))\le-\delta\,m_y+(d-m_y)\,\epsilon<-\delta/2 .
\]
If $y\in Z^{\an}$, the other points of $\pi^{-1}(\pi(y))$ lie in $Z'^{\an}$, which does not meet $V_x\subseteq O_1$; so they lie in $A$, where $u^\sharp=0$, and $h(\pi(y))=m_y\max\{u'(y),-\delta\}\le m_y\max\{-c\delta',-\delta\}<0$.
\end{proof}

The inductions in \cite[Proposition~6.14 and Theorem~6.3]{Wu26}, the latter starting from the smooth locus, are replaced by the following Noetherian argument.

\begin{theorem}\label{thm:normal-local-regularization}
Let $Y$ be a normal affine variety over the countable field $k$, $Z\subsetneq Y$ a closed subset, and $u$ a weakly psh function on $Y^{\an}$ of logarithmic growth. If $Y\setminus Z$ is covered by affine open subsets on which $u$ is regularizable, then $u\in\mathcal{R}(Y)$.
\end{theorem}

\begin{proof}
Let $\Sigma\subseteq Y$ be the union of the affine open subsets $U$ with $u|_{U^{\an}}\in\mathcal{R}(U)$, and $W\coloneqq Y\setminus\Sigma$, a closed subset of $Z$. By \cref{prop:gluing}(2) it suffices to show $W=\varnothing$. Suppose not; then $n\coloneqq\dim Y\ge1$. Since $W$ is reduced and $k$ is perfect, $W$ is smooth at some closed point $z$, and removing from $Y$ the singular locus of $W$ and the irreducible components of $W$ not containing $z$, we find a principal open neighbourhood $Y_1$ of $z$ such that $W\cap Y_1$ is smooth and irreducible. \Cref{lem:projection-existence} gives an affine open neighbourhood $U\subseteq Y_1$ of $z$ and a finite surjective generically \'etale morphism $\pi\colon U\to S$ onto a principal open subset of $\mathbf A^n$ that restricts to a closed immersion on the irreducible set $W_U\coloneqq W\cap U$. Put $Z'\coloneqq\pi^{-1}(\pi(W_U))\setminus W_U$, and let $\overline{Z'}$ be its closure in $U$, so that $\overline{Z'}=Z'\cup E$ with $E\coloneqq\overline{Z'}\cap W_U$. Every irreducible component of $\overline{Z'}$ meets $Z'$, so it is not contained in $W_U$, and it has dimension at most $\dim\pi^{-1}(\pi(W_U))=\dim W_U$; hence it does not contain $W_U$. Thus $E$ is a proper closed subset of $W_U$, and $\pi(E)$ is a proper closed subset of the irreducible set $\pi(W_U)$. Choose $g\in k[S]$ vanishing on $\pi(E)$ but not on $\pi(W_U)$, and put $S_1\coloneqq S_g$ and $U_1\coloneqq\pi^{-1}(S_1)$. Then $\pi$ restricts to a finite surjective generically \'etale morphism $U_1\to S_1$ onto a principal open subset of $\mathbf A^n$, the set $W_1\coloneqq W\cap U_1$ is nonempty, $\pi|_{W_1}$ is a closed immersion, and $\pi^{-1}(\pi(W_1))\setminus W_1=Z'\cap U_1=\overline{Z'}\cap U_1$ is closed in $U_1$. Moreover $U_1\setminus W_1\subseteq\Sigma$ is covered by affine open subsets on which $u$ is regularizable, by \cref{lem:regularizable-closure}(4). By \cref{prop:projection-step}, $u|_{U_1^{\an}}\in\mathcal{R}(U_1)$, so $U_1\subseteq\Sigma$, contradicting $W_1\ne\varnothing$.
\end{proof}

\begin{corollary}\label{cor:normal-weak-global}
Let $X$ be a normal projective variety over the countable field $k$ and $L\in\Pic(X)_{\mathbf Q}$ ample. Let $\varphi\colon X^{\an}\to[-\infty,\infty)$ be a function, and suppose that there is a birational morphism $\mu\colon X'\to X$ from an integral projective variety $X'$ with $\varphi\circ\mu^{\an}\in\PSH(X',\mu^*L)$. Then $\varphi\in\PSH(X,L)$.
\end{corollary}

\begin{proof}
This corresponds to \cite[Theorem~6.5]{Wu26}. Replacing $X'$ by its normalization preserves the hypothesis by pullback (\cref{prop:BJ-package}(4)); we may therefore assume that $X'$ is normal. Since $\mu^{\an}$ is proper and surjective, $\varphi$ is upper semicontinuous, hence bounded above. After multiplying $L$ and $\varphi$ by the same positive integer, we may assume that $L$ is a very ample line bundle, with sections $z_0,\ldots,z_N$ defining an embedding. Put $U_i\coloneqq X_{z_i}$. Then $-\log|z_i|=\max_j\log|z_j/z_i|$ on $U_i^{\an}$, and this function lies in $\FS(U_i)$. By \cref{lem:global-to-local}(2), applied on $X'$ to the semiample line bundle $\mu^*L$, the function $u_i\coloneqq\varphi-\log|z_i|$ is weakly psh on $U_i^{\an}$ with witness $\mu$, and it has logarithmic growth.

Let $Z\subsetneq X$ be a closed subset such that $\mu$ is an isomorphism over $X\setminus Z$. We show that every closed point $y\in U_i\setminus Z$ has an affine open neighbourhood $W\subseteq U_i\setminus Z$ with $u_i|_{W^{\an}}\in\mathcal{R}(W)$. Let $y'\in X'$ be the point over $y$. Choose an effective Cartier divisor $D$ on $X'$ such that $\mathcal O_{X'}(D)$ is ample and $\operatorname{Supp}D$ contains $\mu^{-1}(Z\cup V(z_i))$ but not $y'$; for instance the divisor of a section, not vanishing at $y'$, of the ideal of $\mu^{-1}(Z\cup V(z_i))$ twisted by a large multiple of an ample line bundle. Then $W'\coloneqq X'\setminus\operatorname{Supp}D$ is an affine neighbourhood of $y'$, which $\mu$ maps isomorphically onto an affine open neighbourhood $W\subseteq U_i\setminus Z$ of $y$. Let $s_D$ be the canonical section of $\mathcal O_{X'}(D)$, let $H\coloneqq\mu^*L\otimes\mathcal O_{X'}(D)$, which is ample, and $\sigma\coloneqq\mu^*z_i\otimes s_D\in H^0(X',H)$, so that $X'_\sigma=W'$. By \cref{prop:BJ-package}(3),(5), $\Phi\coloneqq\varphi\circ\mu^{\an}+\log|s_D|$ lies in $\PSH(X',H)$, and $\Phi-\log|\sigma|=u_i\circ\mu^{\an}$ on $W'^{\an}$. By \cref{lem:compactification}(1), $u_i\circ\mu^{\an}|_{W'^{\an}}\in\mathcal{R}(W')$, that is, $u_i|_{W^{\an}}\in\mathcal{R}(W)$. These sets $W$ cover $U_i\setminus Z$, since they contain all its closed points.

By \cref{thm:normal-local-regularization}, $u_i\in\mathcal{R}(U_i)$. For rational $a>1$, \cref{lem:compactification}(2) shows that the function $\psi_{i,a}$ equal to $\varphi+(a-1)\log|z_i|$ on $U_i^{\an}$ and to $-\infty$ elsewhere lies in $\PSH(X,aL)$. Since $\max_i\log|z_i|=0$ and the maximum is attained at an index $i$ with the given point in $U_i^{\an}$, we get $\max_i\psi_{i,a}=\varphi$. Thus $\varphi\in\PSH(X,aL)$ for all rational $a>1$, and PSH2, applied to this constant family as $a\downarrow1$, gives $\varphi\in\PSH(X,L)$; note that $\varphi$ is finite at $\vtriv$ because $\varphi\circ\mu^{\an}$ is psh.
\end{proof}

\section{Ample envelopes on normal varieties}\label{sec:normal-ample}

In this section $k$ is again an arbitrary algebraically closed trivially valued field. We prove the envelope property for ample classes on normal varieties. Over a countable field, it follows from the smooth case by descent along the Stein factorization of an alteration: through the birational part by \cref{cor:normal-weak-global}, and through the finite part by \cref{thm:finite-normal-descent}. The general case is reduced to the countable one by extension of scalars.

\begin{theorem}\label{thm:normal-ample-EP}
Let $X$ be a normal projective variety over $k$. Then $\EP_X(\theta)$ holds for every ample class $\theta\in N^1(X)_{\mathbf R}$.
\end{theorem}

\begin{proof}
\emph{Step 1: $k$ countable and $\theta=c_1(L)$ with $L\in\Pic(X)_{\mathbf Q}$ ample.} Choose an alteration $Y\to X$ with $Y$ smooth and projective \cite[Theorem~4.1]{deJong}, and let
\[
Y\xrightarrow{\ \beta\ }T\xrightarrow{\ h\ }X
\]
be its Stein factorization \cite[\href{https://stacks.math.columbia.edu/tag/03H0}{Tag 03H0}]{Stacks}. Thus $h$ is finite, $\beta$ has geometrically connected fibres, and $T$ is the normalization of $X$ in the function field $k(Y)$. Since $Y$ is normal, so is $T$, and since $Y\to X$ is generically finite, $k(T)=k(Y)$; hence $\beta$ is birational. The variety $T$ and the morphism $\beta$ are projective, and $h^*L$ is ample.

We prove $\EP_T(c_1(h^*L))$ using \cref{prop:EP-criterion}. Let $(\varphi_i)$ be a bounded-above increasing net in $\PSH(T,h^*L)$. The class $\beta^*h^*L$ is nef, so \cref{thm:smooth-EP} and \cref{prop:EP-criterion} give $\Phi\in\PSH(Y,\beta^*h^*L)$ with
\[
\Phi(w)=\sup_i\varphi_i(\beta(w))\qquad(w\in Y^{\divv}).
\]
We claim that $\Phi$ is constant on every fibre of $\beta^{\an}$. Fix $t\in T^{\an}$ and an admissible pair $(m,s)$ for $h^*L$ with $t\in T_s^{\an}$. By \cref{lem:global-to-local}(2), applied to the semiample line bundle $\beta^*h^*L$, the function $\Phi-m^{-1}\log|\beta^*s|$ is locally psh on $\beta^{-1}(T_s)^{\an}$. The fibre $(\beta^{\an})^{-1}(t)$ is the analytification of the proper $\mathcal H(t)$-scheme $Y\times_T\operatorname{Spec}\mathcal H(t)$, which is connected since $\beta$ has geometrically connected fibres; hence it is compact, connected and without boundary \cite[\S3.4--3.5]{Berkovich90}. By \cref{lem:local-basic}(2),(4), the restriction of $\Phi-m^{-1}\log|\beta^*s|$ to it is classically psh, and $m^{-1}\log|\beta^*s|$ is constant on it. Unless $\Phi$ is identically $-\infty$ on the fibre, it attains a finite maximum there, the fibre being compact and $\Phi$ upper semicontinuous; the set where the maximum is attained is closed by upper semicontinuity and open by the local maximum principle \cite[Theorem~6.1]{GublerRabinoff25}, hence equal to the fibre. This proves the claim.

Since $\beta^{\an}$ is surjective, $\Phi=\varphi\circ\beta^{\an}$ for a function $\varphi$ on $T^{\an}$, and \cref{cor:normal-weak-global} gives $\varphi\in\PSH(T,h^*L)$. Since $\beta^{\an}$ maps $Y^{\divv}$ bijectively onto $T^{\divv}$, we have $\varphi=\sup_i\varphi_i$ on $T^{\divv}$, which proves $\EP_T(c_1(h^*L))$. \Cref{thm:finite-normal-descent} now gives $\EP_X(c_1(L))$.

\emph{Step 2: $\theta=c_1(L)$ with $L\in\Pic(X)_{\mathbf Q}$ ample.} By \cref{prop:ample-envelopes}(1), it suffices to show that $\mathrm P_L(f)$ is $L$-psh for every rational PL function $f$. After replacing $L$ and $f$ by a common multiple, $f=\phi_{\mathcal L}$ for a test configuration $(\mathcal X,\mathcal L)$ for $(X,L)$ with $\mathcal L$ a line bundle. Choose a finitely generated subfield of $k$ over which $X$, $L$ and $(\mathcal X,\mathcal L)$ are defined, and let $k_0$ be its algebraic closure in $k$; it is countable. We obtain $X_0$, $L_0$ and $(\mathcal X_0,\mathcal L_0)$ over $k_0$ whose scalar extensions are $X$, $L$ and $(\mathcal X,\mathcal L)$. Let $r\colon X^{\an}\to X_0^{\an}$ be the natural map; then $f=f_0\circ r$ with $f_0\coloneqq\phi_{\mathcal L_0}$, since Gauss extension commutes with restriction to $k_0$. Normality and ampleness descend along the faithfully flat extension $k_0\subseteq k$, so Step~1 and \cref{prop:ample-envelopes}(1) show that $\mathrm P_{L_0}(f_0)$ is $L_0$-psh. Let $\mathfrak a_{m,0}$ and $\mathfrak a_m$ be the base ideals of $m\mathcal L_0$ and $m\mathcal L$. Flat base change \cite[\href{https://stacks.math.columbia.edu/tag/02KH}{Tag 02KH}]{Stacks}, together with the exactness of flat pullback applied to the evaluation maps, gives $\mathfrak a_m=\mathfrak a_{m,0}\mathcal O_{\mathcal X}$, hence $\phi_{\mathfrak a_m}=\phi_{\mathfrak a_{m,0}}\circ r$, and \cref{prop:ample-envelopes}(2), applied over $k_0$ and over $k$, gives $\mathrm P_L(f)=\mathrm P_{L_0}(f_0)\circ r$. Since $L_0$ is ample, $\mathrm P_{L_0}(f_0)$ is the limit of a decreasing sequence in $\mathcal{H}^{\mathrm{gf}}_{\mathbf{Q}}(X_0,L_0)$ \cite[Corollary~12.18(iii) and Example~3.9]{BJ22}. Its pullbacks under $r$ lie in $\mathcal{H}^{\mathrm{gf}}_{\mathbf{Q}}(X,L)$ and decrease to the bounded function $\mathrm P_L(f)$, which is therefore $L$-psh directly by \cite[Definition~4.1]{BJ22}.

\emph{Step 3: $\theta$ real.} Choose rational ample classes $\theta_j\to\theta$ with $\theta_j-\theta$ ample, which is possible by openness of the ample cone and density of rational classes. Each $\theta_j$ has the envelope property by Step~2, and \cite[Lemma~5.9]{BJ22} concludes.
\end{proof}

\section{Homogeneous division and arbitrary classes}
\label{sec:siu-envelopes}

In this section $X$ is a normal projective variety over $k$. We pass from ample classes to arbitrary classes in two steps. First, the ample envelope property of \cref{thm:normal-ample-EP} yields homogeneous Siu division (\cref{thm:homogeneous-Siu}): for an effective rational Cartier divisor $E$, a homogeneous psh function bounded by $e_E$ (defined below) is $e_E$ plus a homogeneous psh function. Division transfers increasing limits of homogeneous psh functions from an ample class to an arbitrary class. Second, a psh function is recovered from its homogeneous Legendre slices, and this reduces the envelope property to increasing limits of homogeneous functions (\cref{prop:homogeneous-slices}).

\subsection{Division by an effective divisor}

Keep $X$ normal and projective over $k$. For an effective rational Cartier divisor $E$, let $s_E$ be the canonical section of $\mathcal O_X(mE)$ for $m$ sufficiently divisible, and put $e_E\coloneqq m^{-1}\log|s_E|$. This is a homogeneous Fubini--Study function for $\mathcal O_X(E)$, hence a homogeneous psh function of class $[E]$, and $e_E(v)=-v(E)$. We first divide by $tE$ for a fixed coefficient $t<1$.

\begin{lemma}\label{lem:fixed-coefficient-division}
Let $\eta\in N^1(X)_{\mathbf R}$ and let $E\geq0$ be a rational Cartier divisor. If $\psi\in\PSHhom(X,\eta)$ satisfies $\psi\leq e_E$, then, for every rational $0<t<1$, there is $q_t\in\PSHhom(X,\eta-t[E])$ such that $\psi=te_E+q_t$.
\end{lemma}

\begin{proof}
If $E=0$, take $q_t=\psi$. Otherwise write $E=\sum_\ell a_\ell E_\ell$ with prime divisors $E_\ell$ and $a_\ell>0$. Choose a homogeneous Fubini--Study regularization $\psi_i\downarrow\psi$ with $\psi_i\in\mathcal{H}_{\mathrm{hom}}(X,L_i)$ and $c_1(L_i)\to\eta$. Since $\psi_i(\ord_{E_\ell})\to\psi(\ord_{E_\ell})\le-a_\ell<-ta_\ell$, a common tail satisfies $\psi_i(\ord_{E_\ell})\le-ta_\ell$ for all $\ell$. Write $\psi_i=m^{-1}\max_j\log|s_j|$ with $m$ so divisible that $mtE$ is Cartier. Every $s_j$ vanishes to order at least $mta_\ell$ along each $E_\ell$, so, $X$ being normal, it is divisible by the canonical section of $mtE$. This gives $q_i\in\mathcal{H}_{\mathrm{hom}}(X,L_i-tE)$ with $\psi_i=te_E+q_i$. The $q_i$ decrease on $X^{\divv}$, where $e_E$ is finite, hence on $X^{\an}$ by divisorial comparison, and $q_i(\vtriv)=0$. Since the $q_i$ are themselves homogeneous Fubini--Study functions and their classes tend to $\eta-t[E]$, \cite[Definition~4.1]{BJ22} gives $q_t\coloneqq\lim_iq_i\in\PSHhom(X,\eta-t[E])$, and $\psi=te_E+q_t$ in the limit.
\end{proof}

\begin{theorem}\label{thm:homogeneous-Siu}
Let $\theta\in N^1(X)_{\mathbf R}$ and let $E\geq0$ be a rational Cartier divisor. If $\psi\in\PSHhom(X,\theta+[E])$ satisfies $\psi\leq e_E$, then $\psi=e_E+\varphi$ for some $\varphi\in\PSHhom(X,\theta)$.
\end{theorem}

\begin{proof}
Put $u\coloneqq\psi-e_E$ on $X^{\divv}$, where both terms are finite; thus $u\le0$. Choose a very ample Cartier divisor $A$ with $\theta+[A]$ ample, and let $A_1,\ldots,A_r$ be the divisors of finitely many sections generating $\mathcal O_X(A)$. At every center one generator is a unit, so
\begin{equation}\label{eq:auxiliary-divisor-maximum}
\max_je_{A_j}=0\qquad\text{on }X^{\an}.
\end{equation}
Fix $j$ and a rational $0<\epsilon<1$. Since $\psi+e_{A_j}\in\PSHhom(X,\theta+[A]+[E])$ is at most $e_E$, \cref{lem:fixed-coefficient-division} gives, for rational $0<t<1$, functions in $\PSHhom(X,\theta+[A]+(1-t)[E])$ with divisorial values $u+e_{A_j}+(1-t)e_E$. As $t\uparrow1$ these increase to $u+e_{A_j}$, and for $t$ close to $1$ the class $\epsilon[A]-(1-t)[E]$ is ample, so they lie in the ample class $\theta+(1+\epsilon)[A]$. By \cref{thm:normal-ample-EP,prop:EP-criterion}, there is $\chi_j\in\PSHhom(X,\theta+(1+\epsilon)[A])$ with divisorial values $u+e_{A_j}$; thus $\chi_j\le e_{A_j}$. \Cref{lem:fixed-coefficient-division}, applied to $\chi_j$ and $A_j$ with coefficient $1-\epsilon$, gives $q_j\in\PSHhom(X,\theta+2\epsilon[A])$ with divisorial values $u+\epsilon e_{A_j}$. By \eqref{eq:auxiliary-divisor-maximum}, $\varphi\coloneqq\max_jq_j\in\PSHhom(X,\theta+2\epsilon[A])$ has divisorial values $u$. By divisorial comparison $\varphi$ does not depend on $\epsilon$, and $\varphi(\vtriv)=0$, so \cref{lem:homogeneous-decreasing-stability}, applied to this constant net as $\epsilon\downarrow0$, gives $\varphi\in\PSHhom(X,\theta)$. Finally $e_E+\varphi$ and $\psi$ are psh and agree on $X^{\divv}$, hence everywhere.
\end{proof}

\begin{lemma}\label{lem:decreasing-finite-maxima}
Let $X$ be an integral projective variety and $\theta\in N^1(X)_{\mathbf R}$. Fix $\Psi\in\PSHhom(X,\theta)$. For each $n\ge1$, let $r_n$ be a positive integer and choose functions $q_{n,a}\in\PSHhom(X,\theta)$ and constants $c_{n,a}\in\mathbf R$, $1\le a\le r_n$, with $q_{n,a}\le\Psi$. Set
\[
 U_n\coloneqq\max\Bigl\{\Psi-n,\ \max_{1\le a\le r_n}(q_{n,a}+c_{n,a})\Bigr\}.
\]
If $U_n$ decreases pointwise to $u\not\equiv-\infty$, then $u\in\PSH(X,\theta)$.
\end{lemma}

\begin{proof}
\emph{Step 1: Fubini--Study approximants.} Put $q_{n,0}\coloneqq\Psi$ and $c_{n,0}\coloneqq-n$. By \cref{prop:BJ-package}(3), choose homogeneous Fubini--Study approximants $p_{n,a,j}\downarrow q_{n,a}$ for $0\le a\le r_n$. As in the proof of \cref{prop:BJ-package}(5), promote these finitely many approximants to a common rational line bundle $B_{n,j}$, with $c_1(B_{n,j})\to\theta$ as $j$ increases. Define
\[
 F_{n,j}\coloneqq\max_{0\le a\le r_n}(p_{n,a,j}+c_{n,a}),
 \qquad Q_{n,j}\coloneqq\max_{0\le a\le r_n}p_{n,a,j}.
\]
Then $F_{n,j}$ is Fubini--Study with real weights, $Q_{n,j}$ is homogeneous Fubini--Study, and
\[
 F_{n,j}\downarrow U_n,\qquad Q_{n,j}\ge p_{n,0,j}\ge\Psi.
\]
For rational $0<b<1$ and $\epsilon>0$, consider the functions
\[
 G(v)\coloneqq F_{n,j}(bv)+\epsilon.
\]
Each is Fubini--Study in $bB_{n,j}$ and lies above $U_n$, since the $p_{n,a,j}$ are homogeneous and nonpositive.

\emph{Step 2: simultaneous comparison.} Given two such functions $G_k(v)=F_{n_k,j_k}(b_kv)+\epsilon_k$, $k=1,2$, we construct another one below both, with arbitrarily large $n$ and class arbitrarily close to $\theta$. Choose $n\ge\max\{n_1,n_2\}$, rational $\max\{b_1,b_2\}<b<1$, and $0<\epsilon<\min\{\epsilon_1,\epsilon_2\}$. We may take $b$ arbitrarily close to $1$ and $\epsilon$ arbitrarily small. Write $Q_k\coloneqq Q_{n_k,j_k}$ and choose rational numbers $b_k<\beta_k<b$. Since $q_{n,a}\le\Psi\le Q_k$, Step~1 of \cref{lem:homogeneous-decreasing-stability} gives, for all sufficiently large $j$,
\[
 b\,p_{n,a,j}\le\beta_kQ_k
 \qquad(0\le a\le r_n,\ k=1,2).
\]
The constants in the finite maxima are fixed, so there is $C>0$, independent of $j$, such that
\[
 G\le\beta_kQ_k+C,\qquad G_k\ge b_kQ_k-C.
\]
Since $\beta_k>b_k$, these bounds give $G\le G_k$ on $\{Q_k<-R\}$ for a sufficiently large $R$, independent of $j$. This also holds where $Q_k=-\infty$, because the domination inequalities make every $p_{n,a,j}$ equal to $-\infty$ there.

On the compact set $K_k\coloneqq\{Q_k\ge-R\}$, the function $G_k$ is finite and continuous. Since $U_n$ decreases in $n$ and is nonincreasing along each scaling orbit,
\[
 U_n(bv)\le U_{n_k}(bv)\le U_{n_k}(b_kv)
 \le F_{n_k,j_k}(b_kv).
\]
Thus $G$ decreases, as $j$ increases, to a function strictly below $G_k$ on $K_k$, because $\epsilon<\epsilon_k$. Upper semicontinuity and compactness give $G<G_k$ throughout $K_k$ for all sufficiently large $j$. A common such $j$ proves $G\le G_1,G_2$ everywhere; increasing $j$ further makes $c_1(B_{n,j})$ arbitrarily close to $\theta$.

\emph{Step 3: passage to the limit.} Fix a norm on $N^1(X)_{\mathbf R}$ and order the choices $(n,j,b,\epsilon)$ by increasing $n$, decreasing $G$, and decreasing
\[
 e\coloneqq\|c_1(B_{n,j})-\theta\|+1-b+\epsilon.
\]
Step~2 shows that this is a directed preorder and that $e\to0$. The resulting Fubini--Study net decreases, and its classes $b\,c_1(B_{n,j})$ tend to $\theta$. Its infimum is $u$. Indeed, given $v\in X^{\an}$ and a real number $A>u(v)$, first choose $n$ with $U_n(v)<A$, then $j$ with $F_{n,j}(v)<A$, and finally $b$ close to $1$ and $\epsilon$ small enough that $G(v)<A$. The last choice is possible because $p_{n,a,j}(bv)=b\,p_{n,a,j}(v)$, so $F_{n,j}(bv)\to F_{n,j}(v)$ as $b\uparrow1$, also at a value $-\infty$. Hence \cite[Definition~4.1]{BJ22} gives $u\in\PSH(X,\theta)$.
\end{proof}

\subsection{Envelopes from homogeneous slices}

Let $X$ be an integral projective variety over $k$, let $\theta\in N^1(X)_{\mathbf R}$, and let $u\in\PSH(X,\theta)$. Put $b\coloneqq\sup u=u(\vtriv)$. For $\lambda\in\mathbf R$, the \emph{homogeneous Legendre slice} of $u$ at $\lambda$ is
\[
  \widehat u_\lambda\coloneqq\inf_{t>0}\,\{t\mathbin\cdot u-t\lambda\},
\]
with the scaling action \eqref{eq:scaling-action}. The slices $\widehat u_\lambda$, $\lambda<b$, lie in $\PSHhom(X,\theta)$ without an appeal to general decreasing stability. For rational $\lambda$, apply the homogeneous projection described after \cref{prop:BJ-package} to the rational Fubini--Study approximants of $u-\lambda$. Their projections are homogeneous Fubini--Study in the same respective classes and decrease to $\widehat u_\lambda$, since the two infima commute; they vanish at $\vtriv$ because $b>\lambda$. For real $\lambda<b$, take rational parameters increasing to $\lambda$ and apply \cref{lem:homogeneous-decreasing-stability}. The slices decrease in $\lambda$. The one-variable convex duality calculation in \cite[Theorem~6.24]{BJ22} gives
\begin{equation}\label{eq:Legendre-reconstruction}
  u=\sup_{\lambda<b}\bigl(\widehat u_\lambda+\lambda\bigr).
\end{equation}
This reconstruction is a pointwise identity for the convex radial profiles $t\mapsto u(tv)$; the preceding direct approximation supplies the psh assertion for its slices. Neither the endpoint continuity assertion of \cite[Lemma~6.23]{BJ22} nor \cite[Theorem~6.21]{BJ22} is used.

\begin{proposition}\label{prop:homogeneous-slices}
Let $X$ be an integral projective variety and $\theta\in N^1(X)_{\mathbf R}$. Suppose that for every increasing net $(\psi_i)$ in $\PSHhom(X,\theta)$ there is $\psi\in\PSHhom(X,\theta)$ with $\psi=\sup_i\psi_i$ on $X^{\divv}$. Then $\EP_X(\theta)$ holds.
\end{proposition}

\begin{proof}
We verify the criterion of \cref{prop:EP-criterion}. Let $(\varphi_i)$ be a bounded-above increasing net in $\PSH(X,\theta)$ and put $f\coloneqq\sup_i\varphi_i$ on $X^{\divv}$. After subtracting a constant we may assume that $b_i\coloneqq\sup\varphi_i$ increases to $0$.

\emph{Step 1: limits of slices.} Let $\lambda<0$. For $i$ with $b_i>\lambda$, the slices $\widehat\varphi_{i,\lambda}$ lie in $\PSHhom(X,\theta)$, and they increase with $i$, since $\varphi\mapsto\widehat\varphi_\lambda$ is monotone. By hypothesis there is $\psi_\lambda\in\PSHhom(X,\theta)$ with $\psi_\lambda=\sup_i\widehat\varphi_{i,\lambda}$ on $X^{\divv}$. The functions $\psi_\lambda$ decrease in $\lambda$, because the slices do. By \eqref{eq:Legendre-reconstruction}, after interchanging two suprema,
\[
 f(v)=\sup_{\lambda<0}\{\psi_\lambda(v)+\lambda\}\qquad(v\in X^{\divv}).
\]
Finally, the hypothesis applied to the increasing sequence $(\psi_{-n})_{n\ge1}$ gives $\Psi\in\PSHhom(X,\theta)$ with $\Psi\ge\psi_\lambda$ for all $\lambda<0$.

\emph{Step 2: a decreasing approximation of $f$.} For $n\ge1$, let $A_n\coloneqq\{-n,-n+2^{-n},\ldots,-2^{-n}\}$ be the set of left endpoints of the partition of $[-n,0]$ into intervals of length $2^{-n}$, and put
\[
 U_n\coloneqq\max\Bigl\{\Psi-n,\ \max_{a\in A_n}(\psi_a+a+2^{-n})\Bigr\}.
\]
Each $U_n$ is $\theta$-psh and vanishes at $\vtriv$. Passing from $n$ to $n+1$, each old interval is subdivided, which does not increase its contribution because $\psi_\lambda$ decreases in $\lambda$, and the new intervals in $[-n-1,-n]$ and the new term $\Psi-n-1$ contribute at most $\Psi-n$. Hence $U_{n+1}\le U_n$. All its homogeneous constituents are bounded above by $\Psi$, and the term $\Psi-n$ is included, so \cref{lem:decreasing-finite-maxima} gives $U\coloneqq\lim_nU_n\in\PSH(X,\theta)$; the limit is zero at $\vtriv$. For $v\in X^{\divv}$, the values $f(v)$ and $\Psi(v)$ are finite, and
\[
 f(v)\leq U_n(v)\leq\max\{\Psi(v)-n,\ f(v)+2^{-n}\}.
\]
The lower bound holds because the parameters $\lambda<-n$ contribute at most $\Psi-n$ and each $\lambda\in[a,a+2^{-n}]$ contributes at most $\psi_a+a+2^{-n}$; the upper bound holds because $\psi_a+a\le f$ on $X^{\divv}$. Hence $U=f$ on $X^{\divv}$, which is the criterion of \cref{prop:EP-criterion}.
\end{proof}

\begin{corollary}\label{cor:normal-EP}
Let $X$ be a normal projective variety over the algebraically closed trivially valued field $k$. Then $\EP_X(\theta)$ holds for every $\theta\in N^1(X)_{\mathbf R}$.
\end{corollary}

\begin{proof}
We verify the hypothesis of \cref{prop:homogeneous-slices}. Let $(\psi_i)$ be an increasing net in $\PSHhom(X,\theta)$. Choose an effective Cartier divisor $E$ such that $\eta\coloneqq\theta+[E]$ is ample. The functions $\psi_i+e_E\in\PSHhom(X,\eta)$ increase and are bounded above by $0$. By \cref{thm:normal-ample-EP,prop:EP-criterion}, there is $\chi\in\PSHhom(X,\eta)$ with $\chi=\sup_i(\psi_i+e_E)$ on $X^{\divv}$. Thus $\chi\le e_E$ on $X^{\divv}$, hence everywhere. By \cref{thm:homogeneous-Siu}, $\chi=e_E+\psi$ with $\psi\in\PSHhom(X,\theta)$, and $\psi=\sup_i\psi_i$ on $X^{\divv}$.
\end{proof}

\section{Descent through the normalization}\label{sec:normalization}

Let $X$ be an integral unibranch projective variety and $\nu\colon Y\coloneqq X^\nu\to X$ its normalization. We prove that psh functions descend through $\nu$ (\cref{thm:normalization-reflection}) and deduce the main theorem. We first descend homogeneous Fubini--Study functions and then homogeneous psh functions. Reconstruction from homogeneous slices gives the general case: in positive characteristic a Frobenius power of every section upstairs descends, and in characteristic zero the additive powers of Roberts--Singh descend after a twist independent of the power.

\begin{lemma}\label{lem:unibranch-universal}
The normalization $\nu$ is a finite birational universal homeomorphism, and $\nu^{\an}\colon Y^{\an}\to X^{\an}$ is a homeomorphism.
\end{lemma}

\begin{proof}
The normalization is finite, varieties being Nagata \cite[\href{https://stacks.math.columbia.edu/tag/035B}{Tag 035B} and \href{https://stacks.math.columbia.edu/tag/035S}{Tag 035S}]{Stacks}, and birational. It is bijective, since the points of $Y$ over $x\in X$ correspond to the branches of $\mathcal O_{X,x}$ \cite[\href{https://stacks.math.columbia.edu/tag/0C37}{Tag 0C37} and \href{https://stacks.math.columbia.edu/tag/0C3B}{Tag 0C3B}]{Stacks}, and $X$ is unibranch \cite[\href{https://stacks.math.columbia.edu/tag/0BQ2}{Tag 0BQ2}]{Stacks}. The diagonal $Y\to Y\times_XY$ is a closed immersion, $\nu$ being finite; were it not surjective, its open complement would be a nonempty scheme of finite type over $k$, hence would have a $k$-point, that is, two distinct $k$-points of $Y$ over one point of $X$, contradicting bijectivity. Thus $\nu$ is radicial \cite[\href{https://stacks.math.columbia.edu/tag/01S4}{Tag 01S4}]{Stacks}, hence a universal homeomorphism \cite[\href{https://stacks.math.columbia.edu/tag/04DF}{Tag 04DF}]{Stacks}. In particular $\nu$ is bijective with purely inseparable residue field extensions. A real valuation extends uniquely through a purely inseparable extension, so $\nu^{\an}$ is a continuous bijection of compact Hausdorff spaces, hence a homeomorphism.
\end{proof}

We henceforth identify $X^{\an}$ and $Y^{\an}$ through $\nu^{\an}$. An extension of rings $A\subseteq B$ is \emph{subintegral} if it is integral, $\operatorname{Spec}B\to\operatorname{Spec}A$ is bijective, and the induced residue-field maps are isomorphisms. Equivalently, every finite subset of $B$ lies in a finite tower $A\subseteq A[e_1]\subseteq\cdots\subseteq A[e_1,\ldots,e_r]\subseteq B$ with $e_i^2,e_i^3\in A[e_1,\ldots,e_{i-1}]$ \cite[\href{https://stacks.math.columbia.edu/tag/0CND}{Tag 0CND}]{Stacks}. In characteristic zero we use the additive-power criterion of Roberts--Singh in the following form.

\begin{lemma}\label{lem:normalization-local-sandwich}
Let $A\subseteq B$ be a subintegral extension of $\mathbf Q$-algebras and let $b\in B$. There are integers $N\geq d\geq 0$ and elements $c_1,\ldots,c_d\in B$ such that
\[
 a_m\coloneqq b^m+\sum_{i=1}^{\min\{d,m\}}\binom mi c_i b^{m-i}\qquad(m\geq0)
 \quad\text{satisfies}\quad a_m\in A\quad(m\geq N).
\]
Moreover, whenever $N\geq d\geq0$ and $c_1,\ldots,c_d\in B$ satisfy this, $b^{n+2d}\in(a_{n+d},\ldots,a_{n+2d})B$ for every integer $n\geq0$. This ideal-membership assertion is in $B$ and does not require the displayed $a_m$ to belong to $A$.
\end{lemma}

\begin{proof}
The first assertion follows from \cite[Theorem~4.17, Proposition~4.2 and Corollary~4.3]{RobertsSingh}, after enlarging $N$. For the second, fix $m\geq d$, put $c_0\coloneqq1$ and $P(T)\coloneqq\sum_{i=0}^{d}\binom Ti c_i b^{m+d-i}\in B[T]$. Then $P$ has degree at most $d$, $P(m+r)=b^{d-r}a_{m+r}$ for $0\leq r\leq d$, and $P(0)=b^{m+d}$. Lagrange interpolation over $\mathbf Q$ expresses $P(0)$ as a linear combination of $P(m),\ldots,P(m+d)$, so $b^{m+d}\in(a_m,\ldots,a_{m+d})B$; take $m=n+d$.
\end{proof}

\begin{theorem}\label{thm:normalization-reflection}
For every $\theta\in N^1(X)_{\mathbf{R}}$, under the identification of $X^{\an}$ with $Y^{\an}$, one has $\PSH(Y,\nu^*\theta)\subseteq\PSH(X,\theta)$.
\end{theorem}

\begin{proof}
First let $\psi\in\PSHhom(Y,\nu^*\theta)$. By homogeneous regularization and \cref{lem:finite-class-promotion}, write $\psi_i\downarrow\psi$ with $\psi_i\in\mathcal H_{\mathrm{hom}}(Y,\nu^*L_i)$ and $c_1(L_i)\to\theta$, presented by sections of $\nu^*(mL_i)$ with $mL_i$ a line bundle on $X$. Class promotion preserves homogeneity. By \cref{lem:homogeneous-decreasing-stability}, it suffices to descend each homogeneous Fubini--Study function to a psh function in its original class. We prove this assertion in the two characteristics.

\emph{Positive characteristic $p=\operatorname{char}k$.} Write $\psi=\frac1m\max_j\{\log|s_j|+\lambda_j\}$ with $s_j\in H^0(Y,\nu^*(mL))$, and fix a finite affine open cover of $X$ trivializing $mL$. For $U$ in the cover, $B_U\coloneqq\mathcal O_Y(\nu^{-1}U)$ is a finitely generated algebra over $A_U\coloneqq\mathcal O_X(U)$. Since $\nu$ is a universal homeomorphism (\cref{lem:unibranch-universal}), each generator has a $p$-power in $A_U$ by \cite[\href{https://stacks.math.columbia.edu/tag/0CNF}{Tag 0CNF}]{Stacks}; since Frobenius is a ring map, a single $p$-power $q$ satisfies $b^q\in A_U$ for every $b\in B_U$ and every $U$. Let $f_{U,j}\in B_U$ be the coefficient of $s_j$ in the pulled-back frame of $mL$ over $U$, so that $f_{U,j}^{\,q}\in A_U$. On overlaps the coefficients change by the inverse transition function of $mL$, so their $q$-th powers change by that of $qmL$; as $\mathcal O_X\to\nu_*\mathcal O_Y$ is injective, the $f_{U,j}^{\,q}$ glue to $s_j^{[q]}\in H^0(X,qmL)$ with $\nu^*s_j^{[q]}=s_j^{\otimes q}$. Hence $\psi$ coincides with $\frac1{qm}\max_j\{\log|s_j^{[q]}|+q\lambda_j\}\in\mathcal{H}^{\mathrm{gf}}_{\mathbf{Q}}(X,L)$, which is contained in $\PSH(X,c_1(L))$ by \cref{prop:BJ-package}(3).

\emph{Characteristic zero.} By \cref{prop:BJ-package}(5), it suffices to show that for every line bundle $M$ on $X$ and every nonzero $s\in H^0(Y,\nu^*M)$, the function $u\coloneqq\log|s|\leq0$ lies in $\PSH(X,c_1(M))$.

Choose a finite cover of $X$ by affine opens $U_\alpha=D(\tau_\alpha)$ trivializing $M$, with all $\tau_\alpha$ sections of one ample line bundle $A_0$, fix frames $e_\alpha$ of $M$ over $U_\alpha$, and write $s=f_\alpha\nu^*e_\alpha$ with $f_\alpha\in B_\alpha\coloneqq\mathcal O_Y(\nu^{-1}U_\alpha)$. Since $\operatorname{char}k=0$, $\nu$ induces isomorphisms on residue fields by \cref{lem:unibranch-universal}, so $A_\alpha\coloneqq\mathcal O_X(U_\alpha)\subseteq B_\alpha$ is subintegral \cite[\href{https://stacks.math.columbia.edu/tag/0CND}{Tag 0CND}]{Stacks}. Apply \cref{lem:normalization-local-sandwich} to $f_\alpha$, obtaining $c_{\alpha i}$ and $a_{\alpha m}$; enlarging $d$ (with new coefficients $0$, which changes no $a_{\alpha m}$) and $N$, we take $d$ and $N$ independent of $\alpha$. Put $\mathcal Q\coloneqq\nu_*\mathcal O_Y/\mathcal O_X$, and fix an integer $R$ with the following properties, each of which holds for all large $R$: each $\tau_\alpha^Rc_{\alpha i}(\nu^*e_\alpha)^{\otimes i}$ extends to a global section $\widetilde c_{\alpha i}$ of $\nu^*(iM+RA_0)$; multiplication by $\tau_\alpha^R$ kills every section of a twist of $\mathcal Q$ by a line bundle that is supported on $V(\tau_\alpha)$ (locally, by stabilization of the kernels of multiplication by powers of $\tau_\alpha$ on the coherent sheaf $\mathcal Q$); and $RA_0+jM$ is globally generated for $0\leq j\leq d$.

For $m\geq N$, the section
\[
 \widetilde a_{\alpha m}\coloneqq\tau_\alpha^Rs^{\otimes m}+\sum_{i=1}^{d}\binom mi\widetilde c_{\alpha i}s^{\otimes(m-i)}\in H^0\bigl(Y,\nu^*(mM+RA_0)\bigr)
\]
equals $\tau_\alpha^Ra_{\alpha m}(\nu^*e_\alpha)^{\otimes m}$ over $\nu^{-1}U_\alpha$, so its image in $H^0(X,\mathcal Q\otimes(mM+RA_0))$ is supported on $V(\tau_\alpha)$ and killed by $\tau_\alpha^R$. Tensoring $0\to\mathcal O_X\to\nu_*\mathcal O_Y\to\mathcal Q\to0$ by $mM+2RA_0$ and using left exactness of global sections, we obtain a unique $t_{\alpha m}\in H^0(X,mM+2RA_0)$ with $\nu^*t_{\alpha m}=\tau_\alpha^R\widetilde a_{\alpha m}$.

Fix $n\geq\max\{1,N-d\}$. Multiplying each $t_{\alpha m}$ with $n+d\leq m\leq n+2d$ by each member of a fixed generating system of $RA_0+(n+2d-m)M$ gives a finite system $\mathcal S_n$ of sections of $(n+2d)M+3RA_0$. For $v\in X^{\an}$, evaluating in local frames at the center, we claim
\begin{equation}\label{eq:normalization-global-sandwich}
 n\,v(s)\leq\min_{t\in\mathcal S_n}v(t)\leq(n+2d)\,v(s).
\end{equation}
The first inequality holds because $m-i\geq n$ in the definition of $\widetilde a_{\alpha m}$, so that each $\nu^*t$ with $t\in\mathcal S_n$ is $s^{\otimes n}$ times a global section of a line bundle on $Y$. For the second, choose $\alpha$ with $c(v)\in U_\alpha$. Then $\tau_\alpha$ and a member of each generating system are units at $c(v)$, so $\min_tv(t)\leq v(a_{\alpha m})$ for $n+d\leq m\leq n+2d$, and since $v\geq0$ on $B_\alpha$, \cref{lem:normalization-local-sandwich} gives $\min_mv(a_{\alpha m})\leq(n+2d)v(f_\alpha)=(n+2d)v(s)$. At $\vtriv$ both bounds vanish, so $\mathcal S_n$ contains a nonzero section, and the function
\[
 F_n\coloneqq\frac1{n+2d}\max_{t\in\mathcal S_n}\log|t|\in\mathcal{H}_{\mathrm{hom}}\Bigl(X,M+\frac{3R}{n+2d}A_0\Bigr)
\]
satisfies $u\leq F_n\leq\frac n{n+2d}\,u$ by \eqref{eq:normalization-global-sandwich}.

Choose an integer $n_0>2d$ with $n_0\geq N-d$, and put $n_j\coloneqq2^jn_0$ and $G_j\coloneqq(1-2d/n_j)F_{n_j}$. For $n=n_j$,
\[
 G_{j+1}\leq\frac{n-d}{n+d}\,u\leq\frac{n-2d}{n}\,u\leq G_j,
\]
since $u\leq0$ and $\frac{n-d}{n+d}-\frac{n-2d}n=\frac{2d^2}{n(n+d)}\geq0$. Thus $G_j$ decreases to $u$ pointwise, also where $u=-\infty$, its classes $(1-2d/n_j)\bigl(c_1(M)+\frac{3R}{n_j+2d}c_1(A_0)\bigr)$ tend to $c_1(M)$, and $G_j(\vtriv)=0$. The $G_j$ are homogeneous Fubini--Study functions, so \cite[Definition~4.1]{BJ22} gives $u\in\PSH(X,c_1(M))$ directly.

This proves reflection for homogeneous functions. Now let $\psi\in\PSH(Y,\nu^*\theta)$ be arbitrary and subtract its supremum, so $\sup\psi=0$. Its slices $\widehat\psi_\lambda$, $\lambda<0$, are homogeneous psh. By \cref{cor:normal-EP,prop:EP-criterion}, the increasing sequence $\widehat\psi_{-n}$ has a homogeneous psh upper envelope $\widetilde\Psi$ in the class $\nu^*\theta$, dominating every slice. Homogeneous reflection gives functions $q_\lambda,\Psi\in\PSHhom(X,\theta)$ corresponding to $\widehat\psi_\lambda,\widetilde\Psi$. The homeomorphism $\nu^{\an}$ commutes with scaling, so $q_\lambda\le\Psi$ and the Legendre reconstruction identity is unchanged. Form $U_n$ from $q_\lambda$ and $\Psi$ using the dyadic formula in Step~2 of \cref{prop:homogeneous-slices}. These functions decrease and satisfy, pointwise under the identification of the analytifications,
\[
 \psi\le U_n\le\max\{\Psi-n,\psi+2^{-n}\}.
\]
Since $\Psi\le0$, the right-hand side tends to $\psi$. Thus $U_n\downarrow\psi$, and \cref{lem:decreasing-finite-maxima} proves $\psi\in\PSH(X,\theta)$. Adding back its supremum finishes the proof.
\end{proof}

\begin{proof}[Proof of \cref{thm:intro-main}]
By \cref{lem:unibranch-universal}, $\nu^{\an}$ is a homeomorphism, so pullback commutes with pointwise suprema and usc regularization, and $\nu^*\sups_i\varphi_i=\sups_i\nu^*\varphi_i$. The functions $\nu^*\varphi_i$ lie in $\PSH(Y,\nu^*\theta)$ by \cref{prop:BJ-package}(4). Since $Y$ is normal, the right-hand side lies in $\PSH(Y,\nu^*\theta)$ by \cref{cor:normal-EP}, and hence $\sups_i\varphi_i\in\PSH(X,\theta)$ by \cref{thm:normalization-reflection}.
\end{proof}

\begin{corollary}\label{cor:unibranch-decreasing-stability}
Let $X$ be an integral unibranch projective variety. If $\theta_i\to\theta$ and a net $u_i\in\PSH(X,\theta_i)$ decreases pointwise to $u\not\equiv-\infty$, then $u\in\PSH(X,\theta)$.
\end{corollary}

\begin{proof}
For every $v\in X^{\an}$, the inequalities $u_i(v)\le u_i(\vtriv)$ pass to the limit. Since $u\not\equiv-\infty$, this gives $b\coloneqq u(\vtriv)=\sup u\in\mathbf R$. Subtract $b$, so that $\sup u=0$ and $\sup u_i\ge0$, and restrict to a tail we may assume that $\sup_X u_i$ is bounded from above. For each $\lambda<0$, the homogeneous slices $\widehat u_{i,\lambda}$ decrease and have class $\theta_i$. By \cref{lem:homogeneous-decreasing-stability} their limit $q_\lambda$ belongs to $\PSHhom(X,\theta)$. Commuting the two infima gives
\[
 q_\lambda=\inf_i\inf_{t>0}\{t\mathbin\cdot u_i-t\lambda\}
 =\inf_{t>0}\{t\mathbin\cdot u-t\lambda\}.
\]
The functions $t\mapsto u(tv)$ are decreasing limits of convex nonincreasing functions with a common finite upper bound, and the limit satisfies $u\le0$. So it is either finite and convex or identically $-\infty$. Legendre inversion gives $u=\sup_{\lambda<0}(q_\lambda+\lambda)$.

By \cref{thm:intro-main,prop:EP-criterion}, the increasing sequence $q_{-n}$ has an upper envelope $\Psi\in\PSHhom(X,\theta)$, with $q_\lambda\le\Psi$ for every $\lambda<0$. For $n\ge1$, set $A_n\coloneqq\{-n,-n+2^{-n},\ldots,-2^{-n}\}$, and define
\[
 U_n\coloneqq\max\Bigl\{\Psi-n,\ \max_{a\in A_n}(q_a+a+2^{-n})\Bigr\}.
\]
As in Step~2 of \cref{prop:homogeneous-slices}, refining the partition cannot increase the contributions from the old intervals, since $q_\lambda$ decreases in $\lambda$; the new intervals in $[-n-1,-n]$ contribute at most $\Psi-n$. Thus $U_{n+1}\le U_n$. The reconstruction identity above gives
\[
 u\le U_n\le\max\{\Psi-n,\ u+2^{-n}\}.
\]
Since $\Psi\le0$ and $2^{-n}\to0$, these bounds show that $U_n\downarrow u$ everywhere. Each $U_n$ is a finite maximum of translates of homogeneous $\theta$-psh functions bounded above by $\Psi$, and includes the term $\Psi-n$. Hence \cref{lem:decreasing-finite-maxima} applies and gives $u\in\PSH(X,\theta)$, since $u\not\equiv-\infty$.
\end{proof}

\appendix
\section{The singular Yau--Tian--Donaldson conjecture}
\label{sec:appendix}

\subsection{Introduction}\label{subsec:appendix_introduction}

In this appendix, we apply the envelope property to the weighted uniform Yau--Tian--Donaldson conjecture for projective varieties with klt singularities. We show that, on a klt projective variety, weighted $\mathbb{G}$-uniform K-stability is equivalent to weighted $\mathbb{G}$-coercivity of the Mabuchi functional. The proof follows the strategy of \cite[Appendix~A]{T26App}, with a slight generalization to this setting. The envelope property allows us to reduce the problem on a singular variety to one on a smooth resolution.

Let $X$ be a projective variety over $\mathbb{C}$ with at most klt singularities, and let $L_X$ be an ample $\mathbb{Q}$-line bundle on $X$. Let $T$ be an $r$-dimensional ($r=0$ is allowed) real torus generated by Hamiltonian vector fields on $X$. Let $P$ be the associated convex polytope (when $r=0$, $P$ is a single point), and let $\Aut_T(X)$ denote the connected subgroup of $\Aut(X)$ commuting with $T$. Let $K \subset \Aut_T(X)$ be a compact subgroup containing a maximal torus of $\Aut_T(X)$. Set $\mathbb{G}=K_\mathbb{C}$, and $\mathbb{T}$ as the center of $\mathbb{G}$. Fix $\mathrm{v}\in C^\infty(P,\mathbb{R}_{>0})$ and $\mathrm{w}\in C^\infty(P,\mathbb{R})$. We refer to \cite[Section 1]{HL25} for these definitions. Let $\pi\colon Y\to X$ be a $\mathbb{G}$-equivariant resolution of $X$, set $L_Y\coloneqq\pi^*L_X$ and $n\coloneqq\dim X$.

\begin{theorem}\label{thm:weighted_uniform_YTD}
    With the notation above, the following statements are equivalent:
    \begin{enumerate}
        \item $(X,L_X)$ is $(\mathrm{v},\mathrm{w})$-weighted $\mathbb{G}$-uniformly K-stable.
        \item The weighted Mabuchi functional $\mathbf{M}_{\mathrm{v},\mathrm{w}}$ is $\mathbb{G}$-coercive.
    \end{enumerate}
    If, in addition, $\mathrm{v}$ is log-concave, $\mathrm{w}=\mathrm{w}_0\cdot\mathfrak{l}$, where $\mathrm{w}_0\in C^\infty(P,\mathbb{R}_{>0})$, $\mathfrak{l}$ is an affine function on $P$, and the singularity on $X$ is of Fano type, either statement implies the existence of a $(\mathrm{v},\mathrm{w})$-weighted cscK metric in the class $c_1(L_X)$.
\end{theorem}

\subsection{Preliminaries}\label{subsec:appendix-preliminaries}

Since $X$ is klt, it is normal, and $Y$ is smooth. The envelope property for $(X,L_X)$ and for $(Y,L_Y)$, given by \cref{cor:normal-EP}, together with \cite[Proposition~2.10, Lemma~3.4]{HL25}, gives the following proposition.
\begin{proposition}\label{prop:psh_Y_X}
    With the notation of \cref{subsec:appendix_introduction}, let $(X^\an,L_X^\an)$ and $(Y^\an,L_Y^\an)$ be the Berkovich analytifications of $(X,L_X)$ and $(Y,L_Y)$, respectively. Then
    \begin{enumerate}
        \item Pullback by $\pi$ is a bijection from the $K$-invariant functions in $\PSH(X^\an,L_X^\an)$ onto the $K$-invariant functions in $\PSH(Y^\an,L_Y^\an)$, and any $K$-invariant non-Archimedean Fubini--Study metric on $(Y^\an,L_Y^\an)$ is the pullback of a $K$-invariant non-Archimedean Fubini--Study metric on $(X^\an,L_X^\an)$. 
        \item For any normal model $(\mathcal{Y},\mathcal{L}_Y)$ of $(Y,L_Y)$, there exists a normal model $(\mathcal{X},\mathcal{L}_X)$ of $(X,L_X)$ such that $P_{L_Y}(\phi_{\mathcal{L}_Y})=\pi^*P_{L_X}(\phi_{\mathcal{L}_X})$.
        \item For any big normal model $(\mathcal{Y},\mathcal{L}_Y)$ of $(Y,L_Y)$, let $\mathfrak{a}_m$ be the base ideal of $|m\mathcal{L}_Y|$. Then $\phi_m = \phi_{\mathcal{L}_Y}+\frac{1}{m}\varphi_{\mathfrak{a}_m}$ converges uniformly to $P_{L_Y}(\phi_{\mathcal{L}_Y})$.
        \item For any $\mathbb{G}$-equivariant normal model $(\mathcal{X},\mathcal{L}_X)$,  we have
            \[
                \pi_*\MA_{\mathrm{v}}\biggl(\pi^*P_{L_X}(\phi_{\mathcal{L}_X})\biggr) = \MA_{\mathrm{v}}\biggl(P_{L_X}(\phi_{\mathcal{L}_X})\biggr),
            \]
            where the weighted non-Archimedean Monge--Ampère operator $\MA_{\mathrm{v}}$ is defined in \cite[Definition~3.1]{HL25}.
    \end{enumerate}
\end{proposition}
\begin{remark}\label{rem:big-model-stability}
Passing to a higher model, or replacing $\mathcal{L}_Y$ by $\mathcal{L}_Y+c\,\mathcal{Y}_0$ (which shifts $P_{L_Y}(\phi_{\mathcal{L}_Y})$ by the constant $c$), changes neither $\mathbf{M}^{\NA}_{\mathrm{v},\mathrm{w}}$ nor $\mathbf{J}^{\NA}_{\mathrm{v},\mathbb{T}}$. Hence uniform K-stability for models may be tested on $\mathbb{G}$-equivariant big SNC models.
\end{remark}

We also use the following positive-intersection formula for the non-Archimedean entropy from \cite[Proposition~4.1]{HL25}, which is based on \cite{LiChi21,Li22cscK}.

\begin{proposition}\label{prop:positive-intersection-na-functionals}
Let $(\mathcal{Y},\mathcal{L}_Y)$ be a $\mathbb{G}$-equivariant big SNC model of $(Y,L_Y)$ and $\phi\coloneqq P_{L_Y}(\phi_{\mathcal{L}_Y})$. Then
\[
\mathbf{H}_\mathrm{v}^{\NA}(\phi)=\langle{\mathcal{L}_Y}^n\rangle_\mathrm{v}\cdot K^{\log}_{\mathcal{Y}/Y_{\mathbb{P}^1}},
\]
where the weighted positive intersection $\langle\mathcal{L}_Y^n\rangle_{\mathrm{v}}\cdot D$ with a divisor $D$ on $\mathcal{Y}$ is defined in \cite[Proposition~4.1]{HL25}.
\end{proposition}

By \cite[Proposition~4.3, Corollary~4.4]{HL25}, weighted $\mathbb{G}$-uniform K-stability lifts from $(X,L_X)$ to $(Y,L_Y)$ as follows.
\begin{proposition}
    \label{prop:lift_of_K_stability}
    With the notation of \cref{subsec:appendix_introduction}, $(X,L_X)$ is $(\mathrm{v},\mathrm{w})$-weighted $\mathbb{G}$-uniformly K-stable for models if and only if $(Y,L_Y)$ is $(\mathrm{v},\mathrm{w})$-weighted $\mathbb{G}$-uniformly K-stable for models.
\end{proposition}

\begin{proof}[Proof of \cref{thm:weighted_uniform_YTD}]
Assume (1), and let $(\mathcal{Y},\mathcal{L}_Y)$ be a $\mathbb{G}$-equivariant big SNC model of $(Y,L_Y)$. For sufficiently divisible $m$, let $\mathfrak{a}_m$ be the base ideal of $|m\mathcal{L}_Y|$ and $\phi_m\coloneqq\phi_{\mathcal{L}_Y}+\frac1m\varphi_{\mathfrak{a}_m}$. By \cref{prop:psh_Y_X}(1), $\phi_m$ is the pullback of a $K$-invariant Fubini--Study metric on $(X^{\an},L_X^{\an})$, attached to a $\mathbb{G}$-equivariant test configuration of $(X,L_X)$ with the same values of $\mathbf{M}^{\NA}_{\mathrm{v},\mathrm{w}}$ and $\mathbf{J}^{\NA}_{\mathrm{v},\mathbb{T}}$ \cite[Proposition~4.3]{HL25}. Hence (1) gives $\gamma>0$, independent of $m$ and of the model, such that
\begin{equation}\label{eq:K_stable_formula}
 \mathbf{M}^{\NA}_{\mathrm{v},\mathrm{w}}(\phi_m)\geq\gamma\,\mathbf{J}^{\NA}_{\mathrm{v},\mathbb{T}}(\phi_m).
\end{equation}
By \cref{prop:psh_Y_X}(3), $\phi_m$ converges uniformly, hence strongly, to $\phi=P_{L_Y}(\phi_{\mathcal{L}_Y})$. Thus $\mathbf{R}^{\NA}_{\mathrm{v}}$ and $\mathbf{E}^{\NA}_{\mathrm{w}}$ converge along $(\phi_m)$, $\mathbf{J}^{\NA}_{\mathrm{v},\mathbb{T}}$ converges by \cite[Lemma~3.10]{HL25}, and $\mathbf{H}^{\NA}_{\mathrm{v}}$ converges by \cref{prop:entropy_approximation}. Letting $m\to\infty$ in \eqref{eq:K_stable_formula} gives $\mathbf{M}^{\NA}_{\mathrm{v},\mathrm{w}}(\phi)\geq\gamma\,\mathbf{J}^{\NA}_{\mathrm{v},\mathbb{T}}(\phi)$. By \cref{rem:big-model-stability}, $(Y,L_Y)$ is $(\mathrm{v},\mathrm{w})$-weighted $\mathbb{G}$-uniformly K-stable for models, hence so is $(X,L_X)$ by \cref{prop:lift_of_K_stability}, and \cite[Theorem~1.2]{HL25} gives (2).

Conversely, assume (2): there are $\gamma,C>0$ such that $\mathbf{M}_{\mathrm{v},\mathrm{w}}(\Phi_t)\geq\gamma\,\mathbf{J}_{\mathrm{v},\mathbb{T}}(\Phi_t)-C$ for every $K$-invariant geodesic ray $\Phi$ in $\mathcal{E}^1_K(X,L_X)$. For every $\mathbb{G}$-equivariant test configuration of $(X,L_X)$, with Fubini--Study metric $\phi$, the slope formulas \cite[Lemma~3.9, Corollary~3.31]{HL25} along the associated $K$-invariant geodesic ray give $\mathbf{M}^{\NA}_{\mathrm{v},\mathrm{w}}(\phi)\geq\gamma\,\mathbf{J}^{\NA}_{\mathrm{v},\mathbb{T}}(\phi)$, which is (1).

Finally, under the additional hypotheses, (2) and \cite[Theorem~A]{PT26} give a singular $(\mathrm{v},\mathrm{w})$-weighted cscK metric in $c_1(L_X)$; the proof there combines the uniform coercivity estimate of \cite[Theorem~A]{BJT26} on the resolution with uniform a priori estimates.
\end{proof}

\subsection{Entropy approximation}
\label{subsec:entropy-approximation}

Let $(\mathcal{Y},\mathcal{L}_Y)$ be a $\mathbb{G}$-equivariant big SNC model of $(Y,L_Y)$ over $\mathbb{P}^1$, and let $\mathfrak{a}_m$ be the base ideal of $|m\mathcal{L}_Y|$. Set $\phi_m = \phi_{\mathcal{L}_Y} + \frac{1}{m}\varphi_{\mathfrak{a}_m}$ and $\phi = P_{L_Y}(\phi_{\mathcal{L}_Y})$.

\begin{proposition}
    \label{prop:entropy_approximation}
    With the notation above, we have
    \begin{equation}\label{eq:prop:entropy_approximation_formula}
        \lim_{m\to\infty}\mathbf{H}^{\NA}_{\mathrm{v}}(\phi_m) = \mathbf{H}^{\NA}_{\mathrm{v}}(\phi).
    \end{equation}
\end{proposition}

\begin{proof}
There exist a $\mathbb{G}$-equivariant projective birational morphism $\rho\colon\mathcal{Y}'\to\mathcal{Y}$ with $\mathcal Y'$ smooth, an ample $\mathbb{Q}$-divisor $A$ and an effective $\mathbb{Q}$-divisor $G$ on $\mathcal{Y}'$, and an effective $\mathbb{Q}$-divisor $D$ in the class of $\rho^*\mathcal{L}_Y$ such that
    \[
        D\sim_{\mathbb{Q}}A+G,
    \]
    and $\Supp G = \mathbf{B}_+(\rho^*\mathcal{L}_Y)$. Unlike in \cite[Appendix~A]{T26App}, $G$ may have horizontal components.

Let $m$ be a sufficiently divisible positive integer. Set $\mathfrak{a}'_m = \mathfrak{a}_m\cdot\mathcal{O}_{\mathcal{Y}'}$, and let $p_m\colon \mathcal{Y}_m\rightarrow \mathcal{Y}'$ be the normalized blowup of $\mathfrak{a}'_m$. Let $r_m\colon \tilde{\mathcal{Y}}_m\rightarrow \mathcal{Y'}$ be a $\mathbb{G}$-equivariant log resolution of $\mathfrak{a}'_m$, which factors through $q_m\colon \tilde{\mathcal{Y}}_m\rightarrow \mathcal{Y}_m$. On $\mathcal{Y}_m$, there exist a semiample divisor $P_m$ and an effective divisor $F_m$ such that
\[
    p_m^*(mD) = P_m+F_m,
\]
where $P_m$ is the moving part, $F_m$ is the fixed part. Let $M_m = \frac{1}{m}P_m$.

By \cref{prop:positive-intersection-na-functionals},
\[
    \mathbf{H}^{\NA}_{\mathrm{v}}(\phi_m) = \left(q_m^*M_m^n\cdot K^{\log}_{\tilde{\mathcal{Y}}_m/Y_{\mathbb{{P}}^1}}\right)_{\mathrm{v}}, \qquad \mathbf{H}^{\NA}_{\mathrm{v}}(\phi) = \langle \mathcal{L}_Y^n\rangle_{\mathrm{v}}\cdot  K^{\log}_{\mathcal{Y}/Y_{\mathbb{P}^1}}.
\]
On the other hand, by \cite[Proposition~4.1]{HL25},
\[
    \lim_{m\to\infty}\left(q_m^*M_m^n\cdot K^{\log}_{{\mathcal{Y}}/Y_{\mathbb{P}^1}}\right)_{\mathrm{v}} = \langle \mathcal{L}_Y^n\rangle_{\mathrm{v}}\cdot K^{\log}_{\mathcal{Y}/Y_{\mathbb{P}^1}}.
\]
Moreover
\[
    0\leq \left(q_m^*M_m^n\cdot K^{\log}_{\tilde{\mathcal{Y}}_m/Y_{\mathbb{{P}}^1}}\right)_{\mathrm{v}}- \left(q_m^*M_m^n\cdot K^{\log}_{{\mathcal{Y}}/Y_{\mathbb{P}^1}}\right)_{\mathrm{v}} \leq \left( q_m^*M_m^n\cdot K_{\tilde{\mathcal{Y}}_m/\mathcal{Y}} \right)_{\mathrm{v}}.
\]
To prove \eqref{eq:prop:entropy_approximation_formula}, it therefore suffices to show that
\[
    \lim_{m\to\infty}\left( q_m^*M_m^n\cdot K_{\tilde{\mathcal{Y}}_m/\mathcal{Y}} \right)_{\mathrm{v}} = 0.
\]
Since $\mathrm{v}$ is positive and bounded, it suffices to show that
\[
    \lim_{m\to\infty} q_m^*M_m^n\cdot K_{\tilde{\mathcal{Y}}_m/\mathcal{Y}}  = 0.
\]
Since 
\[
    \Supp (K_{\mathcal{Y}'/\mathcal{Y}}) \subset \operatorname{Exc}(\rho) \subset \mathbf{B}_+(\rho^*\mathcal{L}_Y),
\]
we have
\[
     \lim_{m\to\infty}q_m^*M_m^n\cdot K_{\mathcal{Y}'/\mathcal{Y}} = \langle (\rho^*\mathcal{L}_Y)^n\rangle\cdot K_{\mathcal{Y}'/\mathcal{Y}} = 0.
\]

It remains to prove that
\[
    \lim_{m\to\infty}q_m^*M_m^n\cdot K_{\tilde{\mathcal{Y}}_m/\mathcal{Y}'} = 0.
\]
Let $\mathfrak{b}_m\coloneqq\mathcal{J}(\|\mathfrak{a}'^{\,m}_{\bullet}\|)$ be the asymptotic multiplier ideal on $\mathcal{Y}'$. By Nadel vanishing and Castelnuovo--Mumford regularity, there are integers $N$ and $m_0$, with $N$ dividing $m_0$ and $ND$, $m_0A$, $m_0G$ integral and $m_0D\sim m_0A+m_0G$, such that $\mathcal{O}_{\mathcal{Y}'}(mD+m_0A)\otimes\mathfrak{b}_m$ is globally generated for all $m\geq1$ divisible by $N$; below $m$ runs through such integers. Since $(m+m_0)D\sim mD+m_0A+m_0G$, this gives $\mathfrak{b}_m\cdot\mathcal{O}_{\mathcal{Y}'}(-m_0G)\subseteq\mathfrak{a}'_{m+m_0}$. For $m>m_0$, as $\mathfrak{b}_m\subseteq\mathfrak{b}_{m-m_0}$ and $\mathfrak{a}'_m\subseteq\mathcal{J}(\mathfrak{a}'_m)\subseteq\mathfrak{b}_m$, we get $\mathfrak{b}_m\cdot\mathcal{O}_{\mathcal{Y}'}(-m_0G)\subseteq\mathfrak{a}'_m\subseteq\mathfrak{b}_m$, hence
\[
 0\leq\ord_B(\mathfrak{a}'_m)-\ord_B(\mathfrak{b}_m)\leq m_0\ord_B(G)
\]
for every prime divisor $B$ over $\mathcal{Y}'$.

Let $E$ be a prime exceptional divisor of $p_m\colon \mathcal{Y}_m\rightarrow \mathcal{Y}'$. By \cite[Proposition~A.1]{T26App}, 
\[
  0\leq  \ord_E(K_{\mathcal{Y}_m/\mathcal{Y}'})\leq (n+1)\Bigl(\ord_E(\mathfrak{a}'_m)-\ord_E\bigl(\mathcal{J}(\mathfrak{a}'_m)\bigr)\Bigr),
\]
where the first inequality follows from the fact that $\mathcal{Y}'$ is smooth. Since $\mathcal{J}(\mathfrak{a}'_m)\subseteq\mathfrak{b}_m$, we have
\[
    0\leq  \ord_E(K_{\mathcal{Y}_m/\mathcal{Y}'})\leq (n+1)\Bigl(\ord_E(\mathfrak{a}'_m)-\ord_E(\mathfrak{b}_m)\Bigr) \leq (n+1)m_0\,\ord_E(G).
\]
By the projection formula,
\[
    q_m^*M_m^n\cdot K_{\tilde{\mathcal{Y}}_m/\mathcal{Y}'} = M_m^n\cdot K_{\mathcal{Y}_m/\mathcal{Y}'}.
\]
Applying the preceding inequality and \cite[Proposition~4.1]{HL25} again, we obtain
\[
    0\leq \varlimsup_{m\to\infty}q_m^*M_m^n\cdot K_{\tilde{\mathcal{Y}}_m/\mathcal{Y}'} \leq (n+1)m_0\lim_{m\to\infty}M_m^n\cdot G = (n+1)m_0\langle (\rho^*\mathcal{L}_Y)^n\rangle\cdot G = 0,
\]
where the last equality follows from $\Supp G = \mathbf{B}_+(\rho^*\mathcal{L}_Y)$. This completes the proof.
\end{proof}

\printbibliography

@article {Li22cscK,
    AUTHOR = {Li, Chi},
     TITLE = {Geodesic rays and stability in the csc{K} problem},
   JOURNAL = {Ann. Sci. \'Ec. Norm. Sup\'er. (4)},
  FJOURNAL = {Annales Scientifiques de l'\'Ecole Normale Sup\'erieure.
              Quatri\`eme S\'erie},
    VOLUME = {55},
      YEAR = {2022},
    NUMBER = {6},
     PAGES = {1529--1574},
      ISSN = {0012-9593,1873-2151},
   MRCLASS = {32Q26 (14E05 32Q15)},
  MRNUMBER = {4517682},
MRREVIEWER = {Ahmed\ Lesfari},
       DOI = {10.24033/asens.2523},
       URL = {https://doi.org/10.24033/asens.2523},
}

@misc{CLD12,
  author       = {Chambert-Loir, Antoine and Ducros, Antoine},
  title        = {Formes diff\'{e}rentielles r\'{e}elles et courants sur les espaces de {Berkovich}},
  date         = {2025-07-25},
  eprint       = {1204.6277v2},
  eprinttype   = {arxiv},
  eprintclass  = {math.AG},
  doi          = {10.48550/arXiv.1204.6277},
  url          = {https://arxiv.org/abs/1204.6277v2},
  urldate      = {2026-09-26},
  note         = {Version 2},
  label        = {CLD}
}

@article{BJ22,
  author        = {Boucksom, S\'ebastien and Jonsson, Mattias},
  title         = {Global pluripotential theory over a trivially valued field},
  journal       = {Annales de la Facult\'e des Sciences de Toulouse. Math\'ematiques},
  fjournal      = {Annales de la Facult\'e des Sciences de Toulouse. Math\'ematiques. S\'erie 6},
  volume        = {31},
  number        = {3},
  pages         = {647--836},
  year          = {2022},
  issn          = {0240-2963,2258-7519},
  mrclass       = {32P05 (14G22 14T20 31C05 32U05)},
  mrnumber      = {4452253},
  mrreviewer    = {Jackson S. Morrow},
  doi           = {10.5802/afst.1705},
  url           = {https://doi.org/10.5802/afst.1705},
  eprint        = {1801.08229},
  archivePrefix = {arXiv},
  primaryClass  = {math.AG},
}

@misc{Stacks,
  author       = {{The Stacks Project Authors}},
  title        = {Stacks Project},
  howpublished = {\url{http://stacks.math.columbia.edu}},
  year         = {2020},
  shorthand    = {Stacks}
}

@article {FN80,
    AUTHOR = {Forn{\ae}ss, John Erik and Narasimhan, Raghavan},
     TITLE = {The {L}evi problem on complex spaces with singularities},
   JOURNAL = {Math. Ann.},
  FJOURNAL = {Mathematische Annalen},
    VOLUME = {248},
      YEAR = {1980},
    NUMBER = {1},
     PAGES = {47--72},
      ISSN = {0025-5831},
   MRCLASS = {32E10},
  MRNUMBER = {569410},
MRREVIEWER = {Peter Pflug},
       DOI = {10.1007/BF01349254},
       URL = {https://doi.org/10.1007/BF01349254},
}

@article {BE21,
    AUTHOR = {Boucksom, S\'{e}bastien and Eriksson, Dennis},
     TITLE = {Spaces of norms, determinant of cohomology and {F}ekete points
              in non-{A}rchimedean geometry},
   JOURNAL = {Adv. Math.},
  FJOURNAL = {Advances in Mathematics},
    VOLUME = {378},
      YEAR = {2021},
     PAGES = {Paper No. 107501, 124},
      ISSN = {0001-8708},
   MRCLASS = {32P05 (14G22 31C15)},
  MRNUMBER = {4192993},
MRREVIEWER = {Mattias Jonsson},
       DOI = {10.1016/j.aim.2020.107501},
       URL = {https://doi.org/10.1016/j.aim.2020.107501},
}

@book{AC,
  title={Algèbre commutative},
  author={Bourbaki, N.},
  year={2007},
  publisher={Springer Science \& Business Media},
  shorthand={AC}
}

@book{AV,
  title={Abelian varieties},
  author={Mumford, D. and Ramanujam, C. P. and Manin, J. I.},
  volume={5},
  year={1974},
  publisher={Oxford university press Oxford}
}

@article{BJaddendum,
  author        = {Boucksom, S\'ebastien and Jonsson, Mattias},
  title         = {Addendum to the article ``Global pluripotential theory over a trivially valued field''},
  journal       = {Annales de la Facult\'e des Sciences de Toulouse. Math\'ematiques},
  fjournal      = {Annales de la Facult\'e des Sciences de Toulouse. Math\'ematiques. S\'erie 6},
  volume        = {33},
  number        = {2},
  pages         = {405--417},
  year          = {2024},
  issn          = {0240-2963,2258-7519},
  mrclass       = {32P05 (14G22 14T20 31C05 32U05)},
  mrnumber      = {4806788},
  doi           = {10.5802/afst.1775},
  url           = {https://doi.org/10.5802/afst.1775},
  eprint        = {2206.07183},
  archivePrefix = {arXiv},
  primaryClass  = {math.AG},
}

@incollection{MustataTest,
  author        = {Musta\c{t}\u{a}, Mircea},
  title         = {The non-nef locus in positive characteristic},
  booktitle     = {A celebration of algebraic geometry},
  series        = {Clay Mathematics Proceedings},
  volume        = {18},
  pages         = {535--551},
  publisher     = {American Mathematical Society},
  address       = {Providence, RI},
  year          = {2013},
  isbn          = {978-0-8218-8983-1},
  mrclass       = {14C20 (13A35 14F17 14F18)},
  mrnumber      = {3114955},
  mrreviewer    = {Karl Schwede},
  eprint        = {1109.3825},
  archivePrefix = {arXiv},
  primaryClass  = {math.AG},
}

@article{Canton,
  author        = {Canton, Eric},
  title         = {Berkovich log discrepancies in positive characteristic},
  journal       = {Pure and Applied Mathematics Quarterly},
  volume        = {16},
  number        = {5},
  pages         = {1465--1532},
  year          = {2020},
  issn          = {1558-8599,1558-8602},
  mrclass       = {14F18 (13A35 14B05)},
  mrnumber      = {4221003},
  doi           = {10.4310/PAMQ.2020.v16.n5.a5},
  url           = {https://doi.org/10.4310/PAMQ.2020.v16.n5.a5},
  eprint        = {1711.03002},
  archivePrefix = {arXiv},
  primaryClass  = {math.AG},
}

@article{RobertsSingh,
  author   = {Roberts, Leslie G. and Singh, Balwant},
  title    = {Subintegrality, invertible modules and the {Picard} group},
  journal  = {Compositio Mathematica},
  volume   = {85},
  number   = {3},
  pages    = {249--279},
  year     = {1993},
  mrnumber = {1214447},
}

@article{deJong,
  label = {DJ},
  shorthand = {DJ96},
  author    = {de Jong, Aise Johan},
  title     = {Smoothness, semi-stability and alterations},
  journal   = {Publications Math\'ematiques de l'{IH\'ES}},
  volume    = {83},
  pages     = {51--93},
  year      = {1996},
  doi       = {10.1007/BF02698644},
}

@book{Hartshorne,
  author    = {Hartshorne, Robin},
  title     = {Algebraic Geometry},
  series    = {Graduate Texts in Mathematics},
  volume    = {52},
  publisher = {Springer-Verlag},
  address   = {New York-Heidelberg},
  year      = {1977},
  pages     = {xvi+496},
  isbn      = {0-387-90244-9},
  mrclass   = {14-01},
  mrnumber  = {463157},
  mrreviewer = {Robert Speiser},
  doi       = {10.1007/978-1-4757-3849-0},
  url       = {https://doi.org/10.1007/978-1-4757-3849-0},
}

@book{LazarsfeldI,
  author    = {Lazarsfeld, Robert},
  title     = {Positivity in Algebraic Geometry. {I}},
  series    = {Ergebnisse der Mathematik und ihrer Grenzgebiete},
  volume    = {48},
  note      = {Classical setting: line bundles and linear series},
  publisher = {Springer-Verlag},
  address   = {Berlin},
  year      = {2004},
  pages     = {xviii+387},
  isbn      = {3-540-22533-1},
  mrclass   = {14-02 (14C20)},
  mrnumber  = {2095471},
  mrreviewer = {Mihnea Popa},
  doi       = {10.1007/978-3-642-18808-4},
  url       = {https://doi.org/10.1007/978-3-642-18808-4},
}

@book{Berkovich90,
  author    = {Berkovich, Vladimir G.},
  title     = {Spectral theory and analytic geometry over non-{A}rchimedean fields},
  series    = {Mathematical Surveys and Monographs},
  volume    = {33},
  publisher = {American Mathematical Society},
  address   = {Providence, RI},
  year      = {1990},
  doi       = {10.1090/surv/033},
}

@article{HaraYoshida,
  author  = {Hara, Nobuo and Yoshida, Ken-ichi},
  title   = {A generalization of tight closure and multiplier ideals},
  journal = {Transactions of the American Mathematical Society},
  volume  = {355},
  number  = {8},
  year    = {2003},
  pages   = {3143--3174},
  issn    = {0002-9947,1088-6850},
  mrclass = {13A35},
  mrnumber = {1974679},
  mrreviewer = {Irena Swanson},
  doi     = {10.1090/S0002-9947-03-03285-9},
  url     = {https://doi.org/10.1090/S0002-9947-03-03285-9},
}

@article{BMS08,
  author  = {Blickle, Manuel and Musta\c{t}\u{a}, Mircea and Smith, Karen E.},
  title   = {Discreteness and rationality of {$F$}-thresholds},
  journal = {Michigan Mathematical Journal},
  volume  = {57},
  year    = {2008},
  pages   = {43--61},
  note    = {Special volume in honor of Melvin Hochster},
  issn    = {0026-2285,1945-2365},
  mrclass = {13A35 (14F18)},
  mrnumber = {2492440},
  mrreviewer = {Tomasz Szemberg},
  doi     = {10.1307/mmj/1220879396},
  url     = {https://doi.org/10.1307/mmj/1220879396},
}

@incollection{SchwedeTucker,
  author    = {Schwede, Karl and Tucker, Kevin},
  title     = {A survey of test ideals},
  booktitle = {Progress in commutative algebra 2},
  pages     = {39--99},
  publisher = {Walter de Gruyter},
  address   = {Berlin},
  year      = {2012},
  isbn      = {978-3-11-027859-0},
  mrclass   = {13A35 (14B05 14F18)},
  mrnumber  = {2932591},
  mrreviewer = {Wenliang Zhang},
}

@article{Swan80,
  author  = {Swan, Richard G.},
  title   = {On seminormality},
  journal = {Journal of Algebra},
  volume  = {67},
  number  = {1},
  year    = {1980},
  pages   = {210--229},
  doi     = {10.1016/0021-8693(80)90318-X},
}

@misc{HL25,
  author  = {Han, Jiyuan and Liu, Yaxiong},
  title   = {{$\mathbb G$}-uniform weighted {K}-stability for models on klt varieties},
  year    = {2025},
  eprint  = {2506.18039},
  archiveprefix = {arXiv},
  primaryclass = {math.DG},
  note    = {Version 2, revised 17 November 2025},
}

@misc{T26App,
  author  = {Trusiani, Antonio},
  title   = {A solution to the {Yau--Tian--Donaldson} Conjecture through Special Fujita Approximations},
  year    = {2026},
  eprint  = {2605.30063},
  archiveprefix = {arXiv},
  note    = {Version 2, revised 17 August 2026; Appendix A by S\'ebastien Boucksom},
}

@article{PT26,
   title={Weighted cscK metrics on Kähler varieties},
   volume={113},
   ISSN={1469-7750},
   url={https://doi.org/10.1112/jlms.70442},
   DOI={10.1112/jlms.70442},
   number={2},
   journal={Journal of the London Mathematical Society},
   publisher={Wiley},
   author={Pan, Chung-Ming and T\^o, Tat Dat},
   year={2026},
   month=Feb }

@article {BJT26,
    AUTHOR = {Boucksom, S\'ebastien and Jonsson, Mattias and Trusiani,
              Antonio},
     TITLE = {Weighted extremal {K}\"ahler metrics on resolutions of
              singularities},
   JOURNAL = {Comm. Pure Appl. Math.},
  FJOURNAL = {Communications on Pure and Applied Mathematics},
    VOLUME = {79},
      YEAR = {2026},
    NUMBER = {4},
     PAGES = {1073--1148},
      ISSN = {0010-3640,1097-0312},
   MRCLASS = {58E11 (32U05 53C55)},
  MRNUMBER = {5030607},
       DOI = {10.1002/cpa.70026},
       URL = {https://doi.org/10.1002/cpa.70026},
}

@incollection{LiChi21,
  author    = {Li, Chi},
  title     = {{K}-stability and {Fujita} approximation},
  booktitle = {Birational Geometry, K\"ahler--Einstein Metrics and Degenerations},
  editor    = {Cheltsov, Ivan and Chen, Xiuxiong and Katzarkov, Ludmil and Park, Jihun},
  series    = {Springer Proceedings in Mathematics \& Statistics},
  volume    = {409},
  publisher = {Springer},
  address   = {Cham},
  year      = {2023},
  pages     = {545--566},
  doi       = {10.1007/978-3-031-17859-7_26},
  url       = {https://doi.org/10.1007/978-3-031-17859-7_26},
  eprint    = {2102.09457},
  eprinttype = {arxiv},
}

@article {GJKM19,
    AUTHOR = {Gubler, Walter and Jell, Philipp and K\"unnemann, Klaus and
              Martin, Florent},
     TITLE = {Continuity of plurisubharmonic envelopes in non-archimedean
              geometry and test ideals},
      NOTE = {With an appendix by Jos\'e{} Ignacio Burgos Gil and Mart\'in
              Sombra},
   JOURNAL = {Ann. Inst. Fourier (Grenoble)},
  FJOURNAL = {Universit\'e{} de Grenoble. Annales de l'Institut Fourier},
    VOLUME = {69},
      YEAR = {2019},
    NUMBER = {5},
     PAGES = {2331--2376},
      ISSN = {0373-0956,1777-5310},
   MRCLASS = {32P05 (14G22 14M25 32U05)},
  MRNUMBER = {4018262},
MRREVIEWER = {Mattias\ Jonsson},
       DOI = {10.5802/aif.3296},
       URL = {https://doi.org/10.5802/aif.3296},
}

@misc{GublerRabinoff25,
  author       = {Gubler, Walter and Rabinoff, Joseph},
  title        = {Classically psh and pluriharmonic functions on {Berkovich} spaces},
  date         = {2025-09-17},
  eprint       = {2506.13548v2},
  eprinttype   = {arxiv},
  eprintclass  = {math.AG},
  doi          = {10.48550/arXiv.2506.13548},
  url          = {https://arxiv.org/abs/2506.13548v2},
  urldate      = {2026-09-26},
  note         = {Version 2},
  label        = {GR}
}

@misc{Wu26,
  author       = {Wu, Lyuhui},
  title        = {On the local theory of non-archimedean psh functions},
  date         = {2026-09-23},
  eprint       = {2609.27610v1},
  eprinttype   = {arxiv},
  eprintclass  = {math.AG},
  url          = {https://arxiv.org/abs/2609.27610v1},
  urldate      = {2026-09-26},
  note         = {Version 1},
  label        = {Wu}
}

Yanbo Fang, \textsc{Aarhus university}\par\nopagebreak
	\textit{Email address:} \texttt{fangyanbo01@gmail.com}

Jiyuan Han, \textsc{Westlake university}\par\nopagebreak
	\textit{Email address:} \texttt{hanjiyuan@westlake.edu.cn}

Mingchen Xia, \textsc{Institute of Geometry and Physics, University of Science and Technology of China}\par\nopagebreak
	\textit{Email address:} \texttt{xiamingchen2008@gmail.com}

\end{document}